\documentstyle[12pt,righttag,amscd]{amsart}

\title{Perfect Crystals and $\boldsymbol{q}$-deformed Fock Spaces}

\author[Kashiwara \and Miwa \and Petersen \and Yung]
{Masaki Kashiwara \and
Tetsuji Miwa \and\\
Jens-Ulrik H. Petersen \and
Chong Ming Yung}

\address{\hspace{-1.4em}
  Research Institute for Mathematical
  Sciences,\newline Kyoto University, Kyoto 606-01, Japan
  \normalshape\medskip
\begin{tabbing}{\it E-mail addresses}:
    \quad\=\tt masaki@@kurims.kyoto-u.ac.jp,\quad\qquad
    \=\tt miwa@@kurims.kyoto-u.ac.jp,\\
    \>\tt petersen@@kurims.kyoto-u.ac.jp, \>\tt yung@@kurims.kyoto-u.ac.jp
  \end{tabbing}}

\date{March, 1996}

\numberwithin{equation}{subsection}
\allowdisplaybreaks[3]
\newcommand{\newsection}{\bigskip\section}
\newtheorem{lemma}{Lemma}[subsection]
\newtheorem{corollary}[lemma]{Corollary}
\newtheorem{sublemma}[lemma]{Sublemma}
\newtheorem{proposition}[lemma]{Proposition}
\newtheorem{theorem}[lemma]{Theorem}

\newcommand{\be}{\begin{equation}}
\newcommand{\en}{\end{equation}}
\newcommand{\bea}{\begin{eqnarray}}
\newcommand{\ena}{\end{eqnarray}}
\newcommand{\bean}{\begin{eqnarray*}}
\newcommand{\enan}{\end{eqnarray*}}
\newcommand{\refeq}[1]{(\ref{eqn:#1})}
\newcommand{\labeleq}[1]{\label{eqn:#1}}

\newcommand{\lb}[1]{\label{eqn:#1}}

\newcommand{\C}{{\Bbb C}}

\newcommand{\Q}{{\Bbb Q}}
\newcommand{\Z}{{\Bbb Z}}

\newcommand{\Hom}{\operatorname{Hom}\,}
\newcommand{\End}{\operatorname{End}\,}

\newcommand{\id}{\operatorname{id}}

\newcommand{\ch}{\operatorname{ch}\,}
\newcommand{\ad}{\operatorname{ad}\,}
\def\com[#1,#2]{\hbox{$[#1,#2]$}}
\newcommand{\wt}{\operatorname{wt}\,}
\newcommand{\Uq}{U_q(\frak{g})}                          

\newcommand{\Phit}{\widetilde{\Phi}}  
\newcommand{\Psit}{\widetilde{\Psi}}  

\newcommand{\Phim}[1]{\mathrel{\mathop{\kern0pt \Phi}\limits^#1}}
\newcommand{\Psim}[1]{\mathrel{\mathop{\kern0pt \Psi}\limits^#1}}
\newcommand{\Phin}[1]{\mathrel{\mathop{\kern0pt \Phit}\limits^#1}}
\newcommand{\Psin}[1]{\mathrel{\mathop{\kern0pt \Psit}\limits^#1}}

\newcommand{\G}{\Gg}
\newcommand{\UG}{U_q(\Gg)}
\newcommand{\SL}{\Gsl}
\newcommand{\UPG}{U'_q(\Gg)}
\newcommand{\VA}{V_\aff}
\newcommand{\BA}{B_\aff}
\newcommand{\deq}{\overset{\text{def}}{=}}
\newcommand{\F}{{\cal F}}

\newcommand{\hb}{\hfill\break}

\newcommand{\qbox}[1]{\hbox{\quad#1\quad}}
\newcommand{\Gsl}{{\frak{sl}}}

\newcommand{\Gg}{{\frak{g}}}
\newcommand{\Gh}{{\frak{h}}}

\newcommand{\BC}{{\Bbb C}}
\newcommand{\BQ}{{\Bbb Q}}
\newcommand{\BZ}{{\Bbb Z}}

\newcommand{\CF}{{\cal F}}
\newcommand{\CB}{{\cal B}}
\newcommand{\lam}{\lambda}
\newcommand{\Lam}{\Lambda}

\newcommand{\lan}{\langle}
\newcommand{\ran}{\rangle}

\newcommand{\gge}{\gg}
\newcommand{\rk}{{\operatorname{rank}}}
\newcommand{\sgn}{{\operatorname{sgn}}}
\newcommand{\cl}{{\operatorname{cl}}}
\newcommand{\vac}[1]{|{#1}\ran}
\newcommand{\Ker}{\operatorname{Ker}}
\newcommand{\akete}{\noalign{\smallskip}}

\newcommand{\te}{\tilde e}
\newcommand{\tf}{\tilde f}
\newcommand{\tC}{C}

\newcommand{\aff}{{\operatorname{aff}}}
\newcommand{\tB}{\tilde B}

\newcommand{\beq}{\begin{eqnarray}}
\newcommand{\beqn}{\begin{eqnarray*}}
\def\endeq{\end{eqnarray}}
\def\endeqn{\end{eqnarray*}}
\newcommand{\nn}{\notag}
\newcommand{\proof}{\noindent{\it Proof.}\ }
\renewcommand{\qed}{\hfill\qedsymbol\par\medskip}

\newcommand{\Lemma}{\begin{lemma}}
\newcommand{\enlemma}{\end{lemma}}
\newcommand{\Corollary}{\begin{corollary}}
\newcommand{\encorollary}{\end{corollary}}
\newcommand{\Sublemma}{\begin{sublemma}}
\newcommand{\ensublemma}{\end{sublemma}}
\newcommand{\Proposition}{\begin{proposition}}
\newcommand{\enproposition}{\end{proposition}}
\newcommand{\Theorem}{\begin{theorem}}
\newcommand{\entheorem}{\end{theorem}}
\newcommand{\bw}{{\mathop{\textstyle\bigwedge}\nolimits}}
\newcommand{\Remark}{\medskip\noindent{\it Remark.}\ }
\newcommand{\AZ}{\BZ[q]}
\newcommand{\KZ}{\BZ[q,q^{-1}]}
\newcommand{\tvac}[1]{{\overline{|{#1}\ran}}}
\newcommand{\tCF}{\overline{\cal F}}
\newcommand{\Fh}{\BQ((q))\otimes_K\CF}
\newcommand{\U}{U_q(\Gg)}
\newcommand{\isomo}{{\buildrel\sim\over\to}}
\newcommand{\btensor}{\overline\otimes}
\newcommand{\vo}{v^\circ}

\newcommand{\tU}{{\overline U}_q(\Gg)}

\newcommand{\g}{{\frak g}}
\newcommand{\N}{{\Bbb N}}
\newcommand{\Zplus}{\Z_{>0}}
\newcommand{\Vaff}{V_\aff}
\newcommand{\Baff}{B_\aff}
\newcommand{\Bmin}{B_{\min}}
\newcommand{\Pcl}{P_{\text{cl}}}
\newcommand{\bo}{b^\circ}
\newcommand{\fock}{{\cal F}}
\newcommand{\dvac}[1]{\langle#1|}
\renewcommand{\v}[1]{v_{#1}}
\renewcommand{\b}[1]{b_{#1}}
\newcommand{\set}[1]{\{ #1 \}}
\newcommand{\bigset}[1]{\left\{#1\right\}}
\newcommand{\generated}[1]{\langle #1 \rangle}
\newcommand{\brackets}[1]{\left[ #1 \right]} 
\newcommand{\range}[1]{\brackets{#1}} 
\newcommand{\pairing}[1]{\generated{#1}}

\newcommand{\tensor}{\otimes}
\newcommand{\mapto}{\rightarrow}

\newcommand{\dynkin}[1]{{\scriptstyle #1}}
\newcommand{\link}{-\!\!\!-\,}
\newcommand{\longlink}{\link\!\!\!\!\!\!\link}
\newcommand{\rightdoublelink}{\Longrightarrow\!}
\newcommand{\leftdoublelink}{\!\Longleftarrow}
\newcommand{\dcox}{{h^\vee}}
\newcommand{\union}{\cup}
\newcommand{\A}{{\cal A}}

\newcommand{\rightcrystal}[1]{\!\overset{#1}{\longrightarrow}\!}

\newcommand{\leftcrystal}[1]{\!\underset{#1}{\longleftarrow}\!}

\newsymbol{\rightleftarrows} 131D
\newcommand{\rightleftcrystal}[2]{\;\overset{#1}
  {\underset{#2}{\rightleftarrows}}\;}

\newcommand{\squeeze}{\!\!&\!\!}
\newcommand{\e}{\varepsilon}
\newcommand{\f}{\varphi}

\newcommand{\qbin}[2]{\fracwithdelims[][0pt]{#1}{#2}}
\newcommand{\DP}{{\text{DP}}}

\begin{document}
\pagestyle{empty}
\maketitle
\vspace{-1em}
\begin{center}
  q-alg/9603025\\
  RIMS-1070
\end{center}
\vspace{-1em}

\begin{abstract}
In~\cite{S,KMS} the semi-infinite wedge construction of level~1
$U_q(A^{(1)}_n)$ Fock spaces and their decomposition into the tensor product of
an irreducible $U_q(A^{(1)}_n)$-module and a bosonic Fock space was given.
Here a general scheme for the wedge construction of
$q$-deformed Fock spaces using
the theory of perfect crystals is presented.

Let $U_q(\g)$ be a quantum affine algebra.  Let $V$ be a finite-dimensional
$U'_q(\g)$-module with a perfect crystal base of level~$l$.  Let $V_\aff\simeq
V\tensor\C[z,z^{-1}]$ be the affinization of $V$, with crystal base
$(L_\aff,B_\aff)$.  The wedge space $V_\aff\wedge V_\aff$ is defined as the
quotient of $V_\aff\tensor V_\aff$ by the subspace generated by the action of
$U_q(\g)[z^a\tensor z^b +z^b\tensor z^a]_{a,b\in\Z}$ on $v\tensor v$ ($v$ an
extremal vector).  The wedge space $\bigwedge^r V_\aff$ ($r\in\N$) is defined
similarly.  Normally ordered wedges are defined by using the energy function
$H:B_\aff\tensor B_\aff\rightarrow\Z$.  Under certain assumptions, it is proved
that normally ordered wedges form a base of $\bigwedge^r V_\aff$.

A $q$-deformed Fock space is defined as the inductive limit of $\bigwedge^r
V_\aff$ as $r\rightarrow\infty$, taken along the semi-infinite wedge associated
to a ground state sequence.  It is proved that normally ordered wedges form a
base of the Fock space and that the Fock space has the structure of an
integrable $U_q(\g)$-module.  An action of the bosons, which commute with the
$U'_q(\g)$-action, is given on the Fock space.  It induces the decomposition of
the $q$-deformed Fock space into the tensor product of an irreducible
$U_q(\g)$-module and a bosonic Fock space.

As examples, Fock spaces for types $A^{(2)}_{2n}$, $B^{(1)}_n$,
$A^{(2)}_{2n-1}$, $D^{(1)}_n$ and $D^{(2)}_{n+1}$ at level~$1$ and $A^{(1)}_1$
at level~$k$ are constructed.  The commutation relations of the bosons in each
of these cases are calculated, using two point functions of vertex operators.
\end{abstract}

\setcounter{tocdepth}{2}
\tableofcontents

\pagestyle{plain}

\section{Introduction}
Let $\G$ be an affine Lie algebra. The construction of integrable highest
weight modules for $\G$ has been studied extensively for more than 15 years,
with applications to problems in mathematical physics like soliton
equations and conformal field theories. More recently, a further item was
added to the list of interactions between representation theory and
integrable systems:
the link between quantum affine algebras, $\UG$, and solvable lattice models
(see~\cite{JM} and references therein).

The link is twofold: (a) the R-matrices, which appear as the Boltzmann
weights of solvable lattice models, are intertwiners of level~$0$
$\UG$-modules, and (b) the irreducible integrable highest weight modules for
$\UG$ appear as the spaces of the eigenvectors of the corner transfer matrices.
This suggests a construction of integrable highest weight modules
by means of semi-infinite tensor products of level~$0$ modules.
In fact, in the crystal limit, such a construction was given for a large class
of representations known as the representations with perfect
crystals~\cite{KMN1}.

The idea of using Fock spaces of bosons or fermions goes back to
earlier works before the above link was found. In fact, the literature
is vast. Let us mention some of the works that are closely related
to the present work.  In~\cite{LW,KKLW}, bosonic Fock spaces were used
to construct some level~$1$ highest weight modules of affine Lie
algebras using the fact that the actions of the principal Heisenberg
subalgebras are irreducible. In~\cite{DJKM} the level~$1$ highest
weight modules of ${\frak{gl}}_\infty$ were constructed in the fermionic
Fock space. By the boson-fermion correspondence one has the action of
bosons on the Fock space. The action of affine Lie algebras such as
$\widehat{\frak{sl}}_n$, as subalgebras of ${\frak{gl}}_\infty$, was
then realized as the commutant of bosons of degree divisible by $n$.
Likewise, level~$1$ highest weight modules of other affine Lie
algebras $\frak{g}$ were constructed by realizing $\frak{g}$ as a
subalgebra of ${\frak{go}}_{\infty}$ (see also~\cite{JY}) or
${\frak{go}}_{2\infty}$.

Under the influence of quantum groups, several developments were
made further in this direction. A $q$-deformed construction of
the fermion Fock space was achieved in~\cite{H}. In~\cite{MM},
this was connected to the crystal base theory of Kashiwara~\cite{K1}.
These works and the developments in solvable lattice models led to the
semi-infinite construction of affine crystals mentioned above.

Very recently, in~\cite{S},
Stern gave a semi-infinite construction of the level~$1$ Fock spaces
for $\UG$ when $\G=\widehat{\SL}_n$. Subsequently, in~\cite{KMS},
the decomposition
of the Fock spaces into the level~$1$ irreducible highest weight modules and
the bosonic Fock space, was given. In the present paper, we give a similar
construction of Fock spaces and their decomposition, for various cases in
the class of representations with perfect crystals. The case in~\cite{S,KMS}
corresponds to the perfect crystal of level~$1$
for $A^{(1)}_n$. Here we treat
\begin{displaymath}
\text{level 1
$A^{(2)}_{2n}$, $B^{(1)}_n$, $A^{(2)}_{2n-1}$, $D^{(1)}_n$,
$D^{(2)}_{n+1}$}\quad
\text{and}\quad \text{level~$k$ $A^{(1)}_1$}.
\end{displaymath}

In order to handle these cases, we not only follow the basic strategy
in \cite{S,KMS}, but also develop some new machinery, where the
R-matrix and crystal bases play an important role.

In the following we recall the basic construction in~\cite{KMS}
and compare it with the newer version developed in this paper,
by taking the examples of level~$l$ $A^{(1)}_1$,
$l=1,2$.

\subsection{The kernel of $R-1$.}
Let $V$ be a finite-dimensional $\UPG$-module, and $\VA=V\otimes\C[z,z^{-1}]$
its affinization. The $r$-th $q$-wedge space is given by
\begin{displaymath}
\bigwedge^r\VA=\VA^{\otimes r}/N_r,
\end{displaymath}
where
\begin{displaymath}
N_r=\sum_{i=0}^{r-2}\VA^{\otimes i}\otimes N\otimes\VA^{\otimes(r-2-i)}
\end{displaymath}
and the space $N$ is a certain subspace of $\VA\otimes\VA$. Namely,
the $q$-wedge space is defined as a quotient of the tensor product of
$\VA$ modulo certain relations of nearest neighbour type.

For the level~$1$ $A^{(1)}_1$ case, the space $V$ is the $2$-dimensional
representation of $U'_q(\widehat{\SL}_2)$, $V=\Q \v{0}\oplus\Q \v{1}$.
In~\cite{S,KMS},
the action of the Hecke algebra generator $T$ was given on $\VA\otimes\VA$,
and the space $N$ was defined by
\begin{displaymath}
N=\Ker(T+1).
\end{displaymath}
It was also noted that $N=U'_q(\widehat{\SL}_2)\cdot \v{0}\otimes\v{0}$.
In this paper, we define, in general,
\begin{equation}
\lb{N}
N=\UPG[z\otimes z,z^{-1}\otimes z^{-1},z\otimes1+1\otimes z]\cdot v\otimes v,
\end{equation}
where $v$ is an extremal vector in $\VA$
(see~\S\ref{subsec:perfect-crystal} for the definition).  For $l=1$,
any $z^n\v{i}$ $(n\in\Z,i=0,1)$ is extremal. For $l=2$, we take
\begin{displaymath}
V=\Q \v{0}\oplus\Q \v{1}\oplus\Q \v{2}.
\end{displaymath}
The extremal vectors are $z^n\v{0}$ and $z^n\v{2}$ $(n\in\Z)$.  For
$l=1$, in the $q=1$ limit, the construction gives rise to ordinary
wedges with anti-commutation relations
\begin{displaymath}
z^m\v{i}\wedge z^n\v{j}+z^n\v{j}\wedge z^m\v{i}=0.
\end{displaymath}
For $l=2$, this is not the case, e.g.~$\v{1}\wedge \v{1}\not=0$,
even in the $q=1$ limit.

The definition \refeq{N} is appropriate for computational use. For theoretical
use, we have the following equivalent definition
\begin{displaymath}
N=\Ker(R-1).
\end{displaymath}
Here $R$ is the R-matrix acting on $\VA\otimes\VA$
(strictly speaking, the image of $R$ belongs to a certain completion of
$\VA\otimes\VA$).

The R-matrix satisfies the Yang-Baxter equation
\begin{displaymath}
R_{12}R_{23}R_{12}=R_{23}R_{12}R_{23},
\end{displaymath}
commutes with the $\UG$-action on $\VA\otimes\VA$, satisfies
\begin{displaymath}
R(z\otimes1)=(1\otimes z)R,\quad R(1\otimes z)=(z\otimes1)R,
\end{displaymath}
and is normalized as
\begin{displaymath}
R(v\otimes v)=v\otimes v,
\end{displaymath}
where $v$ is an extremal vector.

\subsection{Energy function and the normal ordering rules.}
In~\cite{KMS}, it was shown that the $q$-wedge relations give a normal
ordering rule of products of vectors.  Define $u_m$ ($m\in\Z$) by
\begin{equation}
\lb{U}
z^n\v{i}=u_{2n-i}.
\end{equation}
It was shown that the vectors
\begin{displaymath}
u_{m_1}\wedge\cdots\wedge u_{m_r}\quad(m_1<\cdots<m_r)
\end{displaymath}
form a base of $\bigwedge^r\VA$.

To describe the normal ordering rules in the general case, we use the
energy function
\begin{displaymath}
H:\BA\otimes\BA\longrightarrow\Z.
\end{displaymath}
The set $\BA$ is the crystal of $\VA$. For each element $b$ in $\BA$,
we have a corresponding vector $G(b)$ in $\VA$. In this section
we use the same symbol for $b$ and $G(b)$: e.g.~a general element of
$\BA$ for the level~$1$ $A^{(1)}_1$ case and  that of $\VA$
are denoted by $z^n\v{i}$. The energy function $H$ is such that
\begin{displaymath}
R\bigl(G(b_1)\otimes G(b_2)\bigr)
=z^{H(b_1\otimes b_2)}G(b_1)\otimes z^{-H(b_1\otimes b_2)}G(b_2)
\bmod qL(\VA)\otimes L(\VA),
\end{displaymath}
where $L(\VA)$ is the free module generated by $G(b)$ $(b\in\BA)$
over $A\deq\{f\in\Q(q);\text{ $f$ is regular at $q=0$}\}$.

For the level 2 $A^{(1)}_1$ case,
\begin{displaymath}
\BA=\{z^m\v{i};m\in\Z,i=0,1,2\}
\end{displaymath}
and
\begin{displaymath}
H(z^m\v{i}\otimes z^n\v{j})=-m+n+h_{ij}
\end{displaymath}
where the $(h_{ij})_{i,j=0,1,2}$ are given by
\begin{displaymath}
\bordermatrix
{&j=0&\phantom{j\,}1\phantom{=}&\phantom{j}2\phantom{j}\cr
i=0&0&0&0\cr
\phantom{j=\,}1&1&1&0\cr
\phantom{j=\,}2&2&1&0\cr}.
\end{displaymath}

We show that the set of vectors
\begin{displaymath}
G(b_1)\wedge\cdots\wedge G(b_r)
\end{displaymath}
such that
\begin{equation}
H(b_i\otimes b_{i+1})>0\quad(i=1,\dots,r-1)
\lb{NOV}
\end{equation}
is a base of $\bigwedge^r\VA$.

The vectors satisfying \refeq{NOV} are called normally ordered wedges.
To show that the normally ordered wedges span the $q$-wedge space, we
need to write down the basic $q$-wedge relations explicitly.  This
part of the work is technically much involved. We do it case by
case.  The generality in handling examples in this paper is narrower
than that of~\cite{KMN2} because of this limitation.

In~\cite{KMS} the linear independence of the normally ordered wedges
is proved by reduction to the $q=1$ limit.
Since the $q=1$ result is not known for the general case, we prove the
linear independence directly by using the Yang-Baxter equation for $R$
and the crystal base theory.

\subsection{Fock representations}
In~\cite{KMS}, the Fock spaces are constructed by means of an
inductive limit of $\bigwedge^r\VA$. In the case of
level~$1$ $A^{(1)}_1$, we take the sequence
$(u_m)_{m\in\Z}$ as in \refeq{U}.  The Fock space $\F_m$ is defined as
the space spanned by the semi-infinite wedges
\begin{displaymath}
u_{j_1}\wedge u_{j_2}\wedge u_{j_3}\wedge\cdots
\end{displaymath}
such that $j_k=m+k-1$ for sufficiently large $k$.
The action of $U_q(\widehat{\SL}_2)$ on $\F_m$
is defined by using the semi-infinite coproduct.
It was shown that $\F_m$ is the tensor product
\begin{displaymath}
V(\lambda_m)\otimes\C[H_-].
\end{displaymath}
Here $V(\lambda_m)$ is the irreducible highest weight
representation with the highest weight $\lambda_m$,
where
\begin{displaymath}
  \lambda_m=
  \begin{cases}
    \Lambda_1&\text{if $m\equiv0\bmod2$;}\cr
    \Lambda_0&\text{if $m\equiv1\bmod2$,}
  \end{cases}
\end{displaymath}
and $\C[H_-]$ is the Fock space of the Heisenberg algebra
generated by $B_n$ $(n\in\Z\backslash\{0\})$ that
acts on $\F_m$ by
\begin{displaymath}
B_n=\sum_{k=1}^\infty1\otimes\cdots\otimes1\otimes
{\buildrel{\scriptscriptstyle k\atop\scriptscriptstyle\vee}\over{z^n}}
\otimes1\otimes\cdots.
\end{displaymath}

To construct the Fock spaces in the general case, we
use the construction of affine crystals developed
in~\cite{KMN1}. We assume that $V$ has a perfect crystal $B$ of level
$l$.  Then we can choose a sequence $b^\circ_m$ in $\BA$ such that
\begin{gather*}
\langle c,\varepsilon(b^\circ_m)\rangle=l,\cr
\varepsilon(b^\circ_m)=\varphi(b^\circ_{m+1}),\cr
H(b^\circ_m\otimes b^\circ_{m+1})=1
\end{gather*}
(see subsection~\ref{subsec:perfect-crystal} for the definition of
$\varepsilon(b)$ and $\varphi(b)$).  In the case of
level 2 $A^{(1)}_1$, we have
\begin{equation}
\lb{VAC} b^\circ_m=
\begin{cases}
  z^k\v{j}&\text{if $m$ is odd;}\\
  z^{k+1-j}\v{2-j}&\text{if $m$ is even,}
\end{cases}
\end{equation}
for some $k\in\Z$ and $j\in\{0,1,2\}$ independent of $m$.
Then we shall define the Fock space $\F_m$ as a certain quotient of the
space spanned by the semi-infinite wedges
\begin{displaymath}
G(b_1)\wedge G(b_2)\wedge G(b_3)\wedge\cdots
\end{displaymath}
such that $b_n=b^\circ_{m+n-1}$ for sufficiently large $n$.
In particular, the Fock space contains the highest weight vector
\begin{displaymath}
|m\rangle=G(b^\circ_m)\wedge G(b^\circ_{m+1})\wedge
G(b^\circ_{m+2})\wedge\cdots
\end{displaymath}
with the highest weight
\begin{displaymath}
\lambda_m=
\begin{cases}
  j\Lambda_1+(2-j)\Lambda_0&\text{if $m$ is odd;}\cr
  (2-j)\Lambda_1+j\Lambda_0&\text{if $m$ is even.}
\end{cases}
\end{displaymath}

The quotient is such that if
\begin{displaymath}
H(b\otimes b^\circ_m)\le0
\end{displaymath}
we require that
\begin{displaymath}
G(b)\wedge|m\rangle=0.
\end{displaymath}
Here is a significant difference between level~$1$
$A^{(1)}_n$ and other cases. For the former
if $H(b\otimes b^\circ_m)\le0$ then
\begin{displaymath}
G(b)\wedge G(b^\circ_m)\wedge\cdots\wedge G(b^\circ_{m'})=0
\end{displaymath}
for sufficiently large $m'$. But, this is not true in general. The
correct statement is that for any $n$ we can find $m'$
such that the $q$-wedge $G(b)\wedge
G(b^\circ_m)\wedge\cdots\wedge G(b^\circ_{m'})$ is a linear
combination of normally ordered wedges whose
coefficients are $O(q^n)$ at $q=0$.
Therefore, we need to impose the separability of the $q$-adic topology,
taking the quotient by the closure of $\{0\}$.

It is necessary to check that the action of $\UG$ given by the semi-infinite
coproduct, is well-defined. A careful study of the $q$-wedges shows that
\begin{equation}
\lb{FI}
\Delta^{(\infty/2)}(f_i)|m\rangle=G(\tilde f_ib^\circ_m)\wedge|m+1\rangle,
\end{equation}
where
\begin{displaymath}
\Delta^{(\infty/2)}(f_i)=\sum_{n=1}^\infty
1\otimes\cdots\otimes1\otimes
{\buildrel{\scriptscriptstyle n\atop\scriptscriptstyle\vee}\over{f_i}}
\otimes t_i\otimes t_i\otimes\cdots.
\end{displaymath}
In the case in~\cite{KMS}, the action of $\Delta^{(\infty/2)}(f_i)$ on
each vector in $\F_m$ is such that only finitely many terms in the sum
are different from $0$. This is not true in general.  For example,
consider the case $k=1$ and $j=1$ in~\refeq{VAC}. We have
$f_1\v{1}=[2]\v{2}$ $([2]=q+q^{-1})$ and $t_1|m\rangle=q|m\rangle$.
Therefore, we have
\begin{multline*}
\Delta^{(\infty/2)}(f_1)(\v{1}\wedge \v{1}\wedge \v{1}\wedge\cdots)=
q[2](\v{2}\wedge \v{1}\wedge \v{1}\wedge\cdots)
+q[2](\v{1}\wedge \v{2}\wedge \v{1}\wedge\cdots)\cr
+q[2](\v{1}\wedge \v{1}\wedge \v{2}\wedge\cdots)+\cdots\cdots.
\end{multline*}
On the other hand, we have
\begin{displaymath}
\v{1}\wedge \v{2}+q^2\v{2}\wedge \v{1}=0,
\end{displaymath}
and hence
\begin{displaymath}
\Delta^{(\infty/2)}(f_1)(\v{1}\wedge \v{1}\wedge\cdots)
=\v{2}\wedge \v{1}\wedge \v{1}\wedge\cdots,
\end{displaymath}
by summing up
\begin{displaymath}
1+(-q^2)+(-q^2)^2+\cdots={1\over1+q^2}
\end{displaymath}
in the $q$-adic topology.

In general, based on \refeq{FI} we can show the well-definedness
of the $\UG$-action.

The decomposition of the $q$-Fock spaces into the irreducible
$\UG$-modules and the bosonic Fock space goes the same as the
level~$1$ $A^{(1)}_n$ case.  We carry out the computation of the exact
commutation relations of the bosons in each case, by reducing it to
the commutation relations of vertex operators.

The plan of this paper is as follows.  We list the notations in
section 2.  We define the finite $q$-wedges in section 3 and prove
that the normally ordered wedges form a base. In section 4, we define
the $q$-Fock space and the actions of $\UG$ and the Heisenberg
algebra.  We give level~$1$ examples in section 5 for which we check
the conditions assumed in section 3. We compute the level~1 two point
functions in section 6 in order to find the commutation relations of
the bosons. Section 7 is devoted to a higher level example.  We add
four appendices. In Appendix A we prove a proposition on crystal base
which is necessary in this paper but was not proved in~\cite{KMN1}.
Appendix B is a proof that the Serre relations follow from the
integrability of representations. Appendix C is the computation of the
two-point correlation functions of the $q$-vertex operators in the
$D^{(2)}_{n+1}$ case. In Appendix D we consider the $q\rightarrow 1$
limit for the $A^{(2)}_{2n}$ case and compare it to the result in
\cite{JY}.

{\it Acknowledgement.} We thank Eugene Stern for discussions in the
early stage of this work, and Masato Okado for explaining the method
for the computation of the result in Appendix~C.  This work is partially
supported by Grant-in-Aid for Scientific Research on Priority Areas,
the Ministry of Education, Science, Sports and Culture.  J.-U.~H.~P.\
and C.~M.~Y.\ are supported by the Japan Society for the Promotion of
Science.


\newsection{Preliminary}
\subsection{Notations}
\label{2.1}

In this paper we use the following notations.
\medskip
\halign{\hfill$#$\hfill&#\ &#\hfill\cr
\delta(P)&=&
$\begin{cases}
  1& \text{if a statement $P$ is true}\\
  0& \text{if $P$ is false}.
\end{cases}$\cr
\akete
\Gg&:&an affine Lie algebra.\cr
\Gh&:&its Cartan subalgebra with dimension $\rk(\Gg)+1$.\cr
I&:& the index set for simple roots.\cr
\alpha_i&:& a simple root $\in \Gh^*$ corresponding to $i\in I$.\cr
h_i&:&a simple coroot $\in\Gh$ corresponding to $i\in I$.\cr
\akete
&& We assume that the simple roots and the simple
coroots are \cr
&& linearly independent.\cr
\akete
W&:&the Weyl group of $\Gg$.\cr
(\ ,\ )&:& a $W$-invariant non-degenerate bilinear symmetric form on $\Gh^*$\cr
&&such that $(\alpha_i,\alpha_i)\in 2\BZ_{>0}$.\cr
\pairing{\ ,\ }&:& the coupling $\Gh\times\Gh^*\mapto\C$.\cr
P&:&a weight lattice $\subset \Gh^*$.\cr
Q&=&$\sum_i\BZ\alpha_i$ the root lattice.\cr
Q_\pm&=&$\pm\sum_i\BZ_{\ge0}\alpha_i$.\cr
\delta&:&an element of $Q_+$ such that
$\BZ\delta=\{\lam\in Q;\lan h_i\lam\ran=0\}$.\cr
c&:& an element of $\sum_i\BZ_{>0}h_i$ such that
$\BZ c=\{h\in\sum_i\BZ h_i; \lan h,\alpha_i\ran=0\}$.\cr
\akete
\noalign{\quad\hbox{We write}}
\akete
\delta&=&$\sum_ia_i\alpha_i\quad$ and \cr
c&=&$\sum_ia_i^\vee h_i$.\cr
P_{\cl}&=&$P/\BZ\delta$.\cr
\hfill \cl\ &:&$P\to P_{\cl}$.\cr
&&We assume for the sake of simplicity\cr
&&\qquad$P_\cl\buildrel\sim\over\to\Hom_\BZ(\oplus_{i\in I}\BZ h_i,\BZ)$.\cr
&&This implies $\{\lam\in P\,;\lan h_i,\lam\ran=0$ for any $i\in
I\}=\BZ\delta$.\cr
\akete
\Lambda_i&:&a fundamental weight in $P$,\cr
&&i.e. an element of $P$ such that
$\lan h_j,\Lambda_i\ran=\delta_{ij}$.\cr
\Lambda_i^\cl&=&$\cl(\Lambda_i)$, the fundamental weight in $P_\cl$.\cr
&&Note that $\Lambda_i$ is determined modulo $\BZ\delta$.\cr
\akete
P^0&:&the level $0$ part of $P$, i.e.
$\{\lam\in P:\lan c,\lam\ran=0\}$.\cr
P_{\cl}^0&:&the level $0$ part of $P_{\cl}$, i.e.
$\cl(P^0)$.\cr
U_q(\Gg)&:& the quantized universal enveloping algebra\cr
&& with $\{q^h;h\in P^*\}$ as its Cartan part.\cr
U'_q(\Gg)&:& the quantized universal enveloping algebra\cr
&& with $\{q^h;h\in P_{\cl}{}^*\}$ as its Cartan part.\cr
&&Hence $U'_q(\Gg)$ is a subalgebra of $U_q(\Gg)$.\cr
K&=&$\BQ(q)$.\cr
\akete
&&We consider $U_q(\Gg)$ and $U'_q(\Gg)$ over $K$.\cr
\akete
A&=&$\{f\in K$; $f$ has no pole at $q=0\}\,$.\cr
U'_q(\Gg)_\BZ&:&the $\KZ$-subalgebra of $U'_q(\Gg)$ generated by
the divided powers\cr
&& $e_i^{(n)}$, $f_i^{(n)}$, $t_i$ and $\left\{{t_i\atop n}\right\}$.\cr
U_q(\Gg)_\BZ&:&the $\KZ$-subalgebra of $U_q(\Gg)$ generated by
$U'_q(\Gg)_\BZ$\cr
&&and $\left\{{q^h\atop n}\right\}$ ($h\in P^*$).\cr
}

\medskip
The quantized affine algebra
$U_q(\Gg)$ is a $K$-algebra generated by
$e_i,f_i\,(i\in I)$ and $q^h\,(h\in P^*)$ with the commutation relations
\beqn
&&q^h=1\ \hbox{for $h=0$,}\\
&&q^{h+h'}=q^h q^{h'}\ \hbox{for $h,h'\in P^*$,}\\
&&q^he_iq^{-h}=q^{\lan h,\alpha_i\ran} e_i\qbox{and}
q^hf_iq^{-h}=q^{-\lan h,\alpha_i\ran}f_i,\\
&&[e_i,f_j]=\delta_{ij}{{t_i-t_i^{-1}}\over{q_i-q_i^{-1}}},\\
&&\hbox{for $i\not=j\in I$}\\
&&\begin{matrix}
&\quad\sum_k(-1)^ke_i^{(k)}e_je_i^{(-\lan h_i,\alpha_j\ran-k)}=0,\cr
&\quad\sum_k(-1)^kf_i^{(k)}f_jf_i^{(-\lan h_i,\alpha_j\ran-k)}=0.
\end{matrix}
\endeqn
Here
$$q_i=q^{{(\alpha_i,\alpha_i)\over2}}\hbox{ and }
t_i=q^{{{(\alpha_i,\alpha_i)}\over2}h_i}.$$
\subsection{Coproducts}
There are several coproducts of $\U$ used in the literature.
In this paper, we use a coproduct different from
the ones used in~\cite{DJO,JM,K1,K2,KMN1}.
In this subsection, we shall explain the relations
among four coproducts:
\beq
\Delta_+:&&
\left\{
\begin{aligned}
  q^h &\mapsto q^h\otimes q^h\cr
  e_i &\mapsto e_i\otimes 1+t_i\otimes e_i\cr
  f_i &\mapsto f_i\otimes t_i^{-1}+1\otimes f_i
\end{aligned} \right.\\[5pt]
\Delta_-:&&
\left\{
\begin{aligned}
  q^h &\mapsto q^h\otimes q^h\cr
  e_i &\mapsto e_i\otimes t_i^{-1}+1\otimes e_i\cr
  f_i &\mapsto f_i\otimes 1+t_i\otimes f_i
\end{aligned} \right.\\[5pt]
\bar\Delta_+:&&
\left\{
\begin{aligned}
q^h&\mapsto q^h\otimes q^h\cr
e_i&\mapsto e_i\otimes 1+t_i^{-1}\otimes e_i\cr
f_i&\mapsto f_i\otimes t_i+1\otimes f_i
\end{aligned} \right.\\[5pt]
\bar\Delta_-:&&
\left\{
\begin{aligned}
q^h&\mapsto q^h\otimes q^h\cr
e_i&\mapsto e_i\otimes t_i+1\otimes e_i\cr
f_i&\mapsto f_i\otimes 1+t_i^{-1}\otimes f_i
\end{aligned} \right.
\endeq

Their antipodes are given by
\beq
a_+:&&
\left\{
\begin{aligned}
q^h&\mapsto q^{-h}\cr
e_i&\mapsto -t_i^{-1}e_i\cr
f_i&\mapsto -f_it_i
\end{aligned} \right.\\[5pt]
a_-:&&
\left\{
\begin{aligned}
q^h&\mapsto q^{-h}\cr
e_i&\mapsto -e_it_i\cr
f_i&\mapsto -t_i^{-1}f_i
\end{aligned} \right.\\[5pt]
\bar a_+:&&
\left\{
\begin{aligned}
q^h&\mapsto q^{-h}\cr
e_i&\mapsto -t_ie_i\cr
f_i&\mapsto -f_it_i^{-1}
\end{aligned} \right.\\[5pt]
\bar a_-:&&
\left\{
\begin{aligned}
q^h&\mapsto q^{-h}\cr
e_i&\mapsto -e_it_i^{-1}\cr
f_i&\mapsto -t_if_i
\end{aligned} \right.
\endeq

For two $\U$-modules
$M_1$ and $M_2$, let us denote by
$M_1\otimes_+M_2$, $M_1\otimes_-M_2$
$M_1\btensor_+M_2$ and $M_1\btensor_-M_2$
the vector space $M_1\otimes_KM_2$
endowed with the $\U$-module structure
via the coproduct $\Delta_+$,
$\Delta_-$, $\bar\Delta_+$ and $\bar\Delta_-$,
respectively.

We have functorial isomorphisms of $\U$-modules
\beq
M_1\otimes_+ M_2&\isomo&M_2\btensor_-M_1\\
M_1\otimes_- M_2&\isomo&M_2\btensor_+M_1
\endeq
by
$u_1\otimes u_2\mapsto u_2\otimes u_1$.

We have functorial isomorphisms of $\U$-modules
\beq
q^{-(\cdot,\cdot)}:M_1\otimes_+M_2&\isomo&M_1\otimes_-M_2\\
q^{(\cdot,\cdot)}:M_1\btensor_+M_2&\isomo&M_1\btensor_-M_2
\endeq
Here $q^{-(\cdot,\cdot)}$
sends $u_1\otimes_+u_2$ to
$q^{-(\wt(u_1),\wt(u_2))}u_1\otimes_-u_2$
and $q^{(\cdot,\cdot)}$ sends $u_1\btensor_+u_2$
to $q^{(\wt(u_1),\wt(u_2))}u_1\btensor_-u_2$.

The tensor products $\otimes_+$ and $\btensor_-$
behave well under upper crystal bases
and $\otimes_-$ and $\btensor_+$ behave well under lower crystal bases.
Namely, if $(L_j,B_j)$ is an upper crystal base
of an integrable $\U$-module $M_j$ ($j=1,2$),
then
$(L_1\otimes_A L_2, B_1\otimes B_2)$ is an upper crystal base
of $M_1\otimes_+M_2$ and $M_1\btensor_-M_2$.
Similarly, if $(L_j,B_j)$ is a lower crystal base
of $M_j$, then
$(L_1\otimes_A L_2, B_1\otimes B_2)$ is a lower crystal base
of $M_1\otimes_-M_2$ and $M_1\btensor_+M_2$.
If we use $\otimes_+$ or $\otimes_-$, the tensor product of crystal base is
described as follows.
For two crystals $B_1$, $B_2$
and $b_1\in B_1$, $b_2\in B_2$,
\beqn
&&\wt(b_1\otimes b_2)=\wt(b_1)+\wt(b_2),\\
&&\varepsilon_i(b_1\otimes b_2)=
\max(\varepsilon_i(b_1),\varepsilon_i(b_2)-\lan h_i,\wt(b_1)\ran),\\
&&\varphi_i(b_1\otimes b_2)=
\max(\varphi_i(b_1)+\lan h_i,\wt(b_2)\ran,\varphi_i(b_2)),\\
&&\te_i(b_1\otimes b_2)=
\begin{cases}
\te_ib_1\otimes b_2&\text{if $\varphi_i(b_1)\ge\varepsilon_i(b_2)$},\cr
b_1\otimes\te_ib_2&\text{if $\varphi_i(b_1)<\varepsilon_i(b_2)$,}
\end{cases}\\
&&\tf_i(b_1\otimes b_2)=
\begin{cases}
\tf_ib_1\otimes b_2&\text{if $\varphi_i(b_1)>\varepsilon_i(b_2)$},\cr
b_1\otimes\tf_ib_2&\text{if $\varphi_i(b_1)\le\varepsilon_i(b_2)$}.
\end{cases}
\endeqn

If we use the other tensor products $\btensor_+$ or $\btensor_-$,
we have to exchange the first and the second factors in the formulas above.
Namely the tensor product of crystals is given as follows.
\beq
\label{rule}
\wt(b_1\otimes b_2) &=& \wt(b_1)+\wt(b_2),\cr
\akete
\varepsilon_i(b_1\otimes b_2) &=&
\max(\varepsilon_i(b_1)-\lan h_i,\wt(b_2)\ran,\varepsilon_i(b_2)),\cr
\akete
\varphi_i(b_1\otimes b_2) &=&
\max(\varphi_i(b_1),\varphi_i(b_2)+\lan h_i,\wt(b_1)\ran),\cr
\akete
\te_i(b_1\otimes b_2) &=&
\begin{cases}
  \te_ib_1\otimes b_2 & \text{if $\varepsilon_i(b_1)>\varphi_i(b_2)$},\cr
  b_1\otimes\te_ib_2 & \text{if $\varepsilon_i(b_1)\le\varphi_i(b_2)$},
\end{cases}
\cr
\akete
\tf_i(b_1\otimes b_2) &=&
\begin{cases}
  \tf_ib_1\otimes b_2&\text{if $\varepsilon_i(b_1)\ge\varphi_i(b_2)$},\cr
  b_1\otimes\tf_ib_2&\text{if $\varepsilon_i(b_1)<\varphi_i(b_2)$}.
\end{cases}
\endeq

In this article, we mainly use the tensor product
$\btensor_+$ and lower crystal bases.
The rule of the tensor product of crystals is therefore by
(\ref{rule}).
Note that $\otimes_+$ is used in~\cite{DJO,JM} and
$\otimes_-$ in~\cite{K2,KMN1}.

\newsection{Wedge products}\label{3w}
\subsection{Perfect crystal}\label{subsec:perfect-crystal}
Let us take an integrable finite-dimensional representation $V$ of
$U'_q(\Gg)$.
Let $V=\oplus_{\lam\in P^0_{\cl}}V_\lam$ be its weight space decomposition.
Its affinization is
defined by
$$V_{\aff}=\bigoplus_{\lam\in P}(V_{\aff}){}_\lam$$
where $(V_{\aff}){}_\lam=V_{\cl(\lam)}$ for $\lam \in P$.
Let $\cl:(V_{\aff}){}_\lam\to V_{\cl(\lam)}$ denote
the canonical isomorphism.
Then $V_{\aff}$ has a natural structure of a $U_q(\Gg)$-module
such that $\cl:V_\aff\to V$ is $U'_q(\Gg)$-linear
(see~\cite{KMN1}).

Let $z:V_{\aff}\to V_{\aff}$
be the endomorphism of weight $\delta$ given by
\begin{displaymath}
  \begin{CD}
    (V_{\aff}){}_\lam  @>{z}>>  (V_{\aff}){}_{\lam+\delta}\\
    @V{\wr}VV  @V{\wr}VV \\
    V_{\cl(\lam)} @= V_{\cl(\lam+\delta)}
  \end{CD}
\end{displaymath}
The endomorphism $z$ is $U'_q(\Gg)$-linear.

Taking a section of
$\cl:P\to P_\cl$, $V_{\aff}$ may be identified with
$V\otimes \BC[z,z^{-1}]$ (see section \ref{GOITI}).

We assume that
\begin{description}
\item[(P)] {\sl $V$ has a {\sl perfect} crystal base $(L,B)$.}
\end{description}

Let us recall its definition in~\cite{KMN1}.
A crystal base $(L,B)$ is called perfect of level $l\in\BZ_{>0}$
if it satisfies the following axioms (P1)--(P3).
\begin{description}
\item[(P1)]\enspace There is a weight $\lam^\circ\in P_{\cl}^0$ such that
the weights of $V$ are contained in the convex hull of $W\lam^\circ$
and that $\dim V_{w\lam^\circ}=1$ for any $w$ in the Weyl group $W$.
\hb
We call a vector in $V_{w\lam^\circ}$ an {\sl extremal} vector
with extremal weight $w\lam^\circ$.
\item[(P2)]\enspace$B\otimes B$ is connected.
\item[(P3)]\enspace There is a positive integer $l$
satisfying the following conditions.
\smallskip
\begin{enumerate}
\item 
For every $b\in B$,
$\lan c,\varepsilon(b)\ran=\lan c,\varphi(b)\ran\ge l$.
Here we set
\beq
\labeleq{EP}
\varepsilon(b)&=&\sum_{i\in I}\varepsilon_i(b)\Lam_i^\cl\in P_\cl\cr
\akete
\varphi(b)&=&\sum_{i\in I}\varphi_i(b)\Lam_i^\cl\in P_\cl
\endeq
with the fundamental weights $\Lam_i^\cl\in P_\cl$.
\item 
Set $B_{\min}=\{b\in B;\lan c,\varepsilon(b)\ran=l\}$
and $(P_\cl^+){\raise-1pt\hbox{${}_l$}}=\{\lam\in P_\cl\,;\,
\lan c,\lambda\ran=l \hbox{ and }\lan h_i,\lam\ran\ge0
\hbox{ for every $i\in I$}\}$.
Then
\beqn
\varepsilon,\varphi:B_{\min}\to (P_\cl^+){\raise-1pt\hbox{${}_l$}}\qbox{are
bijective.}
\endeqn
\end{enumerate}
\end{description}
Note that (P1) is equivalent to the irreducibility of $V$ (see~\cite{CP}).

Note that the equality
$\lan c,\varepsilon(b)\ran=\lan c,\varphi(b)\ran$
in (P3)~(i) follows from
$$\varphi(b)=\wt(b)+\varepsilon(b)$$
and the fact that
$V$ is a $U'_q(\Gg)$-module of level $0$.

\medskip
\Remark The map $\varepsilon(b)\mapsto \varphi(b)$ ($b\in B_{\min}$)
defines an automorphism of $(P^+_\cl){\raise-1pt\hbox{${}_l$}}$.
In all the examples of perfect crystals that we know,
this automorphism is induced by
a Dynkin diagram automorphism.

\medskip
We have constructed $V_{\aff}$ out of $V$.
Similarly we construct
the crystal base $(L_{\aff},B_{\aff})$
of $V_{\aff}$ out of $(L,B)$.
We define similarly $\cl:B_\aff\to B$
and $z: B_\aff\to B_\aff$.

We assume further that $V$ has a good base $\{G(b)\}_{b\in B}$
in the following sense.

\begin{description}
\item[(G)]{\sl \enspace $V$ has a lower global base
$\{G(b)\}_{b\in B}$.}
\end{description}

This means that the base $\{G(b)\}_{b\in B}$
satisfies the following conditions (cf.~\cite{K2}).
\begin{enumerate}
\item 
$\bigoplus\limits_{b\in B} \KZ G(b)$
is a $U'_q(\Gg)_\BZ$-submodule of $V$.
\item 
$b\equiv G(b)\ \mod L/qL$.
\item 
$e_iG(b)=[\varphi_i(b)+1]_iG(\te_ib)+\sum E^i_{b,b'}G(b')$,
\item 
$f_iG(b)=[\varepsilon_i(b)+1]_iG(\tf_ib)+\sum F^i_{b,b'}G(b')$.
\end{enumerate}
In both cases, the sum ranges over $b'$ that belongs to an $i$-string
strictly longer than that of $b$ ($\Leftrightarrow$
$\varepsilon_i(b')\ge\varepsilon_i(b)$ or $\varphi_i(b')\ge
\varphi_i(b)$ according to (iii) or (iv)).  Moreover the coefficients
satisfy \beq
E^i_{b,b'}&\in& q q_i^{-\varphi_i(b')}\AZ\cup q^{-1}q_i^{\varphi_i(b')}
\Z[q^{-1}]\\
F^i_{b,b'}&\in&q q_i^{-\varepsilon_i(b')}
\AZ\cup q^{-1}q_i^{\varepsilon_i(b')}\Z[q^{-1}].  \endeq

\Remark The reason why we choose a lower global base is explained in
Theorem~\ref{kern} and the remark after Proposition \ref{vertex}.

\medskip
We define the base $\{G(b)\}_{b\in B_\aff}$ of $V_\aff$
by $\cl(G(b))=G(\cl(b))$.
We have $G(z^nb)=z^nG(b)$  for $n\in\BZ$ and $b\in B_\aff$.

\smallskip
\subsection{Energy function}\label{sec:ener}
Let $H$ be an energy function
(see~\cite{KMN1}).
Namely
$H:B_{\aff}\otimes B_{\aff}\to\BZ$ satisfies
\begin{description}
\item[(E1)]  $H(zb_1\otimes b_2)=H(b_1\otimes b_2)-1$.
\item[(E2)] $H(b_1\otimes zb_2)=H(b_1\otimes b_2)+1$.
\item[(E3)] $H$ is constant on every connected component of the crystal graph
$B_{\aff}\otimes B_{\aff}$.
\end{description}

By (E1--3), $H$ is uniquely determined up to a constant. We normalize $H$ by

\begin{description}
\item[(E4)] $H(b\otimes b)=0$ for any (or equivalently some)
extremal $b\in B_{\aff}$ (i.e. $\cl\bigl(\wt(b)\bigr)\in W\lam^\circ$).
\end{description}
We know already its existence and uniqueness (\cite{KMN1}).
The existence is in fact proved by using R-matrix.
Let us explain their relation.
There is a $U_q(\Gg)$-linear endomorphism (R-matrix)
$R$ of $V_\aff\otimes V_\aff$
such that
\begin{align}
  R\circ(z\otimes 1)&=(1\otimes z)\circ R\\
  R\circ(1\otimes z)&=(z\otimes 1)\circ R\\
  \intertext{and normalized by}
  R(u\otimes u)&=u\otimes u \qquad\text{for every extremal $u\in
    V_\aff$.}
\end{align}
Strictly speaking, $R$ is a homomorphism from $V_\aff\otimes V_\aff$
to its completion $V_\aff\widehat\otimes V_\aff$.
It is proved in~\cite{KMN1} that
$R$ sends $L_\aff\otimes L_\aff$ to
$L_\aff{\widehat\otimes}L_\aff$ and
\beq
&&R(G(b_1)\otimes G(b_2))\nn\\
&&\quad\equiv
G(z^{H(b_1\otimes b_2)}b_1)\otimes G(z^{-H(b_1\otimes b_2)}b_2)\ \mod\,
qL_\aff{\widehat\otimes}L_\aff\label{H-fct}\\
\hfill &&\qquad\qquad\text{for every $b_1$, $b_2\in B_\aff$.}\nn
\endeq
We know that $R$ has finitely many poles.
It means that there is a non-zero
$\psi\in K[z\otimes z^{-1},z^{-1}\otimes z]$
such that $\psi R$ sends $V_\aff\otimes V_\aff$ into itself.
We assume that the denominator $\psi$ of $R$ satisfies
the following property.

\begin{description}
\item[(D)]{\sl \enspace $\psi\in A[z\otimes z^{-1}]$
and $\psi=1$ at $q=0$.}
\end{description}
\medskip
We take a linear form $s:P\to\BQ$ such that $s(\alpha_i)=1$ for every $i\in I$,
and define
$$l:B_{\aff}\to\BZ$$
by $l(b)=s(\wt(b))+c$ for some constant $c$.
With a suitable choice of $c$, $l$ is $\BZ$-valued.
It satisfies

\begin{enumerate}
\item 
$l(zb)=l(b)+a$ for any $b\in B_{\aff}$.
Here $a$ is a positive integer independent of $b$.
\item 
$l(\te_ib)=l(b)+1$ if $i\in I$ and $b\in B_{\aff}$
satisfy $\te_ib\not=0$.
\end{enumerate}

We assume that it satisfies

\medskip\noindent
\begin{description}
\item[(L)] {\sl If $H(b_1\otimes b_2)\le 0$, then $l(b_1)\ge l(b_2)$.}
\end{description}

\subsection{Wedge products} \label{subsec:wedge-prod}
\def\tR{\tilde R}
We define $L(V_\aff^{\otimes 2})$ by $L_\aff\otimes_AL_\aff$.
Let us set $\tR=\psi(z\otimes 1,1\otimes z)R
=R\psi(1\otimes z,z\otimes 1)$.
Then it is an endomorphism of $V_\aff^{\otimes2}$
and $L(V_\aff^{\otimes2})$ is stable by $\tR$.
We shall denote by the same letter $\tR$
the endomorphism
of $L(V_\aff^{\otimes2})/qL(V_\aff^{\otimes2})$
induced by $\tR$.
Then by (D) and (\ref{H-fct})
we have the equality in $L(V_\aff^{\otimes2})/qL(V_\aff^{\otimes2})$
\beq
\tR(b_1\otimes b_2)&=&
z^{H(b_1\otimes b_2)}b_1\otimes z^{-H(b_1\otimes b_2)}b_2
\qbox{for every $b_1$,$b_2\in B_\aff$.}
\label{RH}
\endeq
Since $R^2=1$, we have
\beq
\big(\tR-\psi(z\otimes 1,1\otimes z)\big)\circ
\big(\tR+\psi(1\otimes z,z\otimes 1)\big)&=&0.
\endeq

Let us choose an extremal vector $u\in V_{\aff}$. Then we define
$$
N=U_q(\Gg)[z\otimes z,\,z^{-1}\otimes z^{-1},\,z\otimes 1+1\otimes z]
\,(u\otimes u)\,.
$$
This definition does not depend on the choice of $u$,
because
an extremal vector $u$ of
weight $\lam$ satisfies
\beqn
(f_i^{(n)}u)\otimes(f_i^{(n)}u)&=&f_i^{(2n)}(u\otimes u)
\qbox{if $\lan h_i,\lam\ran=n\ge0$,}\\
(e_i^{(n)}u)\otimes(e_i^{(n)}u)&=&e_i^{(2n)}(u\otimes u)
\qbox{if $\lan h_i,\lam\ran=-n\le0$.}
\endeqn
By the definition, we have
\beq
&&f(z\otimes1,1\otimes z)N\subset N\nn\\
&&\phantom{----}
\hbox{for any symmetric Laurent polynomial $f(z_1,z_2)$.}
\label{symm}
\endeq

We put the following postulate.
\begin{description}
\item[(R)]
{\sl For every pair $(b_1,b_2)$ in $B_\aff$
with $H(b_1\otimes b_2)=0$,
there exists $C_{b_1,b_2}\in N$ which has the form
$$C_{b_1,b_2}=G(b_1)\otimes G(b_2)-
\sum\limits_{b'_1,b'_2} a_{b'_1,b'_2}\,G(b'_1)\otimes G(b'_2)\,.$$
Here the sum ranges over $(b'_1,b'_2)$ such that
\beqn
&H(b'_1\otimes b'_2)>0\,,\cr
&l(b_2)\le l(b'_1)<l(b_1)\,,\cr
&l(b_2)<l(b'_2)\le l(b_1)\,,
\endeqn
and the coefficients $ a_{b'_1,b'_2}$ belong to $\KZ$.}
\end{description}

Later in Lemma \ref{lem:s},
we see that $ a_{b'_1,b'_2}$ belong to $q\AZ$.

Since we have normalized the R-matrix by
$R(u\otimes u)=u\otimes u$, we have
\beq
\tR(v)&=&\psi(z\otimes 1,1\otimes z)v
\qbox{for every $v\in N$.
}\label{tr}
\endeq
Hence $\tilde{R}$ sends $N$ to itself.

We set
\beqn
L(N)&=&N\cap L(V_\aff^{\otimes 2}).
\endeqn
Then by (D) and (\ref{tr}), we have
the equality in $L(V_\aff^{\otimes 2})/qL(V_\aff^{\otimes 2})$
\beq
&&\tR(b)=b\qbox{for every $b\in L(N)/qL(N)$.}
\endeq

We define the wedge product by
$$\bw^2 V_{\aff}=V_{\aff}^{\otimes 2}/N\,.$$
For $v_1$,$v_2\in V$, let us denote
by $v_1\wedge v_2$
the element of $\bw^2 V_{\aff}$
corresponding to $v_1\otimes v_2$.
We set
$$
L(\bw^2V_{\aff})=L(V_\aff^{\otimes 2})/L(N)
\subset \bw^2 V_{\aff}.$$

Now we shall study the properties of $\bw^2 V_{\aff}$
under conditions (P), (G), (D), (L) and (R).

We conjecture that (P) and (G) imply the other conditions (D), (L) and (R).

\begin{lemma}\label{indp}
If\/ $\sum_{H(b_1\otimes b_2)>0}a_{b_1,b_2}G(b_1)\otimes G(b_2)$
belongs to $\Ker\big(\tR-\psi(z\otimes 1,1\otimes z)\big)$,
then all\/ $a_{b_1,b_2}$ vanish.
\end{lemma}

\proof
It is enough to show that for $n\in\BZ$
\beq\label{eq:12}
&\kern-30pt
\hbox{if $a_{b_1,b_2}\in q^nA$ for all $b_1,b_2$,
then $a_{b_1,b_2}\in q^{n+1}A$.}
\endeq
By (D), (\ref{tr}) and (\ref{RH}), we obtain
the identity in
$L(V_\aff^{\otimes 2})/qL(V_\aff^{\otimes 2})$,
$$\sum_{H(b_1\otimes b_2)>0}(q^{-n}a_{b_1,b_2})b_1\otimes b_2
=\sum_{H(b_1\otimes b_2)>0}(q^{-n}a_{b_1,b_2})
z^{H(b_1\otimes b_2)}b_1\otimes z^{-H(b_1\otimes b_2)}b_2.$$
Since $H\big(z^{H(b_1\otimes b_2)}b_1\otimes z^{-H(b_1\otimes b_2)}b_2\big)
=-H(b_1\otimes b_2)<0$, we obtain the desired assertion (\ref{eq:12}).
\qed

A similar argument leads to the following result.
\Lemma\label{lem:s}
If $H(b_1\otimes b_2)=0$ and
$G(b_1)\otimes G(b_2)-\sum\limits_{H(b'_1\otimes b'_2)>0}
a_{b'_1,b'_2}G(b'_1)\otimes G(b'_2)$
belongs to $N$, then $a_{b'_1,b'_2}\in qA$.
\end{lemma}

We shall call a pair $(b_1,b_2)$ of elements in $B_{\aff}$
{\it normally ordered}  and
$G(b_1)\wedge G(b_2)$ a {\sl normally ordered wedge}
if $H(b_1\otimes b_2)>0$.
The axiom (R) may be considered as a rule
to write $G(b_1)\wedge G(b_2)$ as a linear combination of
normally ordered wedges when $H(b_1\otimes b_2)=0$.
In order to treat the case
$H(b_1\otimes b_2)=-c<0$, we introduce
an element of $N$ (see (\ref{symm}))
\beq
C'_{b_1,b_2}
&=&(1\otimes z^{-c}+z^{-c}\otimes 1)C_{b_1,z^cb_2}\\
&=&(1\otimes z^c+z^c\otimes 1)C_{z^{-c}b_1,b_2}.\nn
\endeq
Note that
$H(b_1\otimes z^cb_2)=H(z^{-c}b_1\otimes b_2)=0$.

\Lemma\label{lem:2}
If $H(b_1\otimes b_2)\le0$, then
$C'_{b_1,b_2}$ has the form
$$G(b_1)\otimes G(b_2)-\sum\limits_{b'_1,b'_2}
a_{b'_1,b'_2}G(b'_1)\otimes G(b'_2).$$
Here
the sum ranges over $(b'_1,b'_2)$ such that
\beqn
&H(b'_1\otimes b'_2)>H(b_1\otimes b_2)\,,\cr
&l(b_2)\le l(b'_1)<l(b_1)\,,\cr
&l(b_2)<l(b'_2)\le l(b_1)\,.
\endeqn
Moreover $a_{b'_1,b'_2}$ belongs to $\AZ$.
\enlemma

\proof
Assume $H(b_1\otimes b_2)=-c<0$.
Set
$$C'_{z^{-c}b_1,b_2}=G(z^{-c}b_1)\otimes G(b_2)-\sum_{H(b'_1\otimes b'_2)>0}
a_{b'_1,b'_2}G(b'_1)\otimes G(b'_2).$$
Here the sum ranges over
\beqn
&l(b_2)\le l(b'_1)<l(z^{-c}b_1)\,,\cr
&l(b_2)<l(b'_2)\le l(z^{-c}b_1)\,.
\endeqn
Then
\beqn
C'_{b_1,b_2}&=&
G(b_1)\otimes G(b_2)+G(z^{-c}b_1)\otimes G(z^{c}b_2)\hfill\cr
&&\quad-\sum_{H(b'_1\otimes b'_2)>0}
a_{b'_1,b'_2}\big(G(b'_1)\otimes G(z^{c}b'_2)+G(z^{c}b'_1)\otimes G(b'_2)\big).
\endeqn
The desired properties can be easily checked.
\qed

By the repeated use of the proposition above, we obtain the following
result.

\Corollary\label{cor:4}
If $H(b_1\otimes b_2)\le0$ then
$N$ contains an element $ C_{b_1,b_2}$, which has the form
$$G(b_1)\otimes G(b_2)-\sum\limits_{b'_1,b'_2}a_{b'_1,b'_2}
G(b'_1)\otimes G(b'_2).$$
Here
the sum ranges over $(b'_1,b'_2)$ such that
\beqn
&H(b'_1\otimes b'_2)>0\,,\cr
&l(b_2)\le l(b'_1)<l(b_1)\,,\cr
&l(b_2)<l(b'_2)\le l(b_1)\,.
\endeqn
and $a_{b'_1,b'_2}\in \AZ$.
\end{corollary}

By Lemma \ref{indp}, $ C_{b_1,b_2}$ is uniquely determined.
Note that we shall see
$$a_{b'_1,b'_2}(0)=-\delta(b'_1\otimes b'_2=z^{H(b_1\otimes b_2)}b_1
\otimes z^{-H(b_1\otimes b_2)}b_2)$$ (see Lemma~\ref{lem:9}).

The following corollary is a consequence of
the corollary above and Lemma \ref{indp}.

\Lemma
$L(N)$ is a free $A$-module with
$\{ C_{b_1,b_2}\}_{H(b_1\otimes b_2)\le0}$ as its basis.
\end{lemma}

\begin{proposition}
\begin{enumerate}
\item 
The normally ordered wedges form a base of $\bw^2 V_{\aff}$.
\item 
$L(\bw^2V_\aff)$ is a free $A$-module
with the normally ordered wedges as a base.
\end{enumerate}
\end{proposition}

\proof
Lemma \ref{indp} implies the linear independence of
the normally ordered wedges
and Corollary \ref{cor:4} implies that they generate
$\bw^2V_\aff$.

(ii)~follows from~(i) and Corollary \ref{cor:4}.
\qed

\begin{corollary}
$N=\Ker\big(\tR-\psi(z\otimes 1,1\otimes z)\big)$.
\end{corollary}

\proof
We know already that
$N$ is contained in
$\Ker\big(\tR-\psi(z\otimes 1,1\otimes z)\big)$.
Since the normally ordered wedges are
linearly independent in
$V_\aff^{\otimes2}/\Ker\big(\tR-\psi(z\otimes 1,1\otimes z)\big)$
by Lemma \ref{indp},
$\bw^2V_\aff\to
V_\aff^{\otimes2}/\Ker\big(\tR-\psi(z\otimes 1,1\otimes z)\big)$
is injective.
\qed

\medskip
We define
for $n>0$
$$N_n=\sum_{k=0}^{n-2}\big(V_{\aff}^{\otimes k}\otimes N\otimes V_{\aff}
^{\otimes(n-k-2)}\big)\subset V_{\aff}^{\otimes n}$$
and then
$${\textstyle\bigwedge}^nV_{\aff}=V_{\aff}^{\otimes n}/N_n\,.$$
For $u_1,u_2,\dots,u_n\in V_\aff$, we denote by
$u_1\wedge u_2\wedge\cdots\wedge u_n$ the image of
$u_1\otimes u_2\otimes\cdots\otimes u_n$ in $\bw^nV_\aff$.

There is a $U_q(\Gg)$-linear homomorphism
$$\wedge:\bw^nV_\aff\otimes \bw^mV_\aff\to\bw^{n+m}V_\aff.$$
Let us set $L(V_{\aff}^{\otimes n})=L_{\aff}^{\otimes n}$
and let $L(\bigwedge^n V_\aff)$ be the image of
$L(V_{\aff}^{\otimes n})$ in $\bigwedge^n V_\aff$.
We call a sequence $(b_1,b_2,\dots,b_n)$ normally ordered if
its every consecutive pair is normally ordered,
i.e. if $H(b_j\otimes b_{j+1})>0$ for $j=1,\ldots,n-1$.
In this case we call $G(b_1)\wedge\cdots\wedge G(b_n)$
a {\it normally ordered wedge}.
Set
$$
L(N_n)=\sum_{k=0}^{n-2}L(V_\aff)^{\otimes k}\otimes_AL(N)\otimes_A
L(V_\aff)^{\otimes (n-2-k)}\subset L(V_\aff^{\otimes n}).$$

Note that
we have not yet seen
$L(N_n)\supset N_n\cap L(V_\aff^{\otimes n})$,
which will follow from Lemma \ref{lem:11}.
In the formulae below, we have to pay attention
to a difference
between modulo $qL(N_n)$ and modulo $q L(V_\aff^{\otimes n})$.

\Lemma\label{lem:9}
\begin{enumerate}
\item 
If $H(b_1\otimes b_2)=0$ then
$$G(z^ab_1)\wedge G(z^bb_2)\equiv
-G(z^bb_1)\wedge G(z^ab_2)\ \mod\ qL(\bw^2V_\aff).$$
\item 
If $H(b_1\otimes b_2)\le0$ then
\beqn
&& C_{b_1,b_{2}}\equiv b_1\otimes b_{2}+
\delta(H(b_1\otimes b_2)<0)
z^{H(b_1\otimes b_{2})}b_1\otimes z^{-H(b_1\otimes b_{2})}b_{2}\\
&&\phantom{----------------}\mod\,qL(V_\aff^{\otimes 2}).
\endeqn
\item 
If $H(b_j\otimes b_{j+1})=0$ for $j=1,\dots,n-1$, then
for any $\sigma\in S_n$,
\beqn
&& G(z^{a_1}b_1)\wedge G(z^{a_2}b_2)\wedge\cdots\wedge G(z^{a_n}b_n)\\
&&\phantom{-----}\equiv
\sgn(\sigma)G(z^{a_{\sigma(1)}}b_1)\wedge G(z^{a_{\sigma(2)}}b_2)\wedge\cdots
\wedge G(z^{a_{\sigma(n)}}b_n)\hfill\\
&&\phantom{--------------}\mod\ qL(\bw^nV_\aff).
\endeqn
\end{enumerate}
\enlemma
\proof
By Lemma \ref{lem:2}, (i) holds for $a=b=0$.
The general case is obtained by operating $z^a\otimes z^b+z^b\otimes z^a$
on $G(b_1)\otimes G(b_2)\equiv0$.

The other assertions follow from (i).
\qed

\Proposition\label{prop:10}
Let $a,c\in\BZ$ and $n\in\BZ_{>0}$.
Then for $b_1,\dots, b_n\in B_\aff$ with $a\le l(b_j)\le c$,
we have
\beqn
&G(b_1)\otimes\cdots\otimes G(b_n)
\in \sum \AZ G(b'_1)\otimes\cdots\otimes G(b'_n)+L(N_n)
\endeqn
where the sum ranges over normally ordered sequences
$(b'_1,\dots,b'_n)$ with
$a\le l(b'_j)\le c$ and $l(b'_1)\le l(b_1)$.
\enproposition

\proof
We shall prove this
by induction on $n$ and $l(b_1)$.
By the induction hypothesis on $n$, we may assume that $(b_2,\dots,b_n)$
is normally ordered.
If $H(b_1\otimes b_2)>0$, then we are done.
Assume that $H(b_1\otimes b_2)\le0$.
Then by Corollary \ref{cor:4},
we can write
\beqn
G(b_1)\otimes G(b_2)
&\equiv\sum\limits_{b'_1,b'_2} a_{b'_1,b'_2}G(b'_1)\otimes G(b'_2)\ \mod\ L(N)
\endeqn
with
$H(b'_1\otimes b'_2)>0$
and $l(b_2)\le l(b'_1)<l(b_1)$
and $l(b_2)<l(b'_2)\le l(b_1)$.
Then we have
\beqn
&&G(b_1)\otimes G(b_2)\otimes\cdots\otimes G(b_n)\\
&&\quad\equiv\sum a_{b'_1,b'_2}G(b'_1)\otimes G(b'_2)\otimes
G(b_3)\otimes\cdots\otimes G(b_n)
\ \mod\ L(N_n).
\endeqn
Since $a\le l(b_2)\le l(b'_1)<l(b_1)$, the induction proceeds.
\qed

This proposition says in particular that $\bw^nV_\aff$
is generated by the normally ordered wedges.
In order to see their linear independence,
we need the compatibility of the relations, which follow from the Yang--Baxter
equation for $R$.

\Lemma
Assume $H(b_1\otimes b_2)=H(b_2\otimes b_3)=0$.
Then for $a\ge b\ge c$, we have
\halign{\quad#&&$#$\cr
&(1+\delta_{a,b})&\tC_{z^ab_1,z^bb_2}\otimes G(z^cb_3)\hfill\cr
&&+(1+\delta_{b,c})\tC_{z^bb_1,z^cb_2}\otimes G(z^ab_3)+
(1+\delta_{a,c})\tC_{z^ab_1,z^cb_2}\otimes G(z^bb_3)\cr
&\equiv
(1+\delta_{b,c})&G(z^ab_1)\otimes\tC_{z^bb_2,z^cb_3}\hfill\cr
&&+(1+\delta_{a,c})G(z^bb_1)\otimes\tC_{z^ab_2,z^cb_3}+
(1+\delta_{a,b})G(z^cb_1)\otimes\tC_{z^ab_2,z^bb_3}\cr
\noalign{\vskip2pt}
&&\hfill\mod\ qL(N_3).\quad\cr}
\enlemma
\proof
We have the Yang-Baxter equation
$$\tR_{12}\circ \tR_{23}\circ\tR_{12}=
\tR_{23}\circ \tR_{12}\circ\tR_{23}.
$$
Here $\tR_{ij}$ is the action
of $\tR$ on the $i,j$-th components on $V_\aff^{\otimes 3}$.
Set $\psi_{21}=\psi(1\otimes z\otimes 1,z\otimes 1\otimes 1)$, etc.
Since $\tR+\psi(1\otimes z,z\otimes1)$ sends $L(V_\aff^{\otimes 2})$
to $L(N)$, $R_{ij}+\psi_{ji}$ sends $L(V_\aff^{\otimes 3})$
to $L(N_3)$.
Also we have
\beqn
&&\big(\tR+\psi(1\otimes z,z\otimes 1)\big)(G(z^ab_1)\otimes G(z^bb_2))\\
&&\phantom{---}\equiv G(z^bb_1)\otimes G(z^ab_2)+G(z^ab_1)\otimes
G(z^bb_2)\nn\\
&&\phantom{---}\equiv(1+\delta_{a,b}) C_{z^ab_1,z^bb_2}
\phantom{-----}\mod\,q\,L(V_\aff^{\otimes2}).
\endeqn
Since $L(N)=N\cap L(V_\aff^{\otimes2})$,
the above congruence is also true modulo $qL(N)$.
Since we have
\begin{align*}
  &\tR_{23}\circ\tR_{12}\big(G(z^ab_1)\otimes G(z^bb_2)\otimes G(z^cb_3)\big)
  \\
  &\phantom{---}\equiv G(z^bb_1)\otimes G(z^cb_2)\otimes
  G(z^ab_3)\qquad\qquad \mod\ qL(V_\aff^{\otimes 3}),
\end{align*}
etc., we have
\begin{align*}
&(\tR_{12}+\psi_{21})\circ\tR_{23}\circ\tR_{12}
\big(G(z^ab_1)\otimes G(z^bb_2)\otimes G(z^cb_3)\big)\cr
&\quad\equiv (\tR_{12}+\psi_{21})
\big( G(z^bb_1)\otimes G(z^cb_2)\otimes G(z^ab_3)\big)\\
&\quad\equiv (1+\delta_{b,c})\tC_{z^bb_1,z^cb_2}\otimes G(z^ab_3)
\phantom{--------} \mod\ qL(N_3),\\
\intertext{and similarly}
&(\tR_{23}+\psi_{23})\circ\tR_{12}
\big(G(z^ab_1)\otimes G(z^bb_2)\otimes G(z^cb_3)\big)\cr
&\quad\equiv (1+\delta_{a,c})G(z^bb_1)\otimes\tC_{z^ab_2,z^cb_3}
\phantom{--------}\mod\ qL(N_3).
\end{align*}
They imply
\beqn
&&\tR_{12}\circ \tR_{23}\circ\tR_{12}
\big(G(z^ab_1)\otimes G(z^bb_2)\otimes G(z^cb_3)\big)\cr
&&\equiv
(1+\delta_{b,c})\tC_{z^bb_1,z^cb_2}\otimes G(z^ab_3)\\
&&\phantom{------}-\psi_{21} \tR_{23}\circ\tR_{12}
\big(G(z^ab_1)\otimes G(z^bb_2)\otimes G(z^cb_3)\big)\cr
&&\equiv
(1+\delta_{b,c})\tC_{z^bb_1,z^cb_2}\otimes G(z^ab_3)
-(1+\delta_{a,c})\psi_{21}G(z^bb_1)\otimes\tC_{z^ab_2,z^cb_3}\cr
&&\phantom{------}+\psi_{21}\psi_{32}
\tR_{12}\big(G(z^ab_1)\otimes G(z^bb_2)\otimes G(z^cb_3)\big)
\cr
&&\equiv
(1+\delta_{b,c})\tC_{z^bb_1,z^cb_2}\otimes G(z^ab_3)
-(1+\delta_{a,c})G(z^bb_1)\otimes\tC_{z^ab_2,z^cb_3}\cr
&&\phantom{-----}+(1+\delta_{a,b})\tC_{z^ab_1,z^bb_2}\otimes G(z^cb_3)\\
&&\phantom{--------}
-\psi_{21}\psi_{32}\psi_{31}G(z^ab_1)\otimes G(z^bb_2)\otimes G(z^cb_3).
\endeqn
\smallskip
Here $\equiv$ is taken modulo $qL(N_3)$.
Similarly we have
\begin{align*}
  &\tR_{23}\circ \tR_{12}\circ\tR_{23}
  \big(G(z^ab_1)\otimes G(z^bb_2)\otimes G(z^cb_3)\big)\cr
  &\ \equiv
  (1+\delta_{a,b})G(z^cb_1)\otimes\tC_{z^ab_2,z^bb_3}
  -(1+\delta_{a,c})\tC_{z^ab_1,z^cb_2}\otimes G(z^bb_3)\cr
  &\,\phantom{-}+(1+\delta_{b,c})G(z^ab_1)\otimes\tC_{z^bb_2,z^cb_3}
  -\psi_{32}\psi_{31}\psi_{21}G(z^ab_1)\otimes G(z^bb_2)\otimes G(z^cb_3).
\end{align*}
Comparing these two identities, we obtain the desired result.
\qed

\Lemma\label{lem:11}
The $\BQ$-vector space
$L(N_n)/qL(N_n)$ is generated by
$G(b_1)\otimes\cdots\otimes G(b_{i-1})
\otimes  C_{b_i,b_{i+1}}
\otimes G(b_{i+2})\otimes\cdots\otimes G(b_{n})$
where $(b_1,\dots, b_n)$ ranges over the elements in $B_\aff^n$
such that $(b_{i+1},\dots, b_n)$ is normally ordered and
$H(b_i\otimes b_{i+1})\le0$.
\enlemma
\proof
$L(N_n)$ is generated by
$G(b_1)\otimes\cdots\otimes G(b_{i-1})
\otimes  C_{b_i,b_{i+1}}
\otimes G(b_{i+2})\otimes\cdots\otimes G(b_{n})$.
Here $H(b_i\otimes b_{i+1})\le0$ but
$(b_{i+1},\dots, b_n)$ is not necessarily normally ordered.
We shall prove that
such a vector can be written
as a $\BQ$-linear combination of
vectors satisfying the conditions as in the lemma,
by induction on
$n$ and descending induction on $i$.
Arguing by induction on $n$, we may assume $i=1$.
Write $b_k=z^{a_k}\tilde b_k$ with
$H(\tilde b_k\otimes\tilde b_{k+1})=0$.
Then $a_1\ge a_2$.
By Lemma \ref{lem:9} (iii), we may assume that
$a_3<a_4<\cdots<a_n$.
If $a_2<a_3$, there is nothing to prove.
Assume $a_2\ge a_3$. Then
the preceding lemma implies
\begin{align*}
  (&1+\delta_{a_1,a_2})
  \tC_{z^{a_1}\tilde b_1,z^{a_2}\tilde b_2}
  \otimes G(z^{a_3}\tilde b_3)\cr
  &\equiv\ -(1+\delta_{a_2,a_3})
  \tC_{z^{a_2}\tilde b_1,z^{a_3}\tilde b_2}\otimes G(z^{a_1}\tilde b_3)
  -(1+\delta_{{a_1},{a_3}})
  \tC_{z^{a_1}\tilde b_1,z^{a_3}\tilde b_2}\otimes G(z^{a_2}\tilde b_3)\cr
  &\quad+(1+\delta_{{a_2},{a_3}})G(z^{a_1}\tilde b_1)\otimes
  \tC_{z^{a_2}\tilde b_2,z^{a_3}\tilde b_3}
  +(1+\delta_{{a_1},{a_3}})
  G(z^{a_2}\tilde b_1)\otimes\tC_{z^{a_1}\tilde b_2,z^{a_3}\tilde b_3}\cr
  &\quad+(1+\delta_{{a_1},{a_2}})
  G(z^{a_3}\tilde b_1)\otimes\tC_{z^{a_1}\tilde b_2,z^{a_2}\tilde b_3}
  \phantom{-------}\mod\ qL(N_3).\quad
\end{align*}
Note that $a_3$ is the smallest among $(a_1,\dots,a_n)$.
After tensoring $G(z^{a_4}\tilde b_4)\otimes
\cdots\otimes G(z^{a_n}\tilde b_n)$,
the first two terms can be written in the desired form
by Lemma \ref{lem:9} (iii),
and the last three terms can be written in the desired form
by the hypothesis of induction on $i$.
\qed

\Theorem
The normally ordered wedges
form a base of $\bigwedge^nV_\aff$.
\entheorem
\proof
The normally ordered wedges generate $\bigwedge^nV_\aff$
by Proposition \ref{prop:10}.
We shall show that
any linear combination of normally ordered tensors
in $N_n$ vanishes.
Let $C$ be such a linear combination.
Since
$\bigcap_kq^kL(N_n)\subset\bigcap_kq^kL(\bigwedge^nV_\aff)=0$,
it is enough to show that $C\in L(N_n)$ implies $C\in qL(N_n)$.
By the preceding lemma, we can write
\begin{multline*}
  C\equiv
  \sum_{i=1}^{n-1}\sum_{(b_1,\dots,b_n)\in K_i}a_i(b_1,\dots,b_n)\
  G(b_1)\otimes\cdots\otimes G(b_{i-1})
  \otimes C_{b_i,b_{i+1}}\cr
  \phantom{---------}\otimes G(b_{i+2})\otimes\cdots\otimes G(b_{n})
  \quad\mod\ qL(N_n).
\end{multline*}
Here
the coefficents $a_i(b_1,\dots,b_n)$
belong to $\BQ$ and
$(b_{i+1},\dots,b_n)$ is normally ordered
for $(b_1,\dots,b_n)$$\in K_i$.
In order to show
the vanishing of $a_i(b_1,\dots,b_n)$,
let us calculate $C$ modulo
$\ qL(V_\aff^{\otimes n})$.
\beqn
&&C\equiv
\sum\limits_{i=1}^{n-1}\sum\limits_{(b_1,\dots,b_n)\in K_i}
a_i(b_1,\dots,b_n)\
b_1\otimes\cdots\otimes b_{i-1}\\
&&\phantom{------------}
\otimes
 C_{b_i,b_{i+1}}
\otimes b_{i+2}\otimes\cdots\otimes b_{n}
\phantom{-}\mod\,q L(V_\aff^{\otimes n}).
\endeqn
Since Lemma \ref{lem:9} (ii) implies
$$ C_{b_i,b_{i+1}}\equiv b_i\otimes b_{i+1}+
\delta(H(b_i\otimes b_{i+1})<0)
z^{H(b_i\otimes b_{i+1})}b_i\otimes z^{-H(b_i\otimes b_{i+1})}b_{i+1},$$
we have
\begin{align}
C&\equiv
\sum_{i=1}^{n-1}\sum_{(b_1,\dots,b_n)\in K_i}
a_i(b_1,\dots,b_n)\,
b_1\otimes\cdots\otimes b_{i-1}\nonumber\\
&\ \otimes
\big(b_i\otimes b_{i+1}+
\delta(H(b_i\otimes b_{i+1})<0)
z^{H(b_i\otimes b_{i+1})}b_i\otimes z^{-H(b_i\otimes b_{i+1})}b_{i+1}\big)\nn\\
&\phantom{-}\otimes b_{i+2}\otimes\cdots\otimes b_{n}
\phantom{-------}\mod\,q L(V_\aff^{\otimes n}).
\label{ll}
\end{align}
We shall show $a_i(b_1,\dots,b_n)=0$ by the decending induction on $i$.
Assume that $a_k(b_1,\dots,b_n)=0$ for $k>i$.
Note that $H(b_i\otimes b_{i+1})\le 0$, and
$H\big(z^{H(b_i\otimes b_{i+1})}b_i\otimes
z^{-H(b_i\otimes b_{i+1})}b_{i+1}\big)>0$ when $H(b_i\otimes b_{i+1})<0$.
We also note that
$(b_i,\dots,b_n)$ is not normally ordered
for $(b_1,\dots,b_n)\in K_i$
but it is normally ordered
for $(b_1,\dots,b_n)\in K_k$ with $k<i$.
By these observations, for $(b_1,\dots,b_n)\in K_i$, the coefficient
of $b_1\otimes b_2\otimes\cdots\otimes b_n$
in the right hand side of (\ref{ll})
is $a_i(b_1,\dots,b_n)$
and $b_1\otimes b_2\otimes\cdots\otimes b_n$
does not appear in $C$.
Hence $a_i(b_1,\dots,b_n)$ must vanish.
\qed

\Corollary\label{cor:13}
$L(\bw^nV_\aff)$ is a free $A$-module with the normally ordered
wedges as a base.
\encorollary

In fact, the normally ordered wedges
generate $L(\bw^nV_\aff)$ by Proposition \ref{prop:10}
and are linearly independent by the theorem above.


Let $B(\bigwedge^n V_\aff)$
be the set of normally ordered sequences.
Let us regard $B(\bigwedge^n V_\aff)$ as a subset of
$B_\aff^{\otimes n}$. Since it is invariant by $\te_i$ and
$\tf_i$,
we can endow $B(\bigwedge^n V_\aff)$ with the structure of crystal
induced by $B_\aff^{\otimes n}$.
We regard $B(\bigwedge^n V_\aff)$
as a basis of
$L(\bigwedge^n V_\aff)/qL(\bigwedge^n V_\aff)$.
Then we have

\Proposition
$\big(L(\bigwedge^n V_\aff),B(\bigwedge^n V_\aff)\big)$
is a crystal base of $\bigwedge^n V_\aff$.
\enproposition

The following lemma follows immediately from (\ref{symm}).

\Lemma\label{bos}
Let $f(z_1,\dots,z_n)$ be a symmetric Laurent polynomial.
Then $f(z\otimes1\otimes\cdots\otimes1,
1\otimes z\otimes1\otimes\cdots\otimes1,\dots,
1\otimes\cdots\otimes1\otimes z)$
induces an endomorphism of $\bw^nV_\aff$.
\enlemma


\newsection{Fock space}
\subsection{Ground state sequence}
In this section we shall introduce a $q$-deformed Fock space
in a similar way to the $A^{(1)}_n$--case(\cite{KMS}).

We continue the discussion
on the perfect crystal $B$ of level $l$.
Let us take a sequence $\{b^\circ_m\}_{m\in\BZ}$ in $B_{\aff}$ such that
\begin{align*}
  \lan c,\varepsilon(b^\circ_m)\ran&=l,\cr
  \varepsilon(b^\circ_m)&=\varphi(b^\circ_{m+1})\cr
  \text{and}\qquad
  H(b^\circ_m\otimes b^\circ_{m+1})&=1\,.
\end{align*}
We call $(\cdots,b^\circ_{-1},b^\circ_{0},b^\circ_{1},\dots)$
a ground state sequence.
If we give one of $b_m^\circ$,
then the other members of a ground state sequence
are uniquely determined.

Since $B$ is a finite set, there exist
a positive integer $N$ and an integer $c$ such that
\begin{align}
  b^\circ_{k+N}&=z^cb^\circ_k\qbox{for every $k$.}\label{peor}\\
\intertext{Take weights $\lam_m\in P$ of level $l$ satisfying}
  \lam_m&=\wt(b^\circ_m)+\lam_{m+1} \notag\\
  \text{and}\qquad
  \cl(\lam_m)&=\varphi(b^\circ_{m})=\varepsilon(b^\circ_{m-1}). \notag
\end{align}
Set $\vo_m=G(b^\circ_m)\in V_{\aff}$.

\subsection{Definition of Fock space}
For $m\in \BZ$, let us define first a (fake) $q$-deformed Fock space
$\tCF_m$
as the inductive limit ($k\to \infty$) of
$\bw^{k-m}V_{\aff}$,
where $\bw^{k-m}V_{\aff}\to\bw^{k+1-m}V_{\aff}$
is given by $u\mapsto u\wedge \vo_k$.
Intuitively
$\tCF_m$ is the subspace of $\bw^\infty V_{\aff}$
generated by the vectors of the form $u_m\wedge u_{m+1}\wedge\cdots$
with $u_k=\vo_{k}$ for $k\gge m$.
Similarly we define $L(\tCF_m)$ as the inductive limit of
$L(\bw^{k-m}V_{\aff})$.
We define the vacuum vector $\tvac{m}=\vo_m\wedge \vo_{m+1}\wedge\cdots\in
\tCF_m$.
Then any vector can be written as
$v\wedge\tvac{m+r}$ for some positive integer $r$
and $v\in\bw^rV_\aff$.
Note that $v\wedge\tvac{m+r}=0$ if and only if
$v\wedge \vo_{m+r}\wedge\cdots\wedge \vo_{m+s}=0$
for some $s>r$.

Then we introduce the true ($q$-deformed) Fock space by
$$\CF_m=\tCF_m/\big(\bigcap_{n>0}q^nL(\tCF_m)\big)\,.$$
Let $L(\CF_m)\subset \CF_m$ be the image of $L(\tCF_m)$,
and $\vac{m}$ the image of $\tvac{m}$.
\hb
We have the homomorphism
$$\wedge\,:\, \bw^rV_\aff\otimes \CF_{m+r}\to \CF_{m}.$$

For a normally ordered sequence $(b_{m},b_{m+1},\dots)$
in $B_\aff$ such that $b_k=b^\circ_{k}$ for $k\gge m$,
we call $G(b_m)\wedge G(b_{m+1})\wedge\cdots\in \CF_m$
a {\it normally ordered wedge}.

\Theorem\label{main}
The normally ordered wedges  form a base of $\CF_m$.
\entheorem

In order to prove this theorem, we need some preparations.

\Lemma\label{crt}
If $l(b)>l(b_m^\circ)$, then $H(b\otimes b_{m+1}^\circ)\le 0$.
\enlemma

\proof
If $l(b)\gge0$, then the assertion holds.
Let us prove it by descending induction on $l(b)$.
Assume that there is $i\in I$ such that
$\te_i(b\otimes b_{m+1}^\circ)=(\te_ib)\otimes b_{m+1}^\circ
\not=0$.
Then $l(b)<l(\te_ib)$
and hence
$H(b\otimes b_{m+1}^\circ)=H(\te_ib\otimes b_{m+1}^\circ)\le 0$
by the hypothesis of induction.
Hence we may assume that there is no such $i$.
Then $\varepsilon_i(b)\le\varphi_i(b_{m+1}^\circ)$ for any $i$,
and hence $b=z^ab_m^\circ$ for some $a\in\BZ$.
Since $l(b)>l(b_m^\circ)$, we have $a>0$.
Therefore $H(b\otimes b_{m+1}^\circ)=1-a\le0$.
\qed

\Proposition\label{vv}
Assume $H(b\otimes b^\circ_m)\le 0$.
Then for every n we can find $m_1\ge m$ such that
$$G(b)\wedge \vo_m\wedge\cdots\wedge \vo_{m_1}\in q^nL(\bw^{m_1-m+2}
V_\aff).$$
\enproposition

\proof
We shall prove this by induction on $n$ and $H(b\otimes b_m^\circ)$.
\hb
Set $H(b\otimes b_m^\circ)=-c$ and
$$G(b)\wedge \vo_m
=\sum a(b_1,b_2)G(b_1)\wedge G(b_2).$$
Here the sum ranges over normally ordered pairs $(b_1,b_2)$
such that
\beq
&l(b_m^\circ)\le l(b_1)<l(b),\nn\\
&l(b_m^\circ)<l(b_2)\le l(b).\label{eq:20}
\endeq
By the preceding lemma $H(b_2\otimes b_{m+1}^\circ)\le0$.
Lemma~\ref{lem:9} (i) implies
$$a(b_1,b_2)\equiv
-\delta\big(c<0 \hbox{ and }
(b_1,b_2)=(z^{-c}b,z^{c}b_m^\circ)\big)\ \mod\ qA.$$
We have
\beqn
G(b)\wedge \vo_m\wedge\cdots\wedge \vo_{m_1}
&=\sum a(b_1,b_2)G(b_1)\wedge G(b_2)\wedge \vo_{m+1}\wedge\cdots\wedge
\vo_{m_1}.
\endeqn
Since $l(b_2)>l(b_m^\circ)$, we have
$G(b_2)\wedge \vo_{m+1}\wedge\cdots\wedge \vo_{m_1}\in q^{n-1}L(\bw V_\aff)$.
Hence
$a(b_1,b_2)G(b_1)\wedge G(b_2)\wedge \vo_{m+1}\wedge\cdots\wedge \vo_{m_1}$
belongs to $q^nL(\bw V_\aff)$
except $c<0$ and $(b_1,b_2)=(z^{-c}b,z^{c}b_m^\circ)$.

Assume that $c<0$ and $(b_1,b_2)=(z^{-c}b,z^{c}b_m^\circ)$.  Then we have
$0\ge H(z^{c}b_m^\circ\otimes b_{m+1}^\circ)=1-c>H(b\otimes b^\circ_m)$.
Hence $a(b_1,b_2)G(b_1)\wedge G(b_2)\wedge \vo_{m+1}\wedge\cdots\wedge
\vo_{m_1}$ belongs to $q^nL(\bw^{m_1-m+2} V_\aff)$ by the hypothesis of
induction on $H(b\otimes b_{m}^\circ)$.  \qed

\Remark Assume that $c$ in~(\ref{peor}) is positive (or equivalently,
$l(b_m^\circ)$ tends to infinity as $m$ tends to infinity).  Then
$H(b\otimes b_m)\le 0$ implies $G(b)\wedge \vo_m\wedge\cdots\wedge
\vo_{m_1}=0$ for $m_1\gge m$.  In fact by the same argument as above
we have $G(b)\wedge \vo_m\wedge\cdots\wedge \vo_{m_1}\in
\sum_{b'}\bw^{m_1-m+1}V_\aff\wedge G(b')$ where $b'$ satisfies
$l(b^\circ_{m_1})<l(b')\le l(b)$.

Note that, under the condition of the proposition,
$G(b)\wedge \vo_m\wedge \vo_{m+1}\wedge\cdots\wedge \vo_k=0$ for $k\gge m$
is false in general.

\medskip
A similar argument shows the following dual statement.

\Proposition\label{dual}
Assume $H(b^\circ_m\otimes b)\le 0$.
Then for every n we can find $m_1\le m$ such that
$$\vo_{m_1}\wedge\cdots\wedge \vo_{m}\wedge
G(b)\in q^nL(\bw^{m-m_1+2} V_\aff).$$
\enproposition

As an immediate consequence of Proposition~\ref{vv},
we obtain the following result.

\Theorem\label{kern}
For any vector $b\in B_\aff$ such that $H(b\otimes b^\circ_{m})\le0$, we have
the equality in $\F_m$
$$G(b)\wedge\vac{m}=0\,.$$
\entheorem

\noindent
{\sl Proof of Theorem}~\ref{main}.
Any vector in $\CF_m$ can be written in the form $v\wedge\vac{m+r}$
with $v\in\bw^rV_\aff$.  We may assume that $v$ is a normally ordered
wedge $G(b_m)\wedge\cdots\wedge G(b_{m+r-1})$.  If $H(b_{m+r-1}\otimes
b^\circ_{m+r})>0$, then $v\wedge\vac{m+r}$ is a normally ordered wedge
and otherwise $v\wedge\vac{m+r}=0$ by Proposition~\ref{vv}.

The linear independence follows immediately from the corresponding
statement for the wedge space (Corollary~\ref{cor:13}).
\qed
\medskip

By a similar argument, we have
\Proposition\label{free}
$L(\CF_m)$ is a free $A$-submodule of $\CF_m$ generated by the
normally ordered wedges.
\enproposition

\Proposition\label{kern2}
\beqn
\bigcap\limits_{n>0}q^nL(\tCF_m)&=&
\sum\limits_{H(b\otimes b_{m+r}^\circ)\le0}
\bw^{r-1}V_\aff\wedge G(b)\wedge\tvac{m+r}\\
&=&\sum\limits_{l(b)>l( b_{m+r-1}^\circ)}
\bw^{r-1}V_\aff\wedge G(b)\wedge\tvac{m+r}.
\endeqn
\enproposition
\proof
The first equality follows from Theorems~\ref{main} and~\ref{kern}
and the last follows from Lemma~\ref{crt} and~(\ref{eq:20}).
\qed


As a corollary of Theorem \ref{kern} we have the following result concerning
vertex operators.

\Proposition\label{vertex}
Let $V(\lam_m)$ be the irreducible $U_q(\Gg)$-module
with highest weight $\lam_m$
and $u_{\lam_m}$ its highest wei\-ght vector.
Let $\Phi:V_\aff\otimes V(\lam_m)\to V(\lam_{m-1})$
be an intertwiner.
Then for any vector $b\in B_\aff$ such that $H(b\otimes b^\circ_{m})\le0$,
$\Phi(G(b)\otimes u_{\lam_m})=0$.
\enproposition

\proof
As proved in~\cite{DJO}, the intertwiner is unique
up to a constant.
As seen in the next two subsections,
$\CF_m$ has a $U_q(\Gg)$-module structure and
contains $V(\lam_m)$ as a direct summand.
By this embedding, the highest vector
$u_{\lam_m}$ of $V(\lam_m)$
corresponds to
$\vac{m}$.
Therefore
$\Phi$ is given as the composition:
$$V_\aff\otimes V(\lam_m)\to V_\aff\otimes\CF_m
\to\CF_{m-1}\to V(\lam_{m-1}).$$
Now the result follows from Theorem~\ref{kern}.
\qed

\Remark
It is known (see e.g.~\cite{DJO})
that
$\Phi(v\otimes u_{\lam_m})=0$
for $v\in (V_\aff)_{\lam_{m-1}-\lam_m}$
such that $v\in \sum_ie_i^{1+\lan h_i,\lam_{m-1}\ran}V_\aff$.
On the other hand, by the property of the lower global base
(\cite{K2}),
$G(b)$ belongs to $\sum_ie_i^{1+\lan h_i,\lam_{m-1}\ran}V_\aff$
if and only if
$\varphi_i(b)>\lan h_i,\lam_{m-1}\ran$ for some $i$.
Therefore,
$\Phi(G(b)\otimes u_{\lam_m})=0$
for $b\in (B_\aff)_{\lam_{m-1}-\lam_m}$
other than $b_{m-1}^\circ$.

This observation shows that we have to take a lower global base in
order to have Theorem~\ref{kern}.  Theorem~\ref{kern}, as well as
Proposition~\ref{vertex}, does not hold for an arbitrary choice of
base other than the lower global base.  In the course of our
construction of the Fock space, we have not used explicitly the
property of the lower global base.  This is hidden in postulate~(R).
This postulate fails for an arbitrary choice of base.

\subsection{$U_q(\Gg)$-module structure on the Fock space}
Let us define the action of $U_q(\Gg)$ on $\CF_m$.  We define first
the action of the Cartan part of $U_q(\Gg)$ by assigning weights.  We
set $\wt(\vac{m})=\lam_m$ and
$\wt(v\wedge\vac{m+r})=\wt(v)+\wt(\vac{m+r})$ for $v\in \bw^rV_\aff$.
This defines the weight decomposition of the Fock space.


Let $B(\CF_m)$ denote the set of normally ordered sequences
$(b_m,b_{m+1},\dots)$ in $B_\aff$ such that
$b_k=b_k^\circ$ for $k\gge m$.
Then it has a crystal structure as in~\cite{KMN1}.
Moreover $B(\CF_m)$ may be considered as
a base of $L(\CF_m)/qL(\CF_m)$ by Proposition~\ref{free}.
We write $b_m\wedge b_{m+1}\wedge\cdots$ for $(b_m,b_{m+1},\dots)$.

\Proposition\label{char}
\begin{enumerate}
\item 
$\ch(\CF_m)=\ch(V(\lam_m))\prod_{k>0}(1-{\mathrm e}^{-k\delta})^{-1}
$.
\item 
The weights of $\CF_m$
appear as weights of $V(\lam_m)$.
In particular, any weight $\mu$ of $\CF_m$
satisfies $s(\mu)\le s(\lam_m)$
(see the end of~\S\ref{sec:ener} for
$s:P\to\BQ$). Moreover, $s(\mu)=s(\lam_m)$ implies $\mu=\lam_m$.
\item 
For any $\mu\in P$, $\dim (\CF_m)_\mu<\infty$.
\item 
$(\CF_m)_{\lam_m-n\alpha_i}=
\begin{cases}
  KG(\tf_i^nb_m^\circ)\wedge\vac{m+1}
  & \text{if $0\le n\le\lan h_i,\lam_m\ran$},\cr
  0&\text{otherwise}.\cr
\end{cases}$
\item 
If $b\in B_\aff$ satisfies
$\wt(b)=\wt(b_m^\circ)-n\alpha_i$, then
$G(b)\wedge\vac{m+1}=0$ unless $0\le n\le\lan h_i,\lam_m\ran$
and $b=\tf_i^nb_m^\circ$.
\item 
Any highest weight element of
$B(\CF_m)$ has the form
$z^{a_m}b^\circ_m\wedge z^{a_{m+1}}b^\circ_{m+1}\wedge\cdots$
with $a_m\le a_{m+1}\le\cdots$ and $a_k=0$ for $k\gge m$.
\item 
For $b_m\wedge b_{m+1}\wedge\cdots\in B(\CF_m)$,
$b_m=b^\circ_m$ implies
$b_k=b_k^\circ$ for any $k\ge m$.
\end{enumerate}
\enproposition

\proof
By Proposition~4.6.4 in~\cite{KMN1} (see also Appendix~\ref{perf}),
we have
$$\ch(V(\lam_m))={\mathrm e}^{\lam_m}\sum {\mathrm e}^
{\sum_{n\ge m}(\wt(b_n)-\wt(b_n^\circ))}$$
where the sum ranges over the family $\CB_0$ of
sequences $b_m,b_{m+1},\dots$ in $B_\aff$
such that $b_n=b_n^\circ$ for $n\gge m$ and $H(b_n\otimes b_{n+1})=1$ for
any $n\ge m$.
On the other hand, we have
$$\ch(\CF_m)={\mathrm e}^{\lam_m}\sum {\mathrm e}^
{\sum_{n\ge m}(\wt(b_n)-\wt(b_n^\circ))}$$
where the sum ranges over the family $\CB$ of normally ordered
$b_m,b_{m+1},\dots$
such that $b_n=b_n^\circ$ for $n\gge m$.
We have
\beqn
&&\CB=\{(z^{-a_m}b_m,z^{-a_{m+1}}b_{m+1},\dots)\,;\\
&&\phantom{---}
(b_m,b_{m+1},\dots)\in \CB_0,\ a_m\ge a_{m+1}\ge\cdots\hbox{ and
$a_n=0$ for $n\gge m$}\}.
\endeqn
To obtain~(i), it is enough to remark that $z$ has weight $\delta$.

The assertions (ii)--(vi) follow from~(i) and Theorem~\ref{kern}.
The assertion (vii) follows from (vi) and
$$\tf_i(z^{a_m}b^\circ_m\wedge z^{a_{m+1}}b^\circ_{m+1}\wedge\cdots)
=z^{a_m}\tf_ib^\circ_m\wedge z^{a_{m+1}}b^\circ_{m+1}\wedge\cdots.$$
\qed

Now we shall define the action of $e_i$ and $f_i$ on $\CF_m$.

Taking $\{q^nL(\CF_m)\}_n$ as a neighborhood system of $0$,
$\CF_m$ is endowed with a so called
$q$-adic topology.
Since $\bigcap_nq^nL(\CF_m)=0$ by construction,
the $q$-adic topology is separated.
Since we use $K=\BQ(q)$ as a base field,
$\CF_m$ is not complete with respect to this topology.
For any $\mu\in P$, the completion of $(\CF_m)_\mu$
is $\BQ((q))\otimes_K(\CF_m)_\mu$.

\Proposition
For any vectors $u_m,u_{m+1},\dots\in V_\aff$
such that $u_k=\vo_k$ for $k\gge m$,
\beq
&&\sum_{k\ge m}t_i^{-1}(u_m\wedge\cdots\wedge u_{k-1})\wedge
e_iu_k\wedge u_{k+1}\wedge\cdots\label{conve}\\
\noalign{\hbox{and}}
&&\sum_{k\ge m}u_m\wedge\cdots\wedge u_{k-1}\wedge
f_iu_k\wedge t_i(u_{k+1}\wedge\cdots)\label{convf}
\endeq
converge in the $q$-adic topology
to elements of $\BQ((q))\otimes_K\CF_m$.
\enproposition

\proof
First note that
$(e_i\vo_k)\wedge\vac{k+1}=0$ because $\lam_k+\alpha_i$
is not a weight of $\CF_k$.
Hence, only finitely many terms survive in~(\ref{conve}).

In order to prove the convergence of~(\ref{convf}),
we may assume that
$u_k=\vo_k$ for every $k\ge m$.
Then $$\vo_m\wedge\cdots\wedge \vo_{k-1}\wedge
f_i\vo_k\wedge t_i(\vo_{k+1}\wedge\cdots)
=q_i^{\lan h_i,\lam_{k+1}\ran}\vo_m\wedge\cdots\wedge \vo_{k-1}\wedge
f_i\vo_k\wedge\vac{k+1}.$$
Since $\lan h_i,\lam_{k+1}\ran$ takes only finitely many values,
it is enough to show that
$\vo_m\wedge\cdots\wedge \vo_{k-1}\wedge
f_i\vo_k\wedge\vac{k+1}$
converges in the $q$-adic topology.
This follows from the following lemma.
\qed

\Lemma\label{qconv}
Let $C$ be an endomorphism of the $K$-vector
space $V_\aff$ of weight $\mu\neq0$.
Assume that $Cz=zC$.
Then for any $m$,
$\vo_m\wedge\cdots\wedge \vo_{k-1}\wedge C\vo_k\wedge\vac{k+1}$
converges to $0$ in the $q$-adic topology when
$k$ tends to infinity.
\enlemma

\proof
Write
$$C\vo_k=\sum_{\nu}c_{k,\nu}G(b_{k,\nu}).$$
Take $N$ and $c$ as in~(\ref{peor}).
Then we have also the periodicity
$b_{k+N,\nu}=z^cb_{k,\nu}$ and $c_{k+N,\nu}=c_{k,\nu}$.
Hence $c_{k,\nu}$ is bounded with respect to the $q$-adic topology.
Therefore it is enough to show that
$\vo_m\wedge\cdots\wedge \vo_{k-1}\wedge G(b_{k,\nu})\wedge\vac{k+1}$
converges to $0$.
By Proposition~\ref{char}\,(vi),
$(b_m^\circ,\dots,b_{k-1}^\circ,
b_{k,\nu},b_{k+1}^\circ,\dots)$
is not normally ordered. It means that
either $H(b_{k,\nu}\otimes b_{k+1}^\circ)\le 0$
or $H(b_{k-1}^\circ\otimes b_{k,\nu})\le 0$.
If $H(b_{k,\nu}\otimes b_{k+1}^\circ)\le 0$ then
$G(b_{k,\nu})\wedge\vac{k+1}$ vanishes.
If $H(b_{k-1}^\circ\otimes b_{k,\nu})\le 0$,
then $\vo_{m-s}\wedge\cdots \vo_{k-1}\wedge G(b_{k,\nu})\wedge\vac{k+1}$
converges to $0$ when $s$ tends to infinity by Proposition~\ref{dual}.
By shifting the indices,
$\vo_m\wedge\cdots\wedge \vo_{k-1}\wedge G(b_{k,\nu})\wedge\vac{k+1}$
converges to $0$.
\qed

Let us set
\beq
f_i\tvac{m}&=&f_i\vo_m
\wedge t_i\vac{m+1}+\vo_m\wedge f_i\vo_{m+1}\wedge t_i\vac{m+2}
+\cdots
\endeq
Then it is an element of $\Fh_m$.

\Lemma\label{deff}
$f_i\tvac{m}$ belongs to $\CF_m$.
\enlemma

\proof
Let us take $c$ and $N$ as in~(\ref{peor}).
We define the isomorphism
$\psi_m:\CF_m\to \CF_{m+N}$
by $u_m\wedge u_{m+1}\wedge\cdots
\mapsto z^cu_m\wedge z^cu_{m+1}\wedge\cdots$.
Then $f_i\tvac{m}$ satisfies the recurrence relation
\beqn
&&f_i\tvac{m}-\vo_m\wedge \vo_{m+1}\wedge\cdots\wedge \vo_{m+N-1}\wedge
\psi_m(f_i\tvac{m})\\
&&=f_i(\vo_m\wedge \vo_{m+1}\wedge\cdots\wedge \vo_{m+N-1})\wedge t_i\vac{m+N}
\in\CF_m.
\endeqn
Hence the result follows from the following lemma.
\qed

\Lemma
For $\mu\in P\backslash \{\lam_m\}$,
the endomorphism of
$(\CF_m)_\mu$
given by
$w\mapsto w-
\vo_m\wedge \vo_{m+1}\wedge\cdots\wedge \vo_{m+N-1}\wedge\psi_m(w)$
is an isomorphism.
\enlemma
\proof
It is enough to show its injectivity.
We show that
$w=\vo_m\wedge \vo_{m+1}\wedge\cdots\wedge \vo_{m+N-1}\wedge\psi_m(w)$
implies $w=0$.

For $b_m\wedge b_{m+1}\wedge\cdots\in B(\CF_m)_\mu$,
$(b^\circ_{m-1},b_m,b_{m+1},\dots)$
is not normally ordered by
Proposition~\ref{char} (vii), and
hence $H(b^\circ_{m-1}\otimes b_m)\le0$.
Proposition~\ref{dual}
implies that
$\vo_{m-kN}\wedge\cdots\wedge \vo_{m-1}\wedge G(b_m)\wedge
G(b_m+1)\wedge\cdots$
belongs to $qL(\CF_{m-kN})$ for $k\gge0$.
Shifting the indices,
we conclude that
\begin{displaymath}
  \vo_{m}\wedge\cdots\wedge \vo_{m+kN-1}\wedge
\psi_{m+(k-1)N}\cdots\psi_{m+N}\psi_m
\bigl(G(b_m)\wedge G(b_{m+1})\wedge\cdots\bigr)
\end{displaymath}
belongs to $qL(\CF_{m})$ for $k\gge0$.
Therefore
the homomorphism
$$C:w\mapsto \vo_{m}\wedge\cdots\wedge \vo_{m+kN-1}\wedge
\psi_{m+(k-1)N}\cdots\psi_{m+N}\psi_m(w)$$
sends $L(\CF_m)_\mu$ to
$qL(\CF_{m})_\mu$ for $k\gge0$.
This shows the injectivity of the endomorphism
$\id_{(\CF_{m})_\mu} -C$.
\qed

Now we define
\beq
e_i(v\wedge\tvac{m+r})&=&e_iv\wedge\vac{m+r},\nn\\
f_i(v\wedge\tvac{m+r})&=&f_iv\wedge t_i\vac{m+r}
+v\wedge f_i\tvac{m+r}
\endeq
for $v\in\bw^rV_\aff$.
Then $e_i$ and $f_i$
are well-defined homomorphisms from
$\tCF_m$ to $\CF_m$.
They satisfy
\beq
e_i(v\wedge u)&=&e_iv\wedge u+t_i^{-1}v\wedge e_iu\nn\\
f_i(v\wedge u)&=&f_iv\wedge t_iu+v\wedge f_iu\label{eq:17}
\endeq
for $v\in\bw^rV_\aff$ and $u\in\tCF_{m+r}$.
In order to see that they define endomorphisms of $\CF_m$,
we need to show the following proposition
(see Proposition~\ref{kern2}).

\Proposition\label{comp}
Assume that $b\in B_\aff$ satisfies $l(b)>l(b_{m}^\circ)$.
Then we have the equalities in $\CF_m$.
\beq
&&e_i\big(G(b)\wedge\tvac{m+1}\big)=0,\label{vane}\\
&&f_i\big(G(b)\wedge\tvac{m+1}\big)=0.\label{vanf}
\endeq
\enproposition

The first equality
(\ref{vane}) follows from
the fact that $\wt(b)+\alpha_i+\lam_{m+1}$
is not a weight of $\CF_m$
(following Proposition~\ref{char}~(ii)).

Let us prove~(\ref{vanf}).
First note that
the same consideration on the weight implies that
\beq
&\hbox{if $l(b)>l(b_{m}^\circ)$
and $\wt(b)\not=\wt(b_{m}^\circ)+\alpha_i$, then
(\ref{vanf}) holds.}\phantom{--}&\label{naka}
\endeq

Hence in order to prove~(\ref{vanf}),
we  may assume that
\beq
\wt(b)=\wt(b_{m}^\circ)+\alpha_i.\label{cond;2}
\endeq

\Sublemma
Under the condition~(\ref{cond;2}), we can write
\beq
G(b)\wedge \vo_{m+1}\wedge\cdots\wedge \vo_{m+r}
=u+\sum a(b_{0},\dots,b_{r}) G(b_{0})\wedge\cdots\wedge G(b_{r}).
\phantom{--}&&\label{eq:13}
\endeq
Here $u$ satisfies
$f_i(u\wedge\tvac{m+r+1})=0$,
the coefficients
$a(b_{0},\dots,b_{r})$ belong to $q^rA$,
and the sum ranges over
$(b_{0},\dots,b_{r})$ such that
\beq
&&\wt(b_j)=
\begin{cases}
  \wt(b_{m+j}^\circ) &\text{for } 0\le j<r\cr
  \wt(b_{m+r}^\circ)+\alpha_i&\text{for }j=r
\end{cases}.
\endeq
\ensublemma

\proof
We shall prove this by induction on $r$.
Assuming~(\ref{eq:13}) for $r$,
let us prove~(\ref{eq:13}) for $r+1$.
Since $H(b_r\otimes b_{m+r+1}^\circ)\le0$ by Lemma~\ref{crt},
we can write
\beq
G(b_r)\wedge \vo_{m+r+1}
=\sum a_{b',b''}G(b')\wedge G(b'').
\endeq
Here $(b',b'')$ ranges over normally ordered pairs
such that
\beq
l(b_{m+r+1}^\circ)\le l(b')<l(b_r),\nn\\
l(b_{m+r+1}^\circ)<l(b'')\le l(b_r).\label{eq:10}
\endeq
If $\wt(b'')\neq\wt(b_{m+r+1}^\circ)+\alpha_i$,
then we have $G(b'')\wedge\vac{m+r+2}=0$ and
$f_i(G(b'')\wedge\tvac{m+r+2})=0$
by~\eqref{eq:10} and~\eqref{naka}.
Therefore $f_i(G(b_{0})\wedge\cdots\wedge G(b_{r-1})
\wedge G(b')\wedge G(b'')\wedge\tvac{m+r+2})=0$.
If $\wt(b'')=\wt(b_{m+r+1}^\circ)+\alpha_i$,
then
$\wt(b')=\wt(b_{m+r}^\circ)$.
Moreover Lemma~\ref{lem:9}~(i) implies
$a_{b',b''}\in qA$.
Thus the induction proceeds.
\qed

We resume the proof of Proposition~\ref{comp}.
We have
\beqn
f_i(G(b)\wedge\tvac{m+1})&=&
\sum a(b_{0},\dots,b_{r})\\
&&\big(
G(b_{0})\wedge\cdots\wedge G(b_{r})\wedge f_i\tvac{m+r+1}\\
&&\phantom{-}+
\sum\limits_{0\le j\le r}
G(b_{0})\wedge\cdots\wedge
G(b_{j-1})\wedge
f_iG(b_j)\wedge\\
&&\phantom{----}t_iG(b_{j+1})\wedge\cdots\wedge t_i G(b_{r})
\wedge t_i\vac{m+r+1}\big).
\endeqn
There is a constant $s$ such that
$f_iL_\aff\in q^sL_\aff$, and
$f_i\tvac{m+r+1}$ is bounded with respect to the $q$-adic topology.
Moreover,
$t_iG(b_{j+1})\wedge\cdots\wedge t_i G(b_{r})
\wedge t_i\vac{m+r+1}
=q_i^{\lan h_i,\lam_{m+j+1}\ran+2\delta(j<r)}
G(b_{j+1})\wedge\cdots\wedge G(b_{r})
\wedge\vac{m+r+1}$ and
$\lan h_i,\lam_{m+j}\ran$ is bounded from below.
Hence there is a constant $d$ independent of $r$
such that
$$f_i\big(G(b)\wedge\tvac{m+1}\big)\in q^{r+d}L(\CF_m)
\qbox{for every $r$.}$$
This implies the equality~(\ref{vanf}) in $\CF_m$.
This completes the proof of Proposition~\ref{comp}.

\bigskip
Thus we have defined the action of $e_i$
and $f_i$ on $\CF_m$.
Now we shall show the commutation relations between them.

\Proposition\label{prop:12} On $\CF_m$ we have
\beq
[e_i,f_j]=\delta_{ij}(t_i-t_i^{-1})/(q_i-q_i^{-1}).
\label{eq:4.18}
\endeq
\enproposition

\proof
First note that
(\ref{eq:17}) implies
\beq
[e_i,f_j](v\wedge u)&=&
[e_i,f_j]v\wedge t_ju+t_i^{-1}v\wedge[e_i,f_j] u.\label{co}
\endeq
for $v\in\bw^rV_\aff$ and $u\in\CF_{m+r}$.
Hence, it is enough to prove that
the equality~(\ref{eq:4.18}) holds when it is applied to the vacuum vector.
If $i\not=j$ then $[e_i,f_j]\vac{m}=0$ because
$\lam_m+\alpha_i-\alpha_j$ is not a weight of $\CF_m$.

Now we shall show $[e_i,f_i]=\{t_i\}_i$.
Here $\{x\}_i=(x-x^{-1})/(q_i-q_i^{-1})$.

Since  $[e_i,f_i]\vac{m}$ has weight $\lam_m$,
there is $c_m\in K$ such that
$[e_i,f_i]\vac{m}=c_m\vac{m}$.
Then by~(\ref{co}), we have
\beqn
[e_i,f_i]\vac{m}&=&
[e_i,f_i]\vo_m\wedge t_i\vac{m+1}+t_i^{-1}\vo_m\wedge[e_i,f_i]\vac{m+1}\nn\\
&=&\big(q_i^{\lan h_i,\lam_{m+1}\ran}[\lan h_i,\lam_{m}-\lam_{m+1}\ran]_i
+q_i^{\lan h_i,\lam_{m+1}-\lam_{m}\ran}c_{m+1}\big)\vac{m}.
\endeqn
Hence we have a recurrence relation
\beqn
c_m&=&q_i^{\lan h_i,\lam_{m+1}\ran}[\lan h_i,\lam_{m}-\lam_{m+1}\ran]_i
+q_i^{\lan h_i,\lam_{m+1}-\lam_{m}\ran}c_{m+1}.
\endeqn
Solving this, there is a constant
$a\in K$ such that
\beq
c_m&=&[\lan h_i,\lam_{m}\ran]_i+q_i^{-\lan h_i,\lam_{m}\ran}a
\qbox{for every $m$.}
\endeq
Namely we have
$[e_i,f_i](\vac{m})=(\{t_i\}_i+at_i^{-1})\vac{m}$.
Hence for $v\in \bw^{r}V_\aff$
\beqn
[e_i,f_i](v\wedge\vac{m+r})&=&
[e_i,f_i]v\wedge t_i\vac{m+r}+t_i^{-1}v\wedge[e_i,f_i]\vac{m+r}\\
&=&\{t_i\}_iv\wedge t_i\vac{m+r}
+t_i^{-1}v\wedge(\{t_i\}_i+at_i^{-1})\vac{m+r}\\
&=&(\{t_i\}_i+at_i^{-1})(v\wedge\vac{m+r}).
\endeqn
Thus we obtain
\beq
[e_i,f_i]&=&\{t_i\}_i+at_i^{-1}.\label{co1}
\endeq
Let us show the vanishing of $a$.

By induction on $n$ we can see the following commutaion relation
\beq
e_i^{(n)}f_i^{(n)}
&=&\sum_{k=0}^{n}f_i^{(n-k)}e_i^{(n-k)}
{{\prod_{\nu=0}^{k-1}(\{q_i^{-\nu}t_i\}_i+aq_i^{\nu}t_i^{-1})}\over{[k]_i!}}.
\endeq
Setting $c=\lan h_i,\lam_{m}\ran$, we have
$f_i^{(c+1)}\vac{m}=0$ by Proposition~\ref{char}~(iv).
Hence
\beqn
0&=&e_i^{(c+1)}f_i^{(c+1)}\vac{m}\\
&=&{{\prod_{\nu=0}^{c+1}([c-\nu]_i+aq_i^{\nu-c})}\over{[c+1]_i!}}
\vac{m}.
\endeqn
Therefore there is an integer $s$ such that
$a=-q_i^s[s]_i$.
Then the commutation relation~(\ref{co1}) can be rewritten as
$$[e_i,q_i^{-s}f_i]=\{q_i^{-s}t_i\}_i.$$
Hence
$e_i$, $q_i^{-s}f_i$ and $q_i^{-s}t_i$ form
$U_q(\Gsl_2)$. Then the representation theory
of $U_q(\Gsl_2)$ and Proposition~\ref{char}~(iv) implies $s=0$.
In fact, the string containing the weight of $\vac{m}$
(with respect to $q_i^{-s}t_i$)
is $\{c-s-2n;0\le n\le c\}$, and
hence the symmetry of a string under the simple reflection
implies $c-s=-(-c-s)$.
\qed

Thus the actions of $e_i$ and $f_i$ satisfy the commutation relations.
By Proposition~\ref{char}~(ii),
for any $i\in I$ and $\mu\in P$, $\mu+n\alpha_i$ is a weight of $\CF_m$
only for a finitely many integers $n$.
Therefore $\CF_m$ is integrable over the
$U_q(\Gsl_2)_i=\lan e_i,f_i,t_i,t_i^{-1}\ran$.
This implies the Serre relations (see Appendix~\ref{sec:serre}).

Thus we obtain
\Proposition
$\CF_m$ has the structure of an integrable $U_q(\Gg)$-module.
\enproposition

By Proposition~\ref{char}~(i), $\CF_m$ is
a direct sum of $V(\lam_m-k\delta)$'s.
This decomposition is studied in the next subsection through bosons.

Note that
$$\wedge:\bw^rV_\aff\otimes\CF_{m+r}\to\CF_{m}$$
is $U_q(\Gg)$-linear.

\Lemma\label{int}
$$f_i^{(k)}\vac{m}=G(\tf_i^kb_m^\circ)\wedge\vac{m+1}.$$
\enlemma

\proof
If $k>\lan h_i,\lam_m\ran$, then the both side vanish.
Assume that $0\le k\le\lan h_i,\lam_m\ran$.
By Proposition~\ref{char}~(iv),
there is $c\in K$ such that
$$f_i^{(k)}\vac{m}=cG(\tf_i^kb_m^\circ)\wedge\vac{m+1}.$$
We have
$$e_i^{(k)}f_i^{(k)}\vac{m}=
\left[{\lan h_i,\lam_m\ran}\atop k\right]_i\vac{m}.$$
On the other hand, by the repeated use of~(iii) in (G), we have
$$e_i^{(k)}G(\tf_i^kb_m^\circ)=
\left[{\lan h_i,\lam_m\ran}\atop k\right]_i\vo_m+\cdots.$$
Here, $\cdots$ is a linear combination of
global bases other than $\vo_m$, which is annihilated
after being wedged with $\vac{m+1}$ by Proposition~\ref{char} (v).
Hence we have
\beqn
e_i^{(k)}\big(G(\tf_i^kb_m^\circ)\wedge\vac{m+1}\big)
&=&\big(e_i^{(k)}G(\tf_i^kb_m^\circ)\big)\wedge\vac{m+1}\\
&=&\left[{\lan h_i,\lam_m\ran}\atop k\right]_i\vo_m\wedge\vac{m+1}.
\endeqn
Comparing these two identities, we obtain $c=1$.
\qed

\label{subsec:classical-limit}
Let $\CF_m^\BZ$ be the $\KZ$-submodule of $\CF_m$
generated by the normally ordered wedges.
Then $\CF_m^\BZ$ is a module over $U_q(\Gg)_\BZ$
by Lemma~\ref{int}.
Hence by specializing at $q=1$, we obtain
a Fock representation of $U(\Gg)$.
\hb
However, the action of the bosons on $\CF_m$
introduced in the next subsection may have a pole at $q=1$
and it cannot be specialized at $q=1$ in a na\"\i ve way.

\subsection{The action of Bosons}
We shall define the action of the bosons $B_n$ ($n\neq 0$) on $\CF_m$.

\Proposition
For $n\neq 0$ and
any $u_m,u_{m+1},\dots\in V_\aff$
such that $u_k=\vo_k$ for $k\gge m$,
\begin{equation} \label{B_n}
  \begin{aligned}
  &(z^nu_m\wedge u_{m+1}\wedge u_{m+2}\wedge\cdots)\\
  &+(u_m\wedge z^nu_{m+1}\wedge u_{m+2}\wedge\cdots)\\
  &+(u_m\wedge u_{m+1}\wedge z^nu_{m+2}\wedge\cdots)\\
  &+\cdots\cdots
\end{aligned}
\end{equation}
converges in the $q$-adic topology.
\enproposition
\proof
Reducing to the case $u_k=\vo_k$ for every $k\ge m$,
apply Lemma~\ref{qconv}.
\qed

\Lemma
$z^n\vo_m\wedge \vac{m+1}+\vo_m\wedge z^n\vo_{m+1}\wedge\vac{m+2}
+\cdots$
belongs to
$\CF_m$
\enlemma

The proof is similar to the one for Lemma~\ref{deff}.

\medskip
By these lemmas and Lemma~\ref{bos},
(\ref{B_n}) defines a homomorphism
from $\tCF_m$ to $\CF_m$.
Since
$L(\tCF_m)$ is stable by the
correspondence
(\ref{B_n}),
it induces an endomorphism of $\CF_m$.
We denote it by $B_n$.
It is clear that
$B_n$ is a $U'_q(\Gg)$-linear
endomorphism of $\CF_m$
with weight $n\delta$.

By the definition,
we have
\beq\label{copb}
&&B_n(v\wedge u)=z^nv\wedge u+v\wedge B_n(u)
\qbox{for $v\in V_\aff$ and $u\in\CF_{m}$.}
\endeq

\Proposition
There is $\gamma_n\in K$ (independent of $m$) such that
$$[B_n,B_{n'}]=\delta_{n+n',0}\gamma_n.$$
\enproposition

\proof
(\ref{copb}) implies
$$[B_n,B_{n'}](v\wedge u)=v\wedge [B_n,B_{n'}] u.$$
Since
$[B_n,B_{n'}]\vac{m}$ has weight $\lam_m+(n+n')\delta$
and hence it must vanish when $n+n'>0$.
Therefore $[B_n,B_{n'}]=0$ in this case.

Assume $n+n'<0$.
Write $[B_n,B_{n'}]\vac{m}$ as a linear combination of normally ordered wedges:
$$[B_n,B_{n'}]\vac{m}=\sum_\nu c_\nu G(b_{1,\nu})\wedge\cdots.$$
Then $b_{1,\nu}\not= b_m^\circ$.
Take $N$ and $c$ as in~(\ref{peor}).
Then we have
$$[B_n,B_{n'}]\vac{m+jN}=\sum_\nu c_\nu G(z^{jc}b_{1,\nu})\wedge\cdots.$$
We have also $H( b_{m+jN-1}^\circ\otimes z^{jc}b_{1,\nu})
=H( b_{m-1}^\circ\otimes b_{1,\nu})\le 0$.
Hence by Proposition~\ref{dual},
$\vo_m\wedge\cdots\wedge \vo_{m+jN-1}
\wedge G(z^{jc}b_{1,\nu})\wedge \vac{m+jN+1}$
converges to $0$ when $j$ tends to infinity.
Hence
\beqn
[B_n,B_{n'}]\vac{m}&=&
\vo_m\wedge\cdots\wedge\vo_{m+jN-1}\wedge [B_n,B_{n'}]\vac{m+jN}\\
&=&\sum_\nu c_\nu \vo_m\wedge\cdots\wedge \vo_{m+jN-1}\wedge
G(z^{jc}b_{1,\nu}) \wedge\cdots
\endeqn
converges to $0$.
Therefore $[B_n,B_{n'}]\vac{m}$ must vanish.

Now assume that $n+n'=0$.
Since $[B_n,B_{-n}]\vac{m}$ has the same weight as $\vac{m}$,
there is $\gamma_{m,n}$ such that
$[B_n,B_{-n}]\vac{m}=\gamma_{m,n}\vac{m}$.
Since
$$[B_n,B_{-n}]\vac{m}=\vo_m\wedge [B_n,B_{-n}]\vac{m+1}
=\gamma_{m+1,n}\vo_{m}\wedge\vac{m+1} =\gamma_{m+1,n}\vac{m},$$
$\gamma_{m,n}$ does not depend on $m$. Write $\gamma_n$ for
$\gamma_{m,n}$.  Now we have $[B_n,B_{-n}](v\wedge\vac{m})
=v\wedge[B_n,B_{-n}]\vac{m}=\gamma_n v\wedge\vac{m}$.  \qed

Now we shall show that $\gamma_n$ does not vanish.

\Lemma\label{prep}
 Let $n$ be a positive integer.
\begin{enumerate}
\item 
$z^n\vo_k\wedge\vac{k+1}=0$.
\item 
$\vo_m\wedge \vo_{m+1}\wedge\cdots\wedge \vo_{k-1}\wedge z^{-n}\vo_k
\wedge\vac{k+1}\equiv0$
for $k\ge m+n$.
\item 
$z^n\vo_m\wedge \vo_{m+1}\wedge\cdots\wedge \vo_{k-1}
\wedge z^{-n}\vo_k\wedge\vac{k+1}\equiv 0$ for $m<k<m+n$.
\end{enumerate}
Here $\equiv$ is modulo $qL(\CF_m)$.
\enlemma

\proof
(i) follows from Theorem~\ref{kern}.
\hb
In order to prove the other statements,
write $ b_k=z^{-k}b_k^\circ$.
Then $H( b_k\otimes  b_{k+1})=0$.

We have
\beqn
&&\vo_m\wedge \vo_{m+1}\wedge\cdots\wedge \vo_{k-1}\wedge z^{-n}\vo_k
\equiv\\
&&\phantom{-----------}
z^{m} b_{m}\wedge z^{1+m} b_{m+1}\wedge
\cdots\wedge z^{k-1} b_{k-1}\wedge z^{k-n} b_{k}.
\endeqn
Since $m\le k-n\le k-1$, it is zero modulo $qL(\bw V_\aff)$
by Lemma~\ref{lem:9} (iii).

\noindent
The proof of~(iii) is similar to that of~(ii).
We have
\begin{multline*}
  z^n\vo_m\wedge \vo_{m+1}\wedge\cdots\wedge \vo_{k-1}
  \wedge z^{-n}\vo_k\wedge \vo_{k+1}\wedge\cdots\wedge \vo_{n+m}\equiv\\
  z^{n+m} b_{m}\wedge z^{m+1} b_{m+1}\wedge
  z^{k-1} b_{k-1}
  \wedge z^{k-n} b_{k}\wedge z^{k+1} b_{k+1}\wedge\cdots
  \wedge z^{n+m} b_{n+m}.
\end{multline*}
Then it is zero modulo $qL(\bw V_\aff)$
again by Lemma~\ref{lem:9}~(iii).
\qed

\Proposition
For $n\neq0$, $\gamma_n\in K$ has no pole at $q=0$ and
$\gamma_n(0)=n$.
\enproposition

\proof
We may assume $n>0$.
Noting that
$B_{\pm n}$ sends $L(\CF_m)$ to itself,
let us calculate the commutator modulo $qL(\CF_m)$.
We have
$[B_n,B_{-n}]\vac{m}=B_nB_{-n}\vac{m}$.
By Lemma~\ref{prep}~(ii), we have
\beqn
B_{-n}\vac{m}&=&
z^{-n}\vo_m\wedge\vac{m+1}+\vo_m\wedge z^{-n}\vo_{m+1}\wedge\vac{m+2}+\cdots\\
&\equiv&\sum_{m\le k<m+n}
\vo_m\wedge \vo_{m+1}\wedge\cdots\wedge \vo_{k-1}\wedge z^{-n}\vo_k
\wedge\vac{k+1}.
\endeqn
Here $\equiv$ is taken modulo $qL(\CF_m)$.
Hence we have by Lemma~\ref{prep}~(i) and~(iii)
\begin{align*}
  &B_nB_{-n}\vac{m}\cr
  &\ \equiv
  \sum\limits_{m\le k<m+n}
  \Bigl(
  \sum\limits_{m\le j<k}\vo_m\wedge \vo_{m+1}\wedge\cdots\wedge
  z^n\vo_j\wedge\cdots\wedge \vo_{k-1}\wedge z^{-n}\vo_k
  \wedge\vac{k+1}\cr
  &\phantom{-----}+\vac{m}\cr
  &\phantom{-----}+\sum\limits_{j>k}
  \vo_m\wedge \vo_{m+1}\wedge\cdots\wedge \vo_{k-1}\wedge z^{-n}\vo_k
  \wedge \vo_{k+1}\wedge\cdots\wedge z^n\vo_j\wedge\vac{j+1}
  \Bigr)\cr
  &\ \equiv n\vac{m}. 
\end{align*}
\qed

Let $H$ be the Heisenberg algebra generated by
$\{B_n\}_{n\in\BZ\backslash\{0\}}$  with
the defining relations
$[B_n,B_{n'}]=\delta_{n+n',0}\gamma_n$.
Then $H$ acts on the Fock space $\CF_m$
commuting with the action of $U'_q(\Gg)$.
Let $\BQ[H_-]$ be the Fock space for $H$.
Namely, $\BQ[H_-]$ is the $H$-module generated by the vacuum vector
$1$ with the defining relation
$B_n1=0$ for $n>0$.
Since $\vac{m}$ is annihilated by the $e_i$ and
the $B_n$ with $n>0$,
we have an injective $U'_q(\Gg)\otimes H$--linear homomorphism
\beq
\iota_m:V(\lam_m)\otimes \BQ[H_-]\to \CF_m\label{dec}
\endeq
sending $u_{\lam_m}\otimes 1$ to $\vac{m}$.
Comparing their characters
(see Proposition~\ref{char}~(i)),
we obtain

\Theorem
$\iota_m:V(\lam_m)\otimes \BQ[H_-]\to \CF_m$ is an isomorphism.
\entheorem

\subsection{Vertex operator} \label{subsec:vo}
Similarly to the case $A^{(1)}_n$ in~\cite{KMS},
the intertwiner
\begin{align*}
  \mbox{}\qquad\qquad \Omega_m:V_\aff\otimes\CF_{m+1}&\to\CF_{m},\\
  v\otimes u &\mapsto v\wedge u \qquad (v\in\Vaff, u\in\CF_{m+1}),
\end{align*}
induced by the wedge product
is related with vertex operators.
Let us describe it briefly. The proof is similar to
\cite{KMS}.

Take an intertwiner
\beq
\Phi_m:V_\aff\otimes V(\lam_{m+1})\to V(\lam_{m})
\endeq
and normalize it by
$$\Phi_m(\vo_{m}\otimes u_{\lam_{m+1}})=u_{\lam_m}$$
(cf. Appendix~\ref{perf}).

Let
$$\iota_m: V(\lam_m)\otimes \BQ[H_-]{\buildrel\sim\over\longrightarrow}
\CF_m$$
be the isomorphism in~(\ref{dec}).

We define
$$\Omega'_m:V_\aff\otimes V(\lam_{m+1})\otimes\BQ[H_-]
\to V(\lam_{m})\otimes \BQ[H_-]$$
by requiring the commutativity of the following diagram
\begin{equation}
  \begin{CD} \label{DGRM}
    V_\aff\otimes V(\lambda_{m+1})\otimes\BQ[H_-]
    @>{\sim}>{\id \otimes\iota_{m+1}}>
    V_\aff\otimes \CF_{m+1}\\
    @VV{\Omega'_m}V   @VV{\Omega_m}V\\
    V(\lambda_m)\otimes\BQ[H_-] @>\sim>\iota_m> \CF_{m}
  \end{CD}
  \quad .
\end{equation}
We shall write the intertwiners in the form of generating
functions.
Namely, introducing an indeterminate $w$
(of weight $\delta$),
we set for $v\in V_\aff$
\begin{equation}
  \begin{matrix}
\hfill v(w)&=&\sum_nz^{n}v\otimes w^{-n},\nn\cr
\Phi_m(w)(v\otimes u)&=&\Phi_m(v(w)\otimes u)
&=&\sum_n\Phi_m(z^{n}v\otimes u)w^{-n},\nn\cr
\Omega_m(w)(v\otimes u)&=&
\Omega_m(v(w)\otimes u)&=&
\sum_n\Omega_m(z^{n}v\otimes u)w^{-n},\nn\cr
\Omega'_m(w)(v\otimes u)&=&
\Omega'_m(v(w)\otimes u)&=&\sum_n\Omega'_m(z^{n}v\otimes u)w^{-n}.\nn
\end{matrix}
\end{equation}
Here $u\in V(\lam_{m+1})$, $\CF_{m+1}$ or
$V(\lambda_{m+1})\otimes\BQ[H_-]$.

We define the vertex operator for the bosons by
\begin{equation}
\Theta(w)=
\exp\left( \sum_{n\ge1}{B_{-n}w^n\over\gamma_n}\right)
\exp\left( -\sum_{n\ge1}{B_{n}w^{-n}\over\gamma_n}\right). \label{D}
\end{equation}

\Theorem\label{TEIRI}
$\Omega'_m(w)=\Phi_m(w)\otimes \Theta(w)$.
\entheorem

As a corollary of this theorem,
we have the relations of the two-point functions of
the vertex operators and $\gamma_n$ as in~\cite{KMS}.

Set
$$\Phi_m^v(w)(u)=\Phi_m(w)(v\otimes u)$$
for $u\in\CF_{m+1}$.

For $v,v'\in V_\aff$, we
define
$\lan\Phi^v_{m-1}(w_1)\Phi^{v'}_m(w_2)\ran$ to be
the coefficient of $u_{\lam_{m-1}}$ in
$\Phi^v_{m-1}(w_1)\Phi^{v'}_m(w_2)u_{\lam_{m+1}}\in V(\lam_{m-1})$.
We introduce functions by
\beq
\omega_{v,v'}(w_2/w_1)&=&
\lan{m-1}|v(w_1)\wedge v'(w_2)\wedge\vac{m+1} \label{eq:fock-2pt}\\
\phi_{v,v'}(w_2/w_1)&=&
\lan\Phi^v_{m-1}(w_1)\Phi^{v'}_m(w_2)\ran \label{eq:vertex-2pt}\\\\
\noalign{\hbox{and}}
\theta(w_2/w_1)&=&\exp\left(-\sum_{n>0}{{(w_2/w_1)^n}\over{\gamma_n}}\right).
\label{eq:boson-2pt}
\endeq
Here for $u\in \CF_{m-1}$,
$\lan{m-1}|u$ means the coefficient of $\vac{m-1}$ in $u$.
Then the theorem above implies
\Proposition\label{prop:2pt-decomp}
For $v,v'\in V_\aff$, we have
\beq
\omega_{v,v'}(w_2/w_1)=\phi_{v,v'}(w_2/w_1)\theta(w_2/w_1).
\endeq
\enproposition

This formula will be used later to calculate $\gamma_n$.



\newsection{Examples of level~1 Fock spaces}
In this section we give some examples of the theory developed in the
earlier sections.  The case of
level~$1$ type $A^{(1)}_n$ described in~\cite{S,KMS} is first reviewed
in the perfect crystal language.  Then we
present results for types $A^{(2)}_{2n}$, $B^{(1)}_n$,
$A^{(2)}_{2n-1}$, $D^{(1)}_n$ and $D^{(2)}_{n+1}$ at level~$1$,
corresponding to the perfect crystals of~\cite{KMN1}~Table~2.

\subsection{Preliminaries}\label{GOITI}
Define $\range{m,n}:=\set{i\in\Z\mid m\leq i\leq n}$. We label the
simple roots by $I=\range{0,n}$.  We choose $0\in I$ so that $W_\cl$ is
generated by $\set{s_i}_{i\in I\setminus\set{0}}$ and $a_0=1$.

We take fundamental weights $\set{\Lambda_i}_{i\in I}$ such that
\begin{eqnarray*}
\alpha_i&=&\sum_{j\in I}\langle h_j,\alpha_i\rangle\Lambda_j
+\delta_{i,0}\delta.
\end{eqnarray*}
Let $s_0:P^0_\cl\to P^0$ be a section of $\cl:P^0\to P^0_\cl$ such that
\begin{displaymath}
  s_0(P^0_\cl)\subset\sum_{i\in I\setminus\set{0}}\BQ\alpha_i
  =\sum_{i\in I\setminus\set{0}}\BQ(a^\vee_0\Lambda_i-a^\vee_i\Lambda_0).
\end{displaymath}
Then we have
\begin{displaymath}
  s_0(\lambda+\cl(\alpha_i))=
  \begin{cases}
    s_0(\lambda)+\alpha_i &\text{for } i\in I\setminus\set{0},\\
    s_0(\lambda)+\alpha_0 -\delta &\text{for } i=0.
  \end{cases}
\end{displaymath}
We regard $V$ as a subspace of $V_\aff$ by
$V\supset V_\lambda\simeq (V_\aff)_{s_0(\lambda)}\subset V_\aff$.
Then $V_\aff$ is identified with
$\BQ[z,z^{-1}]\otimes V$.
With this identification, the action of $U_q(\Gg)$
on $\BQ[z,z^{-1}]\otimes V$
is given by
\begin{eqnarray*}
e_i(a\otimes v)&=&z^{\delta_{i,0}}a\otimes e_iv,\\
f_i(a\otimes v)&=&z^{-\delta_{i,0}}a\otimes f_iv.
\end{eqnarray*}
Similarly we identify $B$ as a subset of $B_\aff$.

In the examples that we treat in this paper, the action of $U_q(\Gg)$
on the lower global base of $\Vaff$ (respectively $V$) is completely
determined by its crystal structure as we have
\begin{equation} \label{eq:LG}
  \begin{aligned}
  e_iG(b)&=[1+\varphi_i(b)]_iG(\te_ib),\\
  f_iG(b)&=[1+\varepsilon_i(b)]_iG(\tf_ib),\\
  q^hG(b)&= q^{\pairing{h,\wt(b)}}G(b),
  \end{aligned}
\end{equation}
for $b\in\Baff$ ($b\in B$), $i\in I$ and $h\in P^*$ ($h\in P_\cl^*$).

\subsection{Level~1 $A^{(1)}_n$}

\subsubsection{Cartan datum}
The Dynkin diagram for $A^{(1)}_n$ ($n\geq 1$) is
\begin{eqnarray*}
  \dynkin{0} \squeeze
  \!\!\!\longlink\dynkin{1}\longlink\!\!\!\!\squeeze
  \dynkin{2}\\[-6pt]
  | \squeeze \squeeze |\\[-6pt]
  \dynkin{n} \squeeze \squeeze \dynkin{3} \qquad .\\[-6pt]
  | \squeeze \squeeze |\\[-6pt]
  \dynkin{n-1}\!\!\!\squeeze\;\: \link\! \cdots \link \squeeze \dynkin{4}
\end{eqnarray*}
For $A^{(1)}_n$ we have
\begin{eqnarray*}
  \delta &=& \sum_{i\in I}\alpha_i,\\
  c &=& \sum_{i\in I}h_i,\\
  (\alpha_i,\alpha_i) &=& 2 \qquad (i\in I).
\end{eqnarray*}

\subsubsection{Perfect crystal}
Let $J:=\range{0,n}$.  Let $V$ be the
$(n+1)$-dimensional $U'_q(A^{(1)}_n)$-module with the level~$1$ perfect
crystal $B:=\set{\b{i}}_{i\in J}$ with crystal graph:
\begin{eqnarray*}
  \b{0} \squeeze \rightcrystal{1}\; \b{1}\; \rightcrystal{2}
  \squeeze \b{2}\\[-5pt]
  {\scriptstyle 0}\!\! \uparrow\, \squeeze
  \squeeze\downarrow\!\!{\scriptstyle 3}\\[-7pt]
  \b{n}\squeeze \squeeze \hspace{0.375mm}\vdots \qquad .\\[-5pt]
  {\scriptstyle n}\!\! \uparrow\, \squeeze \squeeze
  \downarrow\!\!{\scriptstyle n-3}\\[-5pt]
  \b{n-1}\squeeze \leftcrystal{n-1} \b{n-2}
  \leftcrystal{n-2} \squeeze \b{n-3}
\end{eqnarray*}
The elements of $B$ have the following weights
\begin{displaymath}
  \wt(\b{i}) = \Lambda_{i+1} -\Lambda_{i}
  \qquad (i\in J).
\end{displaymath}
Let $\v{j}:=G(\b{j})$ ($j\in J$).  The action of
$U'_q(A^{(1)}_n)$ on $\v{j}\in V$ obeys~\eqref{eq:LG}.

\subsubsection{Energy function}
The energy function $H$ takes the following values on $B\tensor B$
\begin{eqnarray*}
  H(\b{i}\tensor \b{j})=
  \begin{cases}
    1 &\text{for } i>j,\\
    0 &\text{for } i\leq j.
  \end{cases}
\end{eqnarray*}
Write $H(i,j)$ for $H(\b{i}\tensor\b{j})$ ($i,j\in J$).

The Coxeter number of $A^{(1)}_n$ is $h=n+1=\dim V$.  We take
$l:\Baff\mapto\Z$ to be
\begin{displaymath}
  l(z^m\b{j}) = mh-j\qquad (m\in\Z,j\in J).
\end{displaymath}
The functions $H$ and $l$ satisfy condition~(L) (see end of
\S\ref{sec:ener}).  The map $l$ gives a total ordering of $\Baff$.

\subsubsection{Wedge relations}
We have
\begin{displaymath}
  N:=U_q(A^{(1)}_n){[z\tensor z, z^{-1}\tensor z^{-1},z\tensor1 +
    1\tensor z]}\cdot \v{0}\tensor\v{0}\subset \Vaff\tensor\Vaff.
\end{displaymath}
The following elements are contained in
$U_q(A^{(1)}_n)\cdot\v{0}\tensor\v{0}\subset N$:
\begin{eqnarray*}
  C_{i,i} &=& \v{i}\tensor\v{i}
  \qquad\qquad\qquad\qquad\qquad\qquad (i\in J),\\
  C_{i,j} &=& \v{i}\tensor z^{-H(i,j)}\v{j}\\
  &&\quad +q z^{-H(i,j)} \v{j}\tensor\v{i} \qquad \left((i,j)\in
    J^2\setminus\set{(k,k)}_{k\in J}\right).
\end{eqnarray*}

\begin{proposition}
  Identify $C_{i,j}$ with $C_{\b{i}, z^{-H(i,j)}\b{j}}$.  Then
  $\set{z^m\tensor z^m\cdot C_{i,j}}_{m\in\Z;i,j\in J}$ with the
  function $l$ satisfy condition~(R) of
  subsection~\ref{subsec:wedge-prod}.
\end{proposition}

\subsubsection{Fock space}
For $U_q(A^{(1)}_n)$ we have
\begin{eqnarray*}
  \Bmin&=&B,\\
  (\Pcl^+)_1&=&\set{\Lambda_i^\cl}_{i\in I},
\end{eqnarray*}
with
\begin{displaymath}
  \e(\b{j})= \Lambda_j^\cl, \qquad
  \f(\b{j})= \Lambda_{j+1\bmod h}^\cl\qquad (j\in J).
\end{displaymath}
Since $H(\b{j}\tensor \b{j-1})=1$ ($j\in\range{1,n}$) and
$H(\b{0}\tensor z\b{n})=1$
there is a unique ground state sequence given as follows:  every
$m\in\Z$ fixes uniquely $a\in\Z$ and $j\in J$ such that $m=ah-j$, then
\begin{eqnarray*}
  \bo_m &=& z^a\b{j} \qquad (m\in\Z),\\
  \cl(\lambda_m) &=& \Lambda_{j+1\bmod h} \qquad (m\in\Z).
\end{eqnarray*}
With $\vo_m=G(\bo_m)$, the vacuum vector of $\fock_m$ is then given by
\begin{displaymath}
  \vac{m}= \vo_m\wedge \vo_{m+1}\wedge \vo_{m+2}\wedge\cdots\cdots
\end{displaymath}
with highest weight $\lambda_m$.

\subsection{Level~1 $A^{(2)}_{2n}$}

\subsubsection{Cartan datum}
The Dynkin diagram for $A^{(2)}_{2n}$ ($n\geq1$) is
\begin{displaymath}
  \dynkin{0}\rightdoublelink \dynkin{1} \link \dynkin{2} \link \cdots
  \cdots  \link \dynkin{n-2} \link \dynkin{n-1}
  \rightdoublelink \dynkin{n}\;.
\end{displaymath}
For $A^{(2)}_{2n}$ we have
\begin{eqnarray*}
  \delta &=& \alpha_0 + \sum_{i=1}^n 2\alpha_i,\\
  c &=& (\sum_{i=0}^{n-1}2h_i)+h_n,\\
  (\alpha_i,\alpha_i) &=&
  \begin{cases}
    8 &\text{for } i=0,\\
    4 &\text{for } i\in\range{1,n-1},\\
    2 &\text{for } i=n.
  \end{cases}
\end{eqnarray*}

\subsubsection{Perfect crystal}
Let $J:=\range{-n,n}$.  Let $V$ be the $(2n+1)$-dimensional
$U'_q(A^{(2)}_{2n})$-module with the level~$1$ perfect crystal
$B:=\set{\b{i}}_{i\in J}$ and crystal graph:
\begin{eqnarray*}
  \b{1} \squeeze \rightcrystal{1} \b{2} \rightcrystal{2} \cdots\cdots
  \rightcrystal{n-2} \b{n-1}\rightcrystal{n-1} \squeeze \b{n}\\[-6pt]
  \makebox[1.1ex]{$\uparrow$} \, \squeeze
  \squeeze\downarrow\!\!{\scriptstyle n}\\[-6.5pt]
  {\scriptstyle 0} \vert \;\squeeze \squeeze \b{0}
  \qquad .\\[-6pt]
  | \; \squeeze \squeeze \downarrow\!\!{\scriptstyle n}\\[-6pt]
  \b{-1}\squeeze \leftcrystal{1} \b{-2}
  \leftcrystal{2} \cdots\cdots \leftcrystal{n-2}
  \b{1-n}\leftcrystal{n-1} \squeeze \b{-n}
\end{eqnarray*}

The elements of $B$ have the following weights
\begin{eqnarray*}
  \wt(\b{i}) &=& \sum_{k=i}^{n} \alpha_k
  = (1+\delta_{i,n})\Lambda_i -\Lambda_{i-1}\qquad
  (i\in\range{1,n}),\\
  \wt(\b{0}) &=& 0,\\
  \wt(\b{-i}) &=& -\wt(\b{i}) \qquad (i\in\range{1,n}).
\end{eqnarray*}
Let $\v{j}:=G(\b{j})$ ($j\in J$).  The action of
$U'_q(A^{(2)}_{2n})$ on $\v{j}\in V$ obeys~\eqref{eq:LG}.

\subsubsection{Energy function}
Define the following ordering of $J$
\begin{displaymath}
  1 \succ 2 \succ \cdots \succ n \succ 0 \succ -n \succ 1-n \succ \dots
  \succ -1.
\end{displaymath}

The energy function $H$ takes the following values on $B\tensor B$
\begin{eqnarray*}
  H(\b{i}\tensor\v{j})=
  \begin{cases}
    1 &\text{for } (i,j)\in\set{(i',j')\in J^2\mid i'\prec
j'}\union\set{(0,0)},\\
    0 &\text{for } (i,j)\in\set{(i',j')\in J^2\mid i'\succ j'}\union
    \set{(k,k)}_{k\in J\setminus\set{0}}.
  \end{cases}
\end{eqnarray*}
Write $H(i,j)$ for $H(\b{i}\tensor\b{j})$ ($i,j\in J$).

The Coxeter number of $A^{(2)}_{2n}$ is $h=2n+1=\dim V$.  We take
$l:\Baff\mapto\Z$ to be
\begin{displaymath} \label{eq:a2even-l-func}
  l(z^m\b{j}) =
  \begin{cases}
    hm+n+1-j &\text{for } j\in\range{1,n},\\
    hm &\text{for } j=0,\\
    hm-(n+1+j) &\text{for } j\in\range{-n,-1}.
  \end{cases}
\end{displaymath}
The functions $H$ and $l$ satisfy condition~(L) (see end of
\S\ref{sec:ener}).  The map $l$ gives a total ordering of $\Baff$.

\subsubsection{Wedge relations}
In $\Vaff\tensor\Vaff$ we have
\begin{displaymath}
  N:=U_q(A^{(2)}_{2n}){[z\tensor z, z^{-1}\tensor z^{-1},z\tensor1 +
    1\tensor z]}\cdot \v{1}\tensor\v{1}.
\end{displaymath}
The following elements are contained in
$U_q(A^{(2)}_{2n})\cdot\v{1}\tensor\v{1}\subset N$:
\begin{eqnarray*}
  \tilde{C}_{i,i} &=& \v{i}\tensor\v{i}
  \qquad\qquad\qquad\qquad\qquad\qquad
  \left(i\in J\setminus\set{0}\right),\\
  \tilde{C}_{i,-i} &=& \v{i}\tensor z^{-H(i,-i)} \v{-i} + q^2 \v{i+1}\tensor
  z^{-H(i,-i)} \v{-i-1}\\
  && \quad +q^2z^{-H(i,-i)}\v{-i-1}\tensor\v{i+1}\\
  && \quad + q^4 z^{-H(i,-i)} \v{-i}\tensor\v{i}\qquad\qquad \left(i\in
    J\setminus\set{-1,0,n}\right),\\
  \tilde{C}_{i,j} &=& \v{i}\tensor z^{-H(i,j)}\v{j}\\
  && \quad +q^2z^{-H(i,j)}\v{j}\tensor\v{i} \qquad \left((i,j)\in
    J^2\setminus\set{(k,k),(k,-k)}_{k\in J}\right),\\
  \tilde{C}_{0,0} &=& \v{0}\tensor z^{-1} \v{0} + q^2[2]\v{-n}\tensor
  z^{-1} \v{n} \\
  && \quad+ q^2[2]z^{-1}\v{n}\tensor\v{-n} + q^2
  z^{-1}\v{0}\tensor\v{0}, \\
  \tilde{C}_{n,-n} &=& \v{n}\tensor\v{-n} + q\v{0}\tensor\v{0} + q^4
  \v{-n}\tensor\v{n}, \\
  \tilde{C}_{-1,1} &=& \v{-1}\tensor z^{-1} \v{1} + q^4
  z^{-1} \v{1}\tensor\v{-1}.
\end{eqnarray*}
Notice that each $\tilde{C}_{i,j}$ has $\v{i}\tensor z^{-H(i,j)}\v{j}$
as its first term and a term in $z^{-H(i,j)}\v{j}\tensor\v{i}$.

Define the following elements in $N$.
\begin{displaymath}
  C_{i,j}:=
  \begin{cases}
    \tilde{C}_{i,j} &\text{for } (i,j)\in J^2\setminus\set{(k,-k)}_{k\in J},\\
    \sum_{k=i}^n(-q^2)^{k-i}\tilde{C}_{k,-k} &\text{for }
    (i,j)\in\set{(k,-k)}_{k\in\range{1,n}},\\
    \sum_{k=1}^j(-q^2)^{j-k}\tilde{C}_{-k,k} &\text{for }
    (i,j)\in\set{(-k,k)}_{k\in\range{1,n}},\\
    \tilde{C}_{0,0} -q^2[2]C_{-n,n} &\text{for } (i,j)=(0,0).
  \end{cases}
\end{displaymath}
Explicitly for $(i,j)\in\set{(k,-k)}_{k\in J}$, we have
\begin{eqnarray*}
  C_{j,-j} &=& \v{j}\tensor\v{-j} + q^4 \v{-j}\tensor\v{j}
  +q (-q^2)^{n-j}\v{0}\tensor\v{0} \\
  &&\quad -(1-q^4) \sum_{k=j+1}^n (-q^2)^{k-j} \v{-k}\tensor\v{k} \qquad
  (j\in\range{1,n}),\\
  C_{-j,j} &=& \v{-j}\tensor z^{-1} \v{j} + q^4 z^{-1}\v{j}\tensor
  \v{-j}\\
  &&\quad -(1-q^4)\sum_{k=1}^{j-1} (-q^2)^{j-k} z^{-1}\v{k}\tensor\v{-k}
  \qquad (j\in\range{1,n}),\\
  C_{0,0} &=& \v{0}\tensor z^{-1} \v{0} + q^2 z^{-1} \v{0}\tensor
  \v{0}\\
  &&\quad +q^2[2](1-q^4) \sum_{k=1}^n (-q^2)^{n-k} z^{-1} \v{k}\tensor
  \v{-k}.
\end{eqnarray*}

\begin{proposition}
  Identify $C_{i,j}$ with $C_{\b{i}, z^{-H(i,j)}\b{j}}$.  Then
  $\set{z^m\tensor z^m\cdot C_{i,j}}_{m\in\Z;i,j\in J}$ with the
  function $l$ satisfy condition~(R) of
  subsection~\ref{subsec:wedge-prod}.
\end{proposition}

\subsubsection{Fock space}
We have $\Bmin=\set{\b{0}}$ and $(\Pcl^+)_1=\set{\Lambda_n^\cl}$.
Since $H(\b{0},\b{0})=1$ we have a unique ground state sequence (up
to an overall shift by $z^{k}$ ($k\in\Z$)): $\bo_m=\b{0}$ and
$\lambda_m=\Lambda_n$ ($m\in\Z$).  Therefore the vacuum vector of
the Fock space $\fock_m$ is
\begin{displaymath}
  \vac{m}:= \v{0}\wedge \v{0}\wedge \v{0}\wedge \v{0}\wedge
  \cdots\cdots \qquad (m\in\Z).
\end{displaymath}
We set $\wt(\vac{m})=\Lambda_n$.

As an illustration of the use of the $q$-adic topology, let us check
Proposition~\ref{prop:12} on $\vac{m}$: i.e.\ that ${[e_i,f_i]}\cdot
\vac{m}= \frac{t_i-t_i^{-1}}{q_i-q_i^{-1}}\cdot \vac{m}$. The case
$i\in I\setminus\set{n}$ is trivial.  For $e_n\vac{m}$, consider first
$\v{n}\wedge\vac{m+1}=(-q^2)^j (\v{0})^{\wedge j}\wedge
\v{n}\wedge\vac{m+j+1}$ ($j\in\N$). As $j\mapto\infty$, the vector
vanishes by the $q$-adic topology on $\fock_m$.  Hence
\begin{eqnarray*}
  e_n\cdot\vac{m} &=& \sum_{j=0}^\infty (\v{0})^{\wedge j}
  \wedge (e_n\cdot\v{0})\wedge \vac{m+j+1}\\
  &=& [2] \sum_{j=0}^\infty (\v{0})^{\wedge j}\wedge \v{n}\wedge
  \vac{m+j+1}\\
  &=& 0.
\end{eqnarray*}
For $f_n$ we have
  \begin{eqnarray*}
    f_n\cdot\vac{m} &=& \sum_{j=0}^\infty (\v{0})^{\wedge
      j}\wedge (f_n\cdot \v{0})\wedge t_n\vac{m+j+1}\\
    &=& q[2] \sum_{j=0}^\infty (\v{0})^{\wedge
      j}\wedge \v{-n}\wedge \vac{m+j+1}\\
    &=& q[2] \sum_{j=0}^\infty (-q^2)^j \v{-n}\wedge\vac{m+1}\\
    &=& q[2](1+q^2)^{-1} \v{-n}\wedge\vac{m+1}\\
    &=& \v{-n}\wedge \vac{m+1}.
  \end{eqnarray*}
Then since
\begin{eqnarray*}
  e_n\cdot f_n\cdot\vac{m} &=& e_n\cdot \v{-n}\wedge\vac{m+1}\\
  &=& \v{0}\wedge\vac{m+1}\\
  &=& \vac{m},
\end{eqnarray*}
and $[\pairing{h_n,\Lambda_n}]_n=1$, this completes the check.

\subsection{Level~1 $B^{(1)}_n$}

\subsubsection{Cartan datum}
The Dynkin diagram for $B^{(1)}_n$ ($n\geq 3$) is
\begin{eqnarray*}
  \dynkin{0}\link \!\!\!\! &\dynkin{2}& \!\!\! \link \dynkin{3} \link \cdots
  \cdots \link \dynkin{n-2}  \link \dynkin{n-1}
  \rightdoublelink \dynkin{n}\;.\\[-5pt]
  &|&\\[-5pt]
  &\dynkin{1}&
\end{eqnarray*}
For $B^{(1)}_n$ we have
\begin{eqnarray*}
  \delta &=& \alpha_0+ \alpha_1+ \sum_{i=2}^n 2\alpha_i,\\
  c &=& h_0+ h_1+ (\sum_{i=2}^{n-1}2h_i)+ h_n,\\
  (\alpha_i,\alpha_i) &=&
  \begin{cases}
    4 &\text{for } i\in\range{0,n-1},\\
    2 &\text{for } i=n.
  \end{cases}
\end{eqnarray*}

\subsubsection{Perfect crystal}
Let $J:=\range{-n,n}$.  Let $V$ be the $(2n+1)$-dimensional
$U'_q(B^{(1)}_n)$-module with the level~$1$ perfect crystal
$B:=\set{\b{i}}_{i\in J}$ with crystal graph:
\begin{eqnarray*}
  \b{2}\: \squeeze \rightcrystal{2} \b{3} \rightcrystal{3} \cdots\cdots
  \rightcrystal{n-2} \b{n-1}\rightcrystal{n-1} \squeeze \b{n}\\[-5pt]
  {\scriptstyle 1}\!\!\nearrow\; \nwarrow\!\!{\scriptstyle 0}
  \!\!\!\!\squeeze \squeeze\downarrow\!\!{\scriptstyle n}\\[-5pt]
  \b{1}\;\;\;\;\;\b{-1} \!\!\!\!\!\!\!\! \squeeze \squeeze \b{0}
  \qquad  .\\[-5pt]
  {\scriptstyle 0}\!\!\nwarrow\; \nearrow\!\!{\scriptstyle 1}
  \!\!\!\!\squeeze\squeeze \downarrow\!\!{\scriptstyle
    n}\\[-5pt]
  \b{-2}\squeeze \leftcrystal{2} \b{-3}
  \leftcrystal{3} \cdots\cdots \leftcrystal{n-2}
  \b{1-n}\leftcrystal{n-1} \squeeze \b{-n}
\end{eqnarray*}
The elements of $B$ have the following weights
\begin{eqnarray*}
  \wt(\b{i}) &=& \sum_{k=i}^{n} \alpha_k
  = (1+\delta_{i,n})\Lambda_i -\Lambda_{i-1}
  -\delta_{i,2}\Lambda_0 \qquad (i\in\range{1,n}),\\
  \wt(\b{0}) &=& 0,\\
  \wt(\b{-i}) &=& -\wt(\b{i}) \qquad (i\in\range{1,n}).
\end{eqnarray*}
Let $\v{j}:=G(\b{j})$ ($j\in J$).  The action of
$U'_q(B^{(1)}_n)$ on $\v{j}\in V$ obeys~\eqref{eq:LG}.

\subsubsection{Energy function}
Define the following ordering of $J$
\begin{displaymath}
  1 \succ 2 \succ \cdots \succ n \succ 0 \succ -n \succ 1-n \succ \dots
  \succ -1.
\end{displaymath}
The energy function $H$ takes the following values on $B\tensor B$
\begin{eqnarray*}
  H(\b{i}\tensor\b{j})=
  \begin{cases}
    2 &\text{for } (i,j)=(-1,1),\\
    1 &\text{for } (i,j)\in\bigset{(i',j')\in J^2\setminus\set{(-1,1)} \mid
i'\prec
      j'}\union\set{(0,0)},\\
    0 &\text{for } (i,j)\in\set{(i',j')\in J^2\mid i'\succ
      j'}\union\set{(k,k)}_{k\in J\setminus\set{0}}.
  \end{cases}
\end{eqnarray*}
Write $H(i,j)$ for $H(\b{i}\tensor\b{j})$ ($i,j\in J$).

The Coxeter number of $B^{(1)}_n$ is $h=2n=\dim V-1$.  We take
$l$ to be
\begin{displaymath}
  l(z^m\b{j}) =
  \begin{cases}
    hm+n+1-j &\text{for } j\in\range{1,n},\\
    hm &\text{for } j=0,\\
    hm-(n+1+j) &\text{for } j\in\range{-n,-1}.
  \end{cases}
\end{displaymath}
The functions $H$ and $l$ satisfy condition~(L).  Note that
$l(z^m\b{1})=l(z^{m+1}\b{-1})$ ($m\in\Z$), so the map~$l$ gives a
partial ordering of $\Baff$.

\subsubsection{Wedge relations}
We have
\begin{displaymath}
  N:=U_q(B^{(1)}_n){[z\tensor z, z^{-1}\tensor z^{-1},z\tensor1 +
    1\tensor z]}\cdot \v{1}\tensor\v{1}\subset \Vaff\tensor\Vaff.
\end{displaymath}
The following elements are contained in
$U_q(B^{(1)}_n)\cdot\v{1}\tensor\v{1}\subset N$:
\begin{eqnarray*}
  \tilde{C}_{i,i} &=& \v{i}\tensor\v{i}
  \qquad\qquad\qquad\qquad\qquad\qquad
  \left(i\in J\setminus\set{0}\right),\\
  \tilde{C}_{i,-i} &=& \v{i}\tensor z^{-H(i,-i)} \v{-i} + q^2 \v{i+1}\tensor
  z^{-H(i,-i)} \v{-i-1}\\
  &&\quad +q^2z^{-H(i,-i)} \v{-i-1}\tensor\v{i+1} \\
  &&\quad +q^4 z^{-H(i,-i)} \v{-i}\tensor\v{i} \qquad\qquad \left(i\in
    J\setminus\set{-1,0,n}\right),\\
  \tilde{C}_{i,j} &=& \v{i}\tensor z^{-H(i,j)}\v{j}\\
  &&\quad +q^2 z^{-H(i,j)} \v{j}\tensor\v{i} \qquad \left((i,j)\in
    J^2\setminus\set{(k,k),(k,-k)}_{k\in J}\right),\\
  \tilde{C}_{0,0} &=& \v{0}\tensor z^{-1} \v{0} + q^2[2]\v{-n}\tensor
  z^{-1} \v{n} \\
  &&\quad +q^2[2]z^{-1}\v{n}\tensor\v{-n} + q^2
  z^{-1}\v{0}\tensor\v{0}, \\
  \tilde{C}_{n,-n} &=& \v{n}\tensor\v{-n} + q\v{0}\tensor\v{0} + q^4
  \v{-n}\tensor\v{n}, \\
  \tilde{C}_{-1,1} &=& \v{-1}\tensor z^{-2}\v{1} + q^2
  z^{-1}\v{-2}\tensor z^{-1}\v{2}\\
  &&\quad +q^2 z^{-1}\v{2}\tensor z^{-1}\v{-2} +q^4
  z^{-2}\v{1}\tensor\v{-1}.
\end{eqnarray*}
Each $\tilde{C}_{i,j}$ has $\v{i}\tensor z^{-H(i,j)}\v{j}$
as its first term and a term in $z^{-H(i,j)}\v{j}\tensor\v{i}$.

Define the following elements in $N$.
\begin{displaymath}
  C_{i,j}:=
  \begin{cases}
    \tilde{C}_{i,j} &\text{for } (i,j)\in J^2\setminus\set{(k,-k)}_{k\in J},\\
    \sum_{k=i}^n(-q^2)^{k-i}\tilde{C}_{k,-k} &\text{for }
    (i,j)\in\set{(k,-k)}_{k\in\range{1,n}},\\
    \sum_{k=2}^j(-q^2)^{j-k}\tilde{C}_{-k,k} &\text{for }
    (i,j)\in\set{(-k,k)}_{k\in\range{2,n}},\\
    \tilde{C}_{0,0} -q^2[2] C_{-n,n} &\text{for } (i,j)=(0,0),\\
    \tilde{C}_{-1,1}- q^2 (z^{-1}\tensor z^{-1}) C_{2,-2} &\text{for }
(i,j)=(-1,1).
  \end{cases}
\end{displaymath}
Explicitly for $(i,j)\in\set{(k,-k)}_{k\in J}$, we have
\begin{eqnarray*}
  C_{j,-j} &=& \v{j}\tensor\v{-j} + q^4 \v{-j}\tensor\v{j}
  +q(-q^2)^{n-j}\v{0}\tensor\v{0} \\
  &&\quad  -(1-q^4) \sum_{k=j+1}^n (-q^2)^{k-j} \v{-k}\tensor\v{k} \qquad
  (j\in\range{1,n}),\\
  C_{-j,j} &=& \v{-j}\tensor z^{-1}\v{j}
  +q^4z^{-1}\v{j}\tensor\v{-j}\\
  &&\quad  -(1-q^4)\sum_{k=2}^{j-1} (-q^2)^{j-k}
  z^{-1}\v{k}\tensor\v{-k}\\
  &&\quad  -(-q^2)^{j-1} (z^{-1}\v{1}\tensor\v{-1} +\v{-1}\tensor z^{-1}\v{1})
  \qquad (j\in\range{2,n}),\\
  C_{0,0} &=& \v{0}\tensor z^{-1}\v{0} +q^2z^{-1}\v{0}\tensor
  \v{0}\\
  &&\quad  +q^2[2](1-q^4) \sum_{k=2}^n (-q^2)^{n-k} z^{-1}\v{k}\tensor
  \v{-k}\\
  &&\quad +q^2[2](-q^2)^{n-1} (z^{-1}\v{1}\tensor\v{-1} +\v{-1}\tensor
  z^{-1}\v{1}),\\
  C_{-1,1} &=& \v{-1}\tensor z^{-2}\v{1}
  +q^4z^{-2}\v{1}\tensor\v{-1}\\
  &&\quad +q(-q^2)^{n-1}z^{-1}\v{0}\tensor\v{0}
  -(1-q^4)\sum_{k=2}^n(-q^2)^{k-1} z^{-1}\v{-k}\tensor\v{k}.
\end{eqnarray*}

\begin{proposition}
  Identify $C_{i,j}$ with $C_{\b{i}, z^{-H(i,j)}\b{j}}$.  Then
  $\set{z^m\tensor z^m\cdot C_{i,j}}_{m\in\Z;i,j\in J}$ with the
  function $l$ satisfy condition~(R) of
  subsection~\ref{subsec:wedge-prod}.
\end{proposition}

\subsubsection{Fock space}
For $U_q(B^{(1)}_n)$ we have
\begin{eqnarray*}
  \Bmin&=&\set{\b{1},\b{0},\b{-1}},\\
  (\Pcl^+)_1&=&\set{\Lambda_1^\cl,\Lambda_n^\cl,\Lambda_0^\cl},
\end{eqnarray*}
with
\begin{eqnarray*}
  \e(\b{1})= \Lambda_0^\cl, \quad&\quad
  \e(\b{0})= \Lambda_n^\cl, \quad&\quad
  \e(\b{-1})= \Lambda_1^\cl,\\
  \f(\b{1})= \Lambda_1^\cl, \quad&\quad
  \f(\b{0})= \Lambda_n^\cl, \quad&\quad
  \f(\b{-1})= \Lambda_0^\cl.
\end{eqnarray*}
Since $H(\b{0}\tensor\b{0})=1$, $H(\b{1}\tensor z\b{-1})=1$ and
$H(z\b{-1}\tensor\b{1})=1$, there are two ground state sequences (up
to overall shifts by $z^{k}$ ($k\in\Z$)):
\begin{eqnarray*}
  \bo_m &=& \b{0} \qquad (m\in\Z)\\
  \lambda_m &=& \Lambda_n \qquad (m\in\Z)
\end{eqnarray*}
and
\begin{eqnarray*}
  \bo_m &=&
  \begin{cases}
    \b{1} &\text{for } m\in 2\Z,\\
    z\b{-1} &\text{for } m\in2\Z+1,
  \end{cases}\\
  \lambda_m &=&
  \begin{cases}
    \Lambda_1-\frac{m}{2}\delta &\text{for } m\in 2\Z,\\
    \Lambda_0-\frac{m-1}{2}\delta &\text{for } m\in2\Z+1.
  \end{cases}
\end{eqnarray*}
The vacuum vectors are respectively
\begin{displaymath}
  \vac{m}:= \v{0}\wedge \v{0}\wedge \v{0}\wedge \v{0}\wedge
  \cdots\cdots \qquad (m\in\Z),
\end{displaymath}
with $\wt(\vac{m})=\Lambda_n$, and
\begin{displaymath}
  \vac{m}:=
  \begin{cases}
    \v{1}\wedge z\v{-1}\wedge\v{1}\wedge\cdots\cdots &\text{for } m\in 2\Z,\\
    z\v{-1}\wedge\v{1}\wedge z\v{-1}\wedge\cdots\cdots &\text{for } m\in 2\Z+1,
  \end{cases}
\end{displaymath}
with $\wt(\vac{m})= \Lambda_1-\frac{m}{2}\delta\;
(m:\text{even}),\; \Lambda_0-\frac{m-1}{2}\delta\;
(m:\text{odd})$.

\subsection{Level~1 $A^{(2)}_{2n-1}$}

\subsubsection{Cartan datum}
The Dynkin diagram for $A^{(2)}_{2n-1}$ ($n\geq 3$) is
\begin{eqnarray*}
  \dynkin{0}\link \!\!\!\! &\dynkin{2}& \!\!\! \link \dynkin{3} \link \cdots
  \cdots \link \dynkin{n-2}  \link \dynkin{n-1}
  \leftdoublelink \dynkin{n}\;.\\[-5pt]
  &|&\\[-5pt]
  &\dynkin{1}&
\end{eqnarray*}
For $A^{(2)}_{2n-1}$ we have
\begin{eqnarray*}
  \delta &=& \alpha_0+ \alpha_1+ (\sum_{i=2}^{n-1} 2\alpha_i) + \alpha_n,\\
  c &=& h_0+ h_1+ (\sum_{i=2}^n2h_i),\\
  (\alpha_i,\alpha_i) &=&
  \begin{cases}
    2 &\text{for } i\in\range{0,n-1},\\
    4 &\text{for } i=n.
  \end{cases}
\end{eqnarray*}

\subsubsection{Perfect crystal}
Let $J:=\range{-n,-1}\union\range{1,n}$.  Let $V$ be the
$(2n)$-dimensional $U'_q(A^{(2)}_{2n-1})$-module with the level~$1$
perfect crystal $B:=\set{\b{i}}_{i\in J}$ with crystal graph:
\begin{eqnarray*}
  \b{2}\: \squeeze \rightcrystal{2} \b{3} \rightcrystal{3} \cdots\cdots
  \rightcrystal{n-2} \b{n-1}\rightcrystal{n-1} \squeeze \b{n}\\[-4pt]
  {\scriptstyle 1}\!\!\nearrow\; \nwarrow\!\!{\scriptstyle 0}
  \!\!\!\!\squeeze \squeeze \,| \\[-7pt]
  \b{1}\;\;\;\;\;\b{-1} \!\!\!\!\!\!\!\! \squeeze \squeeze
  \,|{\scriptstyle n} \qquad  .\\[-6.5pt]
  {\scriptstyle 0}\!\!\nwarrow\;\nearrow\!\!{\scriptstyle 1}
  \!\!\!\!\squeeze\squeeze \makebox[1.5ex]{$\downarrow$}\\[-6pt]
  \b{-2}\squeeze \leftcrystal{2} \b{-3}
  \leftcrystal{3} \cdots\cdots \leftcrystal{n-2}
  \b{1-n}\leftcrystal{n-1} \squeeze \b{-n}
\end{eqnarray*}
The elements of $B$ have the following weights
\begin{eqnarray*}
  \wt(\b{i}) &=& \alpha_n/2 +\sum_{k\in[i,n-1]} \alpha_k
  = \Lambda_i -\Lambda_{i-1} -\delta_{i,2}\Lambda_0 \qquad
  (i\in\range{1,n}),\\
  \wt(\b{-i}) &=& -\wt(\b{i}) \qquad\qquad (i\in\range{1,n}).
\end{eqnarray*}
Let $\v{j}:=G(\b{j})$ ($j\in J$).  The action of
$U'_q(A^{(2)}_{2n-1})$ on $\v{j}\in V$ obeys~\eqref{eq:LG}.

\subsubsection{Energy function}
Define the following ordering of $J$
\begin{displaymath}
  1 \succ 2 \succ \cdots \succ n \succ -n \succ 1-n \succ \dots
  \succ -1.
\end{displaymath}
The energy function $H$ takes the following values on $B\tensor B$
\begin{eqnarray*}
  H(\b{i}\tensor \b{j})=
  \begin{cases}
    2 &\text{for } (i,j)=(-1,1),\\
    1 &\text{for } (i,j)\in\bigset{(i',j')\in J^2\setminus\set{(-1,1)} \mid
i'\prec
      j'},\\
    0 &\text{for } (i,j)\in\set{(i',j')\in J^2\mid i'\succ
      j'}\union\set{(k,k)}_{k\in J}.
  \end{cases}
\end{eqnarray*}
Write $H(i,j)$ for $H(\b{i}\tensor\b{j})$ ($i,j\in J$).

The Coxeter number of $A^{(2)}_{2n-1}$ is $h=2n-1=\dim V-1$.  We take
$l$ to be
\begin{displaymath}
  l(z^m\b{j}) =
  \begin{cases}
    hm+n-j &\text{for } j\in\range{1,n},\\
    hm-(n+1+j) &\text{for } j\in\range{-n,-1}.
  \end{cases}
\end{displaymath}
The functions $H$ and $l$ satisfy condition~(L).  Note that
$l(z^m\b{1})=l(z^{m+1}\b{-1})$ ($m\in\Z$), so $l$ gives a partial
ordering of $\Baff$.

\subsubsection{Wedge relations}
We have
\begin{displaymath}
  N:=U_q(A^{(2)}_{2n-1}){[z\tensor z, z^{-1}\tensor z^{-1},z\tensor1 +
    1\tensor z]}\cdot \v{1}\tensor\v{1}\subset \Vaff\tensor\Vaff.
\end{displaymath}
The following elements are contained in
$U_q(A^{(2)}_{2n-1})\cdot\v{1}\tensor\v{1}\subset N$:
\begin{eqnarray*}
  \tilde{C}_{i,i} &=& \v{i}\tensor\v{i}
  \qquad\qquad\qquad\qquad\qquad\qquad (i\in J),\\
  \tilde{C}_{i,-i} &=& \v{i}\tensor z^{-H(i,-i)} \v{-i} + q\v{i+1}\tensor
  z^{-H(i,-i)} \v{-i-1}\\
  &&\quad +qz^{-H(i,-i)} \v{-i-1}\tensor\v{i+1} \\
  &&\quad +q^2 z^{-H(i,-i)} \v{-i}\tensor\v{i} \qquad\qquad \left(i\in
    J\setminus\set{-1,n}\right),\\
  \tilde{C}_{i,j} &=& \v{i}\tensor z^{-H(i,j)}\v{j}\\
  &&\quad +q z^{-H(i,j)} \v{j}\tensor\v{i} \qquad \left((i,j)\in
    J^2\setminus\set{(k,k),(k,-k)}_{k\in J}\right),\\
  \tilde{C}_{n,-n} &=& \v{n}\tensor\v{-n} + q^2\v{-n}\tensor\v{n}, \\
  \tilde{C}_{-1,1} &=& \v{-1}\tensor z^{-2}\v{1} +
  qz^{-1}\v{-2}\tensor z^{-1}\v{2}\\
  &&\quad +q z^{-1}\v{2}\tensor z^{-1}\v{-2} +q^2
  z^{-2}\v{1}\tensor\v{-1}.
\end{eqnarray*}
Each $\tilde{C}_{i,j}$ has $\v{i}\tensor z^{-H(i,j)}\v{j}$
as its first term and a term in $z^{-H(i,j)}\v{j}\tensor\v{i}$.

Define the following elements in $N$.
\begin{displaymath}
  C_{i,j}:=
  \begin{cases}
    \tilde{C}_{i,j} &\text{for } (i,j)\in J^2\setminus\set{(k,-k)}_{k\in J},\\
    \sum_{k=i}^n(-q)^{k-i}\tilde{C}_{k,-k} &\text{for }
    (i,j)\in\set{(k,-k)}_{k\in\range{1,n}},\\
    \sum_{k=2}^j(-q)^{j-k}\tilde{C}_{-k,k} &\text{for }
    (i,j)\in\set{(-k,k)}_{k\in\range{2,n}},\\
    \tilde{C}_{-1,1}- q(z^{-1}\tensor z^{-1}) C_{2,-2} &\text{for }
(i,j)=(-1,1).
  \end{cases}
\end{displaymath}
Explicitly for $(i,j)\in\set{(k,-k)}_{k\in J}$, we have
\begin{eqnarray*}
  C_{j,-j} &=& \v{j}\tensor\v{-j} + q^2 \v{-j}\tensor\v{j}\\
  &&\quad  -(1-q^2) \sum_{k=j+1}^n (-q)^{k-j} \v{-k}\tensor\v{k} \qquad\qquad
  \quad(j\in\range{1,n}),\\
  C_{-j,j} &=& \v{-j}\tensor z^{-1}\v{j}
  +q^2z^{-1}\v{j}\tensor\v{-j}\\
  &&\quad  -(1-q^2)\sum_{k=2}^{j-1} (-q)^{j-k}
  z^{-1}\v{k}\tensor\v{-k}\\
  &&\quad  -(-q)^{j-1} (z^{-1}\v{1}\tensor\v{-1} +\v{-1}\tensor z^{-1}\v{1})
  \qquad (j\in\range{2,n}),\\
  C_{-1,1} &=& \v{-1}\tensor z^{-2}\v{1}
  +q^2z^{-2}\v{1}\tensor\v{-1}\\
  &&\quad -(1-q^2)\sum_{k=2}^n(-q)^{k-1} z^{-1}\v{-k}\tensor z^{-1}\v{k}.
\end{eqnarray*}

\begin{proposition}
  Identify $C_{i,j}$ with $C_{\b{i}, z^{-H(i,j)}\b{j}}$.  Then
  $\set{z^m\tensor z^m\cdot C_{i,j}}_{m\in\Z;i,j\in J}$ with the
  function $l$ satisfy condition~(R) of
  subsection~\ref{subsec:wedge-prod}.
\end{proposition}

\subsubsection{Fock space}
For $U_q(A^{(2)}_{2n-1})$ we have
\begin{eqnarray*}
  \Bmin&=&\set{\b{1},\b{-1}},\\
  (\Pcl^+)_1&=&\set{\Lambda_1^\cl,\Lambda_0^\cl},
\end{eqnarray*}
with
\begin{eqnarray*}
  \e(\b{1})= \Lambda_0^\cl, \quad&\quad
  \e(\b{-1})= \Lambda_1^\cl,\\
  \f(\b{1})= \Lambda_1^\cl, \quad&\quad
  \f(\b{-1})= \Lambda_0^\cl.
\end{eqnarray*}
Since $H(\b{1}\tensor z\b{-1})=1$ and
$H(z\b{-1}\tensor\b{1})=1$ there is one ground state sequence:
\begin{eqnarray*}
  \bo_m &=&
  \begin{cases}
    \b{1} &\text{for } m\in 2\Z,\\
    z\b{-1} &\text{for } m\in2\Z+1,
  \end{cases}\\
  \lambda_m &=&
  \begin{cases}
    \Lambda_1-\frac{m}{2}\delta &\text{for } m\in 2\Z,\\
    \Lambda_0-\frac{m-1}{2}\delta &\text{for } m\in2\Z+1.
  \end{cases}
\end{eqnarray*}
The vacuum vector of $\fock_m$ is
\begin{displaymath}
  \vac{m}:=
  \begin{cases}
    \v{1}\wedge z\v{-1}\wedge\v{1}\wedge\cdots\cdots &\text{for } m\in 2\Z,\\
    z\v{-1}\wedge\v{1}\wedge z\v{-1}\wedge\cdots\cdots &\text{for } m\in 2\Z+1,
  \end{cases}
\end{displaymath}
with $\wt(\vac{m})= \Lambda_1-\frac{m}{2}\delta\;
(m:\text{even}),\; \Lambda_0-\frac{m-1}{2}\delta\;
(m:\text{odd})$.

\subsection{Level~1 $D^{(1)}_n$}

\subsubsection{Cartan datum}
The Dynkin diagram for $D^{(1)}_n$ ($n\geq 4$) is
\begin{eqnarray*}
  \dynkin{0}\link \dynkin{2}\! & \!\!\! \link \dynkin{3} \link \cdots
  \cdots \link \dynkin{n-3} \link &\!\!\!\!\! \dynkin{n-2}  \link
  \dynkin{n}\;.\\[-5pt]
  |\!&&\!\! |\\[-7pt]
  \dynkin{1}\!&&\!\!\!\!\! \dynkin{n-1}
\end{eqnarray*}
For $D^{(1)}_n$ we have
\begin{eqnarray*}
  \delta &=& \alpha_0+ \alpha_1 +(\sum_{i=2}^{n-2} 2\alpha_i)
  +\alpha_{n-1} +\alpha_n,\\
  c &=& h_0 +h_1 +(\sum_{i=2}^{n-2} 2h_i) +h_{n-1} +h_n,\\
  (\alpha_i,\alpha_i) &=& 2 \qquad (i\in I).
\end{eqnarray*}

\subsubsection{Perfect crystal}
Let $J:=\range{-n,-1}\union\range{1,n}$.  Let $V$ be the
$(2n)$-dimensional $U'_q(D^{(1)}_n)$-module with the level~$1$ perfect
crystal $B:=\set{\b{i}}_{i\in J}$ with crystal graph:
\begin{eqnarray*}
  \b{2}\: \squeeze \rightcrystal{2} \b{3} \rightcrystal{3} \cdots\cdots
  \rightcrystal{n-3} \b{n-2}\rightcrystal{n-2} \squeeze \b{n-1}\\[-5pt]
  {\scriptstyle 1}\!\!\nearrow\; \nwarrow\!\!{\scriptstyle 0}
  \!\!\!\!\squeeze \squeeze \!\!\!\!\!\!\!\!\!\! {\scriptstyle
    n-1}\!\!\swarrow\;  \searrow\!\!{\scriptstyle n}\\[-5pt]
  \b{1}\;\;\;\;\;\b{-1} \!\!\!\!\!\!\!\! \squeeze \squeeze
  \!\!\!\!\! \b{n}\;\;\;\;\; \b{-n} \qquad  .\\[-5pt]
  {\scriptstyle 0}\!\!\nwarrow\; \nearrow\!\!{\scriptstyle 1}
  \!\!\!\!\squeeze\squeeze\!\!\!\! {\scriptstyle n}\!\!\searrow\;
  \swarrow \!\!{\scriptstyle n-1}\\[-5pt]
  \b{-2}\squeeze \leftcrystal{2} \b{-3}
  \leftcrystal{3} \cdots\cdots \leftcrystal{n-3}
  \b{2-n}\leftcrystal{n-2} \squeeze \b{1-n}
\end{eqnarray*}
The elements of $B$ have the following weights
\begin{eqnarray*}
  \wt(\b{i}) &=& (\sum_{k=i}^{n-2} \alpha_k) +(\alpha_{n-1}
  +\alpha_n)/2\\
  &=& \Lambda_i -\Lambda_{i-1} +\delta_{i,n-1}\Lambda_n
  -\delta_{i,2}\Lambda_0\qquad (i\in\range{1,n}),\\
  \wt(\b{-i}) &=& -\wt(\b{i}) \qquad\qquad\qquad (i\in\range{1,n}).
\end{eqnarray*}
Let $\v{j}:=G(\b{j})$ ($j\in J$).  The action of
$U'_q(D^{(1)}_n)$ on $\v{j}\in V$ obeys~\eqref{eq:LG}.

\subsubsection{Energy function}
Define the following ordering of $J$
\begin{displaymath}
  1 \succ 2 \succ \cdots \succ n \succ -n \succ 1-n \succ \dots
  \succ -1.
\end{displaymath}
The energy function $H$ takes the following values on $B\tensor B$
\begin{eqnarray*}
  H(\b{i}\tensor\b{j})=
  \begin{cases}
    2 &\text{for } (i,j)=(-1,1),\\
    1 &\text{for } (i,j)\in\bigset{(i',j')\in J^2\setminus\set{(-1,1)} \mid
i'\prec
      j'}\union \set{(n,-n)},\\
    0 &\text{for } (i,j)\in\set{(i',j')\in J^2\setminus\set{(n,-n)} \mid
i'\succ
      j'}\union\set{(k,k)}_{k\in J}.
  \end{cases}
\end{eqnarray*}
Write $H(i,j)$ for $H(\b{i}\tensor\b{j})$ ($i,j\in J$).

The Coxeter number of $D^{(1)}_n$ is $h=n+1=\dim V-2$.  We take
$l$ to be
\begin{displaymath}
  l(z^m\b{j}) =
  \begin{cases}
    hm+n-j &\text{for } j\in\range{1,n},\\
    hm-(n+j) &\text{for } j\in\range{-n,-1}.
  \end{cases}
\end{displaymath}
The functions $H$ and $l$ satisfy condition~(L).  Note that
$l(z^m\b{1})=l(z^{m+1}\b{-1})$ and $l(z^m\b{n})= l(z^m\b{-n})$
($m\in\Z$), so $l$ gives a partial ordering of $\Baff$.

\subsubsection{Wedge relations}
We have
\begin{displaymath}
  N:=U_q(D^{(1)}_n){[z\tensor z, z^{-1}\tensor z^{-1},z\tensor1 +
    1\tensor z]}\cdot \v{1}\tensor\v{1}\subset \Vaff\tensor\Vaff.
\end{displaymath}
The following elements are contained in
$U_q(D^{(1)}_n)\cdot\v{1}\tensor\v{1}\subset N$:
\begin{eqnarray*}
  \tilde{C}_{i,i} &=& \v{i}\tensor\v{i}
  \qquad\qquad\qquad\qquad\qquad\qquad (i\in J),\\
  \tilde{C}_{i,-i} &=& \v{i}\tensor z^{-H(i,-i)} \v{-i} + q\v{i+1}\tensor
  z^{-H(i,-i)} \v{-i-1}\\
  &&\quad +qz^{-H(i,-i)} \v{-i-1}\tensor\v{i+1} \\
  &&\quad +q^2 z^{-H(i,-i)} \v{-i}\tensor\v{i} \qquad\qquad \left(i\in
    J\setminus\set{-1,n}\right),\\
  \tilde{C}_{i,j} &=& \v{i}\tensor z^{-H(i,j)}\v{j}\\
  &&\quad +q z^{-H(i,j)} \v{j}\tensor\v{i} \qquad \left((i,j)\in
    J^2\setminus\set{(k,k),(k,-k)}_{k\in J}\right),\\
  \tilde{C}_{n,-n} &=& \v{n}\tensor z^{-1}\v{-n}
  +q\v{1-n}\tensor z^{-1}\v{n-1}\\
  &&\quad +qz^{-1}\v{n-1}\tensor\v{1-n}
  +q^2z^{-1}\v{-n}\tensor\v{n},\\
  \tilde{C}_{-1,1} &=& \v{-1}\tensor z^{-2}\v{1} +qz^{-1}\v{-2}\tensor
  z^{-1}\v{2}\\
  &&\quad +q z^{-1}\v{2}\tensor z^{-1}\v{-2}
  +q^2z^{-2}\v{1}\tensor\v{-1}.
\end{eqnarray*}
Each $\tilde{C}_{i,j}$ has $\v{i}\tensor z^{-H(i,j)}\v{j}$
as its first term and a term in $z^{-H(i,j)}\v{j}\tensor\v{i}$.

Define the following elements in $N$.
\begin{displaymath}
  C_{i,j}:=
  \begin{cases}
    \tilde{C}_{i,j} &\text{for } (i,j)\in J^2\setminus\set{(k,-k)}_{k\in J},\\
    \sum_{k=i}^{n-1}(-q)^{k-i}\tilde{C}_{k,-k} &\text{for }
    (i,j)\in\set{(k,-k)}_{k\in\range{1,n-1}},\\
    \sum_{k=2}^j(-q)^{j-k}\tilde{C}_{-k,k} &\text{for }
    (i,j)\in\set{(-k,k)}_{k\in\range{2,n}},\\
    \tilde{C}_{n,-n}-q C_{1-n,n-1} &\text{for } (i,j)=(n,-n),\\
    \tilde{C}_{-1,1}- q(z^{-1}\tensor z^{-1}) C_{2,-2} &\text{for }
(i,j)=(-1,1).
  \end{cases}
\end{displaymath}
Explicitly for $(i,j)\in\set{(k,-k)}_{k\in J}$, we have
\begin{eqnarray*}
  C_{j,-j} &=& \v{j}\tensor\v{-j} + q^2 \v{-j}\tensor\v{j}\\
  &&\quad  -(1-q^2) \sum_{k=j+1}^{n-1} (-q)^{k-j} \v{-k}\tensor\v{k}\\
  &&\quad  -(-q)^{n-j}(\v{n}\tensor\v{-n} +\v{-n}\tensor\v{n})
  \qquad\qquad\quad (j\in\range{1,n}),\\
  C_{-j,j} &=& \v{-j}\tensor z^{-1}\v{j}
  +q^2z^{-1}\v{j}\tensor\v{-j}\\
  &&\quad  -(1-q^2)\sum_{k=2}^{j-1} (-q)^{j-k}
  z^{-1}\v{k}\tensor\v{-k}\\
  &&\quad  -(-q)^{j-1} (z^{-1}\v{1}\tensor\v{-1} +\v{-1}\tensor z^{-1}\v{1})
  \qquad\qquad\quad (j\in\range{2,n}),\\
  C_{n,-n} &=& \v{n}\tensor z^{-1}\v{-n} + q^2z^{-1}\v{-n}\tensor\v{n}\\
  &&\quad -(1-q^2)\sum_{k=2}^{n-1} (-q)^{n-k}
  z^{-1}\v{k}\tensor\v{-k}\\
  &&\quad -(-q)^{n-1}(z^{-1}\v{1}\tensor\v{-1} +\v{-1}\tensor z^{-1}\v{1})\\
  C_{-1,1} &=& \v{-1}\tensor z^{-2}\v{1}
  +q^2z^{-2}\v{1}\tensor\v{-1}\\
  &&\quad -(1-q^2)\sum_{k=2}^{n-1}(-q)^{k-1}
  z^{-1}\v{-k}\tensor z^{-1}\v{k}\\
  &&\quad -(-q)^{n-1}(z^{-1}\v{n}\tensor z^{-1}\v{-n} +z^{-1}\v{-n}\tensor
  z^{-1}\v{n}).
\end{eqnarray*}

\begin{proposition}
  Identify $C_{i,j}$ with $C_{\b{i}, z^{-H(i,j)}\b{j}}$.  Then
  $\set{z^m\tensor z^m\cdot C_{i,j}}_{m\in\Z;i,j\in J}$ with the
  function $l$ satisfy condition~(R) of
  subsection~\ref{subsec:wedge-prod}.
\end{proposition}

\subsubsection{Fock space}
For $U_q(D^{(1)}_n)$ we have
\begin{eqnarray*}
  \Bmin&=&\set{\b{1},\b{-1},\b{n},\b{-n}},\\
  (\Pcl^+)_1&=&\set{\Lambda_1^\cl, \Lambda_0^\cl, \Lambda_{n-1}^\cl,
  \Lambda_n^\cl},
\end{eqnarray*}
with
\begin{eqnarray*}
  \e(\b{1})= \Lambda_0^\cl, \quad&\quad
  \f(\b{1})= \Lambda_1^\cl,\\
  \e(\b{-1})= \Lambda_1^\cl, \quad&\quad
  \f(\b{-1})= \Lambda_0^\cl,\\
  \e(\b{n})= \Lambda_n^\cl, \quad&\quad
  \f(\b{n})= \Lambda_{n-1}^\cl,\\
  \e(\b{-n})= \Lambda_{n-1}^\cl, \quad&\quad
  \f(\b{-n})= \Lambda_n^\cl.
\end{eqnarray*}
Since $H(\b{1}\tensor z\b{-1})=1$, $H(z\b{-1}\tensor\b{1})=1$,
$H(\b{n}\tensor\b{-n})=1$ and $H(\b{-n}\tensor\b{n})=1$,
there are two ground state sequences:
\begin{eqnarray*}
  \bo_m &=&
  \begin{cases}
    \b{1} &\text{for } m\in 2\Z,\\
    z\b{-1} &\text{for } m\in2\Z+1,
  \end{cases}\\
  \lambda_m &=&
  \begin{cases}
    \Lambda_1-\frac{m}{2}\delta &\text{for } m\in 2\Z,\\
    \Lambda_0-\frac{m-1}{2}\delta &\text{for } m\in2\Z+1,
  \end{cases}
\end{eqnarray*}
and
\begin{eqnarray*}
  \bo_m &=&
  \begin{cases}
    \b{n} &\text{for } m\in 2\Z,\\
    \b{-n} &\text{for } m\in2\Z+1,
  \end{cases}\\
  \lambda_m &=&
  \begin{cases}
    \Lambda_{n-1} &\text{for } m\in 2\Z,\\
    \Lambda_n &\text{for } m\in2\Z+1.
  \end{cases}
\end{eqnarray*}

The vacuum vector of $\fock_m$ are respectively
\begin{displaymath}
  \vac{m}:=
  \begin{cases}
    \v{1}\wedge z\v{-1}\wedge\v{1}\wedge\cdots\cdots &\text{for } m\in 2\Z,\\
    z\v{-1}\wedge\v{1}\wedge z\v{-1}\wedge\cdots\cdots &\text{for } m\in 2\Z+1,
  \end{cases}
\end{displaymath}
with $\wt(\vac{m})= \Lambda_1-\frac{m}{2}\delta\;
(m:\text{even}),\; \Lambda_0-\frac{m-1}{2}\delta\;
(m:\text{odd})$, and
\begin{displaymath}
  \vac{m}:=
  \begin{cases}
    \v{n}\wedge \v{-n}\wedge\v{n}\wedge\cdots\cdots &\text{for } m\in 2\Z,\\
    \v{-n}\wedge\v{n}\wedge\v{-n}\wedge\cdots\cdots &\text{for } m\in 2\Z+1,
  \end{cases}
\end{displaymath}
with $\wt(\vac{m})= \Lambda_{n-1}\; (m:\text{even}),\;
\Lambda_n\; (m:\text{odd})$.

\subsection{Level~1 $D^{(2)}_{n+1}$}

\subsubsection{Cartan datum}
The Dynkin diagram for $D^{(2)}_{n+1}$ ($n\geq 4$) is
\begin{displaymath}
  \dynkin{0}\leftdoublelink \dynkin{1} \link \dynkin{2} \link \cdots
  \cdots  \link \dynkin{n-2} \link \dynkin{n-1}
  \rightdoublelink \dynkin{n}\;.
\end{displaymath}
For $D^{(2)}_{n+1}$ we have
\begin{eqnarray*}
  \delta &=& \sum_{i\in I}\alpha_i,\\
  c &=& h_0 +(\sum_{i=1}^{n-1} 2h_i) +h_n,\\
  (\alpha_i,\alpha_i) &=&
  \begin{cases}
    2 &\text{for } i\in\set{0,n},\\
    4 &\text{for } i\in I\setminus\set{0,n}.
  \end{cases}
\end{eqnarray*}

\subsubsection{Perfect crystal}
Let $J:=\range{-n,n}\union\set{\phi}$.  Let $V$ be the
$(2n+2)$-dimensional $U'_q(D^{(2)}_{n+1})$-module with the level~$1$
perfect crystal $B:=\set{\b{i}}_{i\in J}$ and crystal graph:
\begin{eqnarray*}
  \b{1} \squeeze \rightcrystal{1} \b{2} \rightcrystal{2} \cdots\cdots
  \rightcrystal{n-2} \b{n-1}\rightcrystal{n-1} \squeeze \b{n}\\[-5pt]
  {\scriptstyle 0}\!\!\uparrow \squeeze \squeeze\downarrow\!\!{\scriptstyle
    n}\\[-5pt]
  \b{\phi}\squeeze \squeeze \b{0} \qquad .\\[-5pt]
  {\scriptstyle 0}\!\!\uparrow \squeeze \squeeze \downarrow\!\!{\scriptstyle
  n}\\[-5pt]
  \b{-1}\squeeze \leftcrystal{1} \b{-2}
  \leftcrystal{2} \cdots\cdots \leftcrystal{n-2}
  \b{1-n}\leftcrystal{n-1} \squeeze \b{-n}
\end{eqnarray*}
The elements of $B$ have the following weights
\begin{eqnarray*}
  \wt(\b{i}) &=& \sum_{k=i}^{n} \alpha_k
  = (1+\delta_{i,n})\Lambda_i-(1+\delta_{i,1})\Lambda_{i-1} \qquad
  (i\in\range{1,n})\\
  \wt(\b{0}) &=& 0\\
  \wt(\b{\phi}) &=& 0\\
  \wt(\b{-i}) &=& -\wt(\b{i}) \qquad (i\in\range{1,n}).
\end{eqnarray*}
Let $\v{j}:=G(\b{j})$ ($j\in J$).  The action of
$U'_q(D^{(2)}_{n+1})$ on $\v{j}\in V$ obeys~\eqref{eq:LG}.

Let $J_0:=J\setminus\set{\phi}$.  Let $V_0$ denote that subspace of
$V$ spanned by $\set{\v{j}}_{j\in J_0}$.  Then, $\Vaff$ decomposes into
two $U_q(D^{(2)}_{n+1})$-modules:
\begin{eqnarray*}
  \Vaff &=&\bigl(V_0\tensor\C{[z^2,z^{-2}]}
    + \v{\phi}\tensor z\C{[z^2,z^{-2}]}\bigr)\\
  &&\quad \oplus \bigl( V_0\tensor z\C{[z^2,z^{-2}]}
    + \v{\phi}\tensor\C[z^2,z^{-2}]\bigr).
\end{eqnarray*}

\subsubsection{Energy function}
Define the following ordering of $J$
\begin{equation} \label{eqn:totord}
  1 \succ 2 \succ \cdots \succ n \succ 0 \succ -n \succ 1-n \succ \dots
  \succ -1\succ \phi.
\end{equation}
The energy function $H$ takes the following values on $B\tensor B$
\begin{eqnarray*}
  H(\b{i}\tensor\b{j})=
  \begin{cases}
    2 &\text{for } (i,j)\in\set{(i',j')\in J_0^2\mid i'\prec j'}\union
    \set{(0,0),(\phi,\phi)},\\
    1 &\text{for } (i,j)\in\set{(k,\phi),(\phi,k)\in J^2\mid k\in
    J\setminus\set{\phi}},\\
    0 &\text{for } (i,j)\in\set{(i',j')\in J_0^2\mid i'\succ j'}\union
    \set{(k,k)}_{k\in J\setminus\set{0,\phi}}.
  \end{cases}
\end{eqnarray*}
Write $H(i,j)$ for $H(\b{i}\tensor\b{j})$ ($i,j\in J$).

The Coxeter number of $D^{(2)}_{n+1}$ is $h=n+1=\frac{1}{2}\dim V$.
We take $l$ to be
\begin{displaymath}
  l(z^m\b{j}) =
  \begin{cases}
    hm+n+1-j &\text{for } j\in\range{1,n},\\
    hm &\text{for } j\in\set{0,\phi},\\
    hm-(n+1+j) &\text{for } j\in\range{-n,-1}.
  \end{cases}
\end{displaymath}
The functions $H$ and $l$ satisfy condition~(L).  Note that
$l(z^m\b{0})=l(z^m\b{\phi})$ and $l(z^m\b{i})=l(z^{m+1}\b{i-h})$
($m\in\Z$ and $i\in\range{1,n}$), so $l$ gives a partial ordering of
$\Baff$.  ($l$ gives a total ordering of each of the crystals of the
two irreducible submodules.)

\subsubsection{Wedge relations}
In $\Vaff\tensor\Vaff$ we have
\begin{displaymath}
  N:=U_q(D^{(2)}_{n+1}){[z\tensor z, z^{-1}\tensor z^{-1},z\tensor1 +
    1\tensor z]}\cdot \v{1}\tensor\v{1}.
\end{displaymath}
The following elements are contained in
$U_q(D^{(2)}_{n+1}){[z\tensor z,z^{-1}\tensor
  z^{-1}]}\cdot\v{1}\tensor\v{1}\subset N$:
\begin{eqnarray*}
  \tilde{C}_{i,i} &=& \v{i}\tensor\v{i}
  \qquad\qquad\qquad\qquad\qquad\qquad
  (i\in J\setminus\set{0,\phi}),\\
  \tilde{C}_{i,-i} &=& \v{i}\tensor z^{-H(i,-i)} \v{-i} + q^2 \v{i+1}\tensor
  z^{-H(i,-i)} \v{-i-1}\\
  && \quad +q^2z^{-H(i,-i)}\v{-i-1}\tensor\v{i+1}\\
  && \quad + q^4 z^{-H(i,-i)} \v{-i}\tensor\v{i}\qquad\qquad \bigl(i\in
    J\setminus\set{-1,0,\phi,n}\bigr),\\
  \tilde{C}_{i,j} &=& \v{i}\tensor z^{-H(i,j)}\v{j}
  +q^2z^{-H(i,j)}\v{j}\tensor\v{i}\\
  && \quad \qquad \bigl((i,j)\in
    J^2\setminus\set{(k,k),(k,-k)}_{k\in J}\bigr)
    ,\\
  \tilde{C}_{0,0} &=& \v{0}\tensor z^{-2} \v{0} + q^2[2]\v{-n}\tensor
  z^{-2} \v{n} \\
  && \quad +q^2[2]z^{-2}\v{n}\tensor\v{-n} + q^2
  z^{-2}\v{0}\tensor\v{0}, \\
  \tilde{C}_{n,-n} &=& \v{n}\tensor\v{-n} + q\v{0}\tensor\v{0} + q^4
  \v{-n}\tensor\v{n}, \\
  \tilde{C}_{-1,1} &=& \v{-1}\tensor z^{-2} \v{1}
  +qz^{-1}\v{\phi}\tensor z^{-1}\v{\phi}
  +q^4z^{-2}\v{1}\tensor\v{-1},\\
  \tilde{C}_{\phi,\phi} &=& \v{\phi}\tensor z^{-2}\v{\phi}
  +q^2[2]z^{-1}\v{1}\tensor z^{-1}\v{-1}\\
  && \quad +q^2[2] z^{-1}\v{-1}\tensor z^{-1}\v{1}
  +q^2z^{-2}\v{\phi}\tensor\v{\phi}.
\end{eqnarray*}
Notice that each $\tilde{C}_{i,j}$ has $\v{i}\tensor z^{-H(i,j)}\v{j}$
as its first term and a term in $z^{-H(i,j)}\v{j}\tensor\v{i}$.

Define the following elements in $N$.
\begin{displaymath}
  C_{i,j}:=
  \begin{cases}
    \tilde{C}_{i,j} &\text{for } (i,j)\in J^2\setminus\set{(k,-k)}_{k\in J},\\
    \sum_{k=i}^n(-q^2)^{k-i}\tilde{C}_{k,-k} &\text{for }
    (i,j)\in\set{(k,-k)}_{k\in\range{1,n}},\\
    \sum_{k=1}^j(-q^2)^{j-k}\tilde{C}_{-k,k} &\text{for }
    (i,j)\in\set{(-k,k)}_{k\in\range{1,n}},\\
    \tilde{C}_{0,0} -q^2[2]C_{-n,n} &\text{for } (i,j)=(0,0),\\
    \tilde{C}_{\phi,\phi} -q^2[2](z^{-1}\tensor z^{-1})C_{1,-1} &\text{for }
    (i,j)=(\phi,\phi).
  \end{cases}
\end{displaymath}
Explicitly for $(i,j)\in\set{(k,-k)}_{k\in J}\union\set{\phi}$, we have
\begin{eqnarray*}
  C_{j,-j} &=& \v{j}\tensor\v{-j} + q^4 \v{-j}\tensor\v{j}
  +q(-q^2)^{n-j}\v{0}\tensor\v{0} \\
  &&\quad -(1-q^4) \sum_{k=j+1}^n (-q^2)^{k-j} \v{-k}\tensor\v{k} \qquad
  (j\in\range{1,n}),\\
  C_{-j,j} &=& \v{-j}\tensor z^{-2} \v{j} + q^4 z^{-2}\v{j}\tensor
  \v{-j} +q(-q^2)^{j-1}z^{-1}\v{\phi}\tensor z^{-1}\v{\phi}\\
  &&\quad -(1-q^4)\sum_{k=1}^{j-1} (-q^2)^{j-k} z^{-2}\v{k}\tensor\v{-k}
  \qquad (j\in\range{1,n}),\\
  C_{0,0} &=& \v{0}\tensor z^{-2}\v{0} +q^2z^{-2}\v{0}\tensor\v{0}
  +q[2](-q^2)^nz^{-1}\v{\phi}\tensor z^{-1}\v{\phi}\\
  &&\quad -[2](1-q^4)\sum_{k=1}^n (-q^2)^{n+1-k} z^{-2}\v{k}\tensor
  \v{-k},\\
  C_{\phi,\phi} &=& \v{\phi}\tensor z^{-2} \v{\phi} + q^2
  z^{-2}\v{\phi}\tensor\v{\phi} + q[2](-q^2)^nz^{-1}\v{0}\tensor
  z^{-1}\v{0}\\
  && -[2](1-q^4) \sum_{k=1}^n (-q^2)^j z^{-1}\v{-j}\tensor
  z^{-1}\v{j}.
\end{eqnarray*}

\begin{proposition}
  Identify $C_{i,j}$ with $C_{\b{i}, z^{-H(i,j)}\b{j}}$.  Then
  $\set{z^m\tensor z^m\cdot C_{i,j}}_{m\in\Z;i,j\in J}$ with the
  function $l$ satisfy condition~(R) of
  subsection~\ref{subsec:wedge-prod}.
\end{proposition}

\subsubsection{Fock space}
For $U_q(D^{(2)}_{n+1})$ we have
\begin{eqnarray*}
  \Bmin&=&\set{\b{0},\b{\phi}},\\
  (\Pcl^+)_1&=&\set{\Lambda_n^\cl,\Lambda_0^\cl},
\end{eqnarray*}
with
\begin{eqnarray*}
  \e(\b{0})= \Lambda_n^\cl, \quad&\quad
  \e(\b{\phi})= \Lambda_0^\cl,\\
  \f(\b{0})= \Lambda_n^\cl, \quad&\quad
  \f(\b{\phi})= \Lambda_0^\cl.
\end{eqnarray*}
Since $H(\b{0}\tensor z^{-1}\b{0})=1$ and $H(\b{\phi}\tensor
z^{-1}\b{\phi})=1$, there are two ground state sequences (see also the
remark at the end of \S\ref{subsec:d2-2pt}):
\begin{eqnarray*}
  \bo_m &=& z^{-m}\b{0} \qquad (m\in\Z),\\
  \cl(\lambda_m) &=& \Lambda_n \qquad (m\in\Z),
\end{eqnarray*}
and
\begin{eqnarray*}
  \bo_m &=& z^{-m}\b{\phi} \qquad (m\in\Z),\\
  \cl(\lambda_m) &=& \Lambda_0 \qquad (m\in\Z).
\end{eqnarray*}
The vacuum vectors of $\fock_m$ are respectively
\begin{displaymath}
  \vac{m}:= z^{-m}\v{0}\wedge z^{-m-1}\v{0}\wedge z^{-m-2}\v{0}\wedge
  \cdots\cdots \qquad (m\in\Z),
\end{displaymath}
with $\wt(\vac{m})=\Lambda_n$, and
\begin{displaymath}
  \vac{m}:= z^{-m}\v{\phi}\wedge z^{-m-1}\v{\phi}\wedge z^{-m-2}\v{\phi}\wedge
  \cdots\cdots \qquad (m\in\Z),
\end{displaymath}
with $\wt(\vac{m})= \Lambda_0$.



\newsection{Level~1 two point functions}

In this section we calculate the boson commutation relations using the
decomposition of the Fock space vertex operator into a product of a
$U_q(\g)$-vertex operator and a bosonic vertex operator
(Theorem~\ref{TEIRI}), for level~1 $A^{(2)}_{2n}$, $B^{(1)}_n$,
$A^{(2)}_{2n-1}$, $D^{(1)}_n$ and $D^{(2)}_{n+1}$.  The two point
functions of the level~$1$ $U_q(\g)$-vertex operators that we use are
due to Date and Okado~\cite{DO} (except for type $D^{(2)}_{n+1}$ which
is given in Appendix~\ref{app:c}).

\subsection{Summary} \label{sec6:summary}
In the following table we list the dual Coxeter number
$\dcox:=\sum_{i\in I}a^\vee_i$,
$p:=q^{(\alpha_0,\alpha_0)/(2a^\vee_0)}$ and $\xi:=(-)^{r-1}p^{\dcox}$
for $\g=X^{(r)}_n$ of types $A^{(1)}_n$, $A^{(2)}_{2n}$, $B^{(1)}_n$,
$A^{(2)}_{2n-1}$ and $D^{(1)}_n$.
\begin{center}
  \begin{tabular}{|c||c|c|c|c|c|c|}
    \hline $\g$ & $A^{(1)}_n$ & $A^{(2)}_{2n}$ & $B^{(1)}_n$ &
    $A^{(2)}_{2n-1}$ & $D^{(1)}_n$ \\ \hline\hline
    $\dcox$ & $n+1$ & $2n+1$ & $2n-1$ & $2n$ & $2n-2$ \\ \hline
    $p$ & $q$ & $q^2$ & $q^2$ & $q$ & $q$ \\ \hline
    $\xi$ & $q^{n+1}$ & $-q^{2(2n+1)}$ & $q^{2(2n-1)}$ & $-q^{2n}$ &
    $q^{2n-2}$ \\ \hline
  \end{tabular}
\end{center}

\begin{proposition}[\cite{KMS}]
  For $A^{(1)}_n$ at level~1, we have
  \begin{displaymath}
    \gamma_m= m\frac{1-\xi^{2m}}{1-q^{2m}}.
  \end{displaymath}
\end{proposition}
See~\cite{KMS}~\S2 for the proof.
\medskip

Let $\g=X^{(r)}_n$ be of type $A^{(2)}_{2n}$, $B^{(1)}_n$,
$A^{(2)}_{2n-1}$ or $D^{(1)}_n$. Let $\vac{\g}$ be
one of the vacuum vectors of the $U_q(\g)$-Fock modules described in
the previous section. For each type, direct calculations of
$B_m\cdot B_{-m}\cdot\vac{\g}$ for small $m$, suggest the following
result.

\begin{theorem} \label{thm:boson-comm}
  For $\g=X^{(r)}_n\in \set{A^{(2)}_{2n}, B^{(1)}_n,
  A^{(2)}_{2n-1}, D^{(1)}_n}$ at level~1, we have
  \begin{displaymath}
    \gamma_m = m\frac{1+ \xi^m}{1-p^{2m}}.
  \end{displaymath}
\end{theorem}
In this section we prove this theorem case by case using
Proposition~\ref{prop:2pt-decomp}. We also give a corresponding result
for level~1 $D^{(2)}_{n+1}$.

For this boson commutation relation, the boson two point
function~\eqref{eq:boson-2pt} is
\begin{equation}\label{eq:boson-2pt-level1-g}
  \theta(w_2/w_1)=
  \frac{(w_2/w_1;\xi^2)_\infty
    (p^2\xi w_2/w_1;\xi^2)_{\infty}}
  {(p^2 w_2/w_1;\xi^2)_\infty (\xi w_2/w_1;\xi^2)_\infty}.
\end{equation}

\medskip
Let us introduce the operator $Z(t,d)\in\End(\Vaff\tensor\Vaff)$
defined by:
\begin{displaymath}
  Z(t,d):= z^t\tensor z^{d-t} +\delta(2t>d)z^{d-t}\tensor z^t
  -\delta(2t<d)z^t\tensor z^{d-t}\qquad (t,d\in\Z).
\end{displaymath}
Note that $Z(t,d)$ is a symmetric Laurent polynomial in
  $z\tensor1$ and $1\tensor z$, so we have
\begin{lemma}
  $Z(t,d)\cdot N\subset N$ ($t,d\in\Z$).
\end{lemma}

\subsection{Type $A^{(2)}_{2n}$}
Recall we have $\lambda_m=\Lambda_n$ and $\bo_m=\b{0}$ ($m\in\Z$).
So the level~1 intertwiner maps
\begin{displaymath}
  \Phi_m:\Vaff\tensor V(\Lambda_n)\mapto V(\Lambda_n) \qquad
  (m\in\Z).
\end{displaymath}
{}From~\cite{DO}, (up to a factor of a constant power in $w_2/w_1$) we
have in our notation~\eqref{eq:vertex-2pt}
\begin{eqnarray*}
  \phi_{\v{0},\v{0}}(w_2/w_1) &=& (1- p^{\dcox+1} w_2/w_1)\\
  &&\quad\times\frac{(p^{2\dcox+2}w_2/w_1;p^{2\dcox})_\infty
    (-p^{3\dcox}w_2/w_1;p^{2\dcox})_\infty}
  {(-p^{\dcox+2} w_2/w_1;p^{2\dcox})_\infty
    (p^{2\dcox}w_2/w_1; p^{2\dcox})_\infty}.
\end{eqnarray*}

Define
\begin{displaymath}
  g_j(t):= \dvac{m-1}z^{t}\v{j}\wedge z^{-t}\v{-j}\wedge\vac{m+1}
  \qquad (j\in J;t\in\Z).
\end{displaymath}
Note that $g_0(0)=1$, $g_{-j}(0)=0$ ($j\in\range{1,n}$) and $g_j(t)=0$
($j\in J$; $t\in\Z_{<0}$).

\begin{proposition} \label{prop:A2-recurrence}
  $g_0(t)$ satisfies the following recurrence relation
  \begin{eqnarray} \label{eq:A2even-recurrence}
    g_0(t) -(p^2-p^\dcox)g_0(t-1)-p^{\dcox+2}g_0(t-2)=\qquad \notag\\
    \delta_{t,0} -(1+p^{\dcox+1})\delta_{t,1}
    +p^{\dcox+1}\delta_{t,2}.
  \end{eqnarray}
\end{proposition}

\begin{pf}
  The proof for $n>1$ goes as follows (the exceptional case $n=1$ is
  similar).
  First note that any element in $N$, that is generated by $C_{j,-j}$
  ($j\in J$), gives rise to a linear relation of some $g_k(t)$
  ($k\in J$, $t\in\Z$).  For example $(z\tensor1 +1\tensor z)\cdot
  C_{0,0}$
  gives
  \begin{displaymath} 
    g_0(1) +q^2g_0(0) +q^2[2](1-q^4)\sum_{k=1}^n (-q^2)^{n-k}g_k(0)
    +g_0(0) =0.
  \end{displaymath}
  From $C_{k,k}$ ($k\in \range{1,n}$) we get
  \begin{displaymath}
    g_k(0) +q(-q^2)^{n-k}g_0(0) =0.
  \end{displaymath}
  Combining these two relations, we get
  \begin{displaymath}
    g_0(1) + (1 -q^4 +q^{4n+2} +q^{4n+4})g_0(0) =0,
  \end{displaymath}
  which is~\eqref{eq:A2even-recurrence} with $t=1$.

  The recurrence relation in the general case ($t\in\N$) comes from
  \begin{align*}
    \A_t = &\bigl(Z(t,1)
    -p^{\dcox+1}Z(t-1,1) \bigr)\cdot C_{0,0}\\
    &\begin{aligned}
        + {[2]}(1-p^2)(-p)^{n+1} \sum_{j=1}^n \bigl((-p)^{-j}&Z(t-1,0)\cdot
    C_{j,-j}\\ -(-p)^j &Z(t-1,1)\cdot
    C_{-j,j} \bigr).
    \end{aligned}
  \end{align*}
  Let $\A_t^\wedge$ denote the image of $\A_t$ in $\Vaff\wedge \Vaff$.
  We have:
  \begin{multline*}
    \dvac{m-1}\A_t^\wedge\wedge\vac{m+1} =\\
  \begin{aligned}
    &g_0(t) -(p^{\dcox+1}-p)g_0(t-1)-p^{\dcox+2}g_0(t-2)\\
    &
    \begin{aligned}
      +{[2]}(1-p^2) \Bigl\{&-\sum_{j=1}^n(-p)^{n+1-j}\bigl(g_j(t-1)
      -(-p)^{\dcox+1}g_j(t-2) \bigr)\\
      & +\sum_{j=1}^n (-p)^{n+1-j}\bigl(g_j(t-1) +
      p^2g_{-j}(t-1)\bigr)\\
      & -q \frac{(-p)^\dcox-(-p)}{1-p^2}g_0(t-1)\\
      & -(1-p^2)p^{n+1}\sum_{j=1}^n (-)^j{[j-1]}_p g_{-j}(t-1)\\
      & -\sum_{j=1}^n (-p)^{n+1+j}\bigl(g_{-j}(t-1)
      +p^2g_j(t-2)\bigr)\\
      & -(1-p^2)p^{\dcox+1}\sum_{j=1}^n (-)^{n-j}{[n-j]}_p
      g_j(t-2)\Bigr\}
    \end{aligned}\\
    & +\delta(2t>1)g_0(1-t)-\delta(2t<1)g_0(t)\\
    & -p^{\dcox+1}\bigl(\delta(2t>3)g_0(2-t) -\delta(2t<3)g_0(t-1)\bigr)=0.
  \end{aligned}
\end{multline*}
All terms in $g_k$ ($k\in J\setminus\set{0}$) cancel
and the proposition follows.
\end{pf}

The two point function of Fock intertwiners~\eqref{eq:fock-2pt} is
given by
\begin{corollary}
  \begin{displaymath}
    \omega_{\v{0},\v{0}}(w_2/w_1)=
    \frac{(1-w_2/w_1)(1-p^{\dcox+1} w_2/w_1)}
    {(1-p^2 w_2/w_1)(1+p^\dcox w_2/w_1)}.
  \end{displaymath}
\end{corollary}

\begin{pf}
  Note that $\omega_{\v{0},\v{0}}(w)= \sum_{t\in\N}w^{t}g_0(t)$.
  Multiplying both sides of~\eqref{eq:A2even-recurrence} by $w^{t}$
  ($w:=w_2/w_1$) and summing over non-negative $t$ we get
  \begin{eqnarray*}
    \left(1-(p^2-p^\dcox)w-p^{\dcox+2}w^2\right) \sum_{t\in\N}
    w^{t}g_0(t) = \qquad\\
    1-(1+p^{\dcox+1})w+p^{\dcox+1}w^2,
  \end{eqnarray*}
  from which the result follows.
\end{pf}
Hence
\begin{displaymath}
  \frac{\omega_{\v{0},\v{0}}(w_2/w_1)}{\phi_{\v{0},\v{0}}(w_2/w_1)}=
  \frac{(-p^{\dcox+2}w_2/w_1;p^{2\dcox})_\infty
    (w_2/w_1;p^{2\dcox})_\infty}
  {(p^2 w_2/w_1;p^{2\dcox})_\infty (-p^\dcox w_2/w_1;p^{2\dcox})_\infty},
\end{displaymath}
in agreement with~\eqref{eq:boson-2pt-level1-g}.

\subsection{Type $B^{(1)}_n$}
Recall that we have two ground state sequences ($\kappa=0,1$)
\begin{align*}
  \bo_m&=
  \begin{cases}
    \b{\kappa} &\text{for }m\in2\Z,\\
    z^\kappa\b{-\kappa} &\text{for } m\in2\Z+1,
  \end{cases}\\
  \intertext{and}
  \lambda_m&=
  \begin{cases}
    (1-\kappa)\Lambda_n+\kappa(\Lambda_1-\frac{m}{2}\delta)
    &\text{for } m\in2\Z,\\
    (1-\kappa)\Lambda_{n}+\kappa(\Lambda_0-\frac{m-1}{2}\delta)
    &\text{for } m\in2\Z+1.
  \end{cases}
\end{align*}

{}From~\cite{DO}, (up to a factor of a constant power in $w_2/w_1$) we
have in our notation~\eqref{eq:vertex-2pt}
\begin{eqnarray*}
  \phi_{\vo_{m-1},\vo_m}(w_2/w_1) &=&
  (1+p^{\dcox+1}w_2/w_1)^{\delta_{\kappa,0}}\\
  &&\quad \times\frac{(p^{2\dcox+2}w_2/w_1;p^{2\dcox})_\infty
    (p^{3\dcox}w_2/w_1;p^{2\dcox})_\infty}
  {(p^{\dcox+2} w_2/w_1;p^{2\dcox})_\infty
    (p^{2\dcox}w_2/w_1; p^{2\dcox})_\infty}.
\end{eqnarray*}

By a diagram automorphism, it is sufficient just to consider the case
when $m$ is even.  Let $m\in2\Z$.  Define
\begin{displaymath}
  g_j(t):= \dvac{m-1} z^{t+\kappa}\v{j} \wedge
z^{-t}\v{-j}\wedge\vac{m+1}\qquad (j\in J; t\in\Z).
\end{displaymath}
Note that $g_{-\kappa}(0)=1$, $g_{-j}(-\kappa)=0$ ($j\in\range{1,n}$)
and $g_j(t)=0$ ($j\in J$; $t\in\Z_{<0}$).

\begin{proposition}
  $g_{-\kappa}(t)$ satisfies the following recurrence relation
  \begin{gather}\label{eq:B1-recurrence}
    g_{-\kappa}(t)
    -(p^2-p^\dcox)g_{-\kappa}(t-1)-p^{\dcox+2}g_{-\kappa}(t-2)=\qquad
    \notag\\
    \qquad \delta_{t,0} -(1-\delta_{\kappa,0}p^{\dcox+1})\delta_{t,1}
    -\delta_{\kappa,0}p^{\dcox+1}\delta_{t,2}.
  \end{gather}
\end{proposition}

\begin{pf}
  The proof is like the proof of Proposition~\ref{prop:A2-recurrence}
  for type $A^{(2)}_{2n}$, but using
  \begin{align*}
    \A_t := &\bigl(Z(t,1+\kappa)+
    p^{1+\dcox\delta_{\kappa,0}}Z(t+\kappa-1,1+\kappa)\bigr)C_{0,0}\\
    &+{[2]} Z(t-1,\kappa)\bigl((-p)^nC_{1,-1}
    +(1-p^2)\sum_{j=2}^n(-p)^{n+1-j}C_{j,-j} \bigr)\\
    & +{[2]}p^{-\dcox\delta_{\kappa,1}}
    \bigl((1-p^2)Z(t+\kappa-1,1+\kappa)
    \sum_{j=2}^n(-p)^{n+j-1}C_{-j,j}\\
    &\qquad\qquad\qquad
    +(-p)^nZ(t+\kappa,2+\kappa)C_{-1,1}\bigr).
  \end{align*}
\end{pf}

The two point function of Fock intertwiners~\eqref{eq:fock-2pt} is
given by
\begin{corollary} Let $m\in2\Z$.
  \begin{displaymath}
    \omega_{\vo_{m-1},\vo_{m}}(w_2/w_1)=
    \frac{(1-w_2/w_1)(1+p^{\dcox+1} w_2/w_1)^{\delta_{\kappa,0}}}
    {(1-p^2 w_2/w_1)(1-p^\dcox w_2/w_1)}.
  \end{displaymath}
\end{corollary}

\begin{pf}
  Note that $\omega_{\vo_{m-1},\vo_{m}}(w)=
  \sum_{t\in\N}w^{t}g_{-\kappa}(t)$.  Multiplying both sides
  of~\eqref{eq:B1-recurrence} by $w^{t}$ ($w:=w_2/w_1$) and summing
  over non-negative $t$ we get
  \begin{eqnarray*}
    \left(1-(p^2-p^{\dcox})w-p^{\dcox+2}w^2\right) \sum_{t\in\N}
    w^k g_{-\kappa}(t)= \qquad\\
    1-(1-\delta_{\kappa,0}p^{\dcox+1})w-\delta_{\kappa,0}p^{\dcox+1}w^2,
  \end{eqnarray*}
  from which the result follows.
\end{pf}
Hence
\begin{displaymath}
  \frac{\omega_{\vo_{m-1},\vo_{m}}(w_2/w_1)}
  {\phi_{\vo_{m-1},\vo_{m}}(w_2/w_1)}=
  \frac{(p^{\dcox+2}w_2/w_1;p^{2\dcox})_\infty
    (w_2/w_1;p^{2\dcox})_\infty} {(p^2 w_2/w_1;p^{2\dcox})_\infty
    (p^\dcox w_2/w_1;p^{2\dcox})_\infty},
\end{displaymath}
in agreement with~\eqref{eq:boson-2pt-level1-g}.

\subsection{Type $A^{(2)}_{2n-1}$}
Recall that we have
\begin{displaymath}
  \bo_m=
  \begin{cases}
    \b{1} &\text{for } m\in 2\Z,\\
    z\b{-1} &\text{for } m\in2\Z+1,
  \end{cases}
  \quad \text{and}\quad
  \lambda_m=
  \begin{cases}
    \Lambda_1-\frac{m}{2}\delta &\text{for } m\in 2\Z,\\
    \Lambda_0-\frac{m-1}{2}\delta &\text{for } m\in2\Z+1.
  \end{cases}
\end{displaymath}

{}From~\cite{DO}, (up to a factor of a constant power in $w_2/w_1$) we
have in our notation~\eqref{eq:vertex-2pt}
\begin{displaymath}
  \phi_{\vo_{m-1},\vo_m}(w_2/w_1) =
  \frac{(q^{2\dcox+2}w_2/w_1;q^{2\dcox})_\infty
    (-q^{3\dcox}w_2/w_1;q^{2\dcox})_\infty}
  {(-q^{\dcox+2} w_2/w_1;q^{2\dcox})_\infty
    (q^{2\dcox}w_2/w_1; q^{2\dcox})_\infty}.
\end{displaymath}

By a diagram automorphism it is sufficient just to consider the case
when $m$ is even.  Let $m\in2\Z$.  Define
\begin{displaymath}
  g_j(t):= \dvac{m-1} z^{t+1}\v{j} \wedge
  z^{-t}\v{-j}\wedge\vac{m+1}\qquad (j\in J; t\in\Z).
\end{displaymath}
Note that $g_{-1}(0)=1$ and $g_j(t)=0$ ($j\in J$; $t\in\Z_{<0}$).

\begin{proposition}
  $g_{-1}(t)$ satisfies the following recurrence relation
  \begin{equation}\label{eq:A2odd-recurrence}
    g_{-1}(t) -(q^2-q^\dcox)g_{-1}(t-1)-q^{\dcox+2}g_{-1}(t-2)=
    \delta_{t,0} -\delta_{t,1}.
  \end{equation}
\end{proposition}

\begin{pf}
  The proof is like the proof of Proposition~\ref{prop:A2-recurrence},
  but using
  \begin{align*}
    \A_t := &\; Z(t,1)(z\tensor z)C_{-1,1}\\
    &+Z(t,2)(1-q^2)\sum_{j=2}^n (-q)^{j-1} C_{-j,j}\\
    & -Z(t-1,1) q^\dcox\bigl(C_{1,-1}
    +(1-q^2)\sum_{j=2}^n(-q)^{1-j}C_{j,-j}\bigr).
  \end{align*}
\end{pf}

The two point function of Fock intertwiners~\eqref{eq:fock-2pt} is
given by
\begin{corollary} Let $m\in2\Z$.
  \begin{displaymath}
    \omega_{\vo_{m-1},\vo_{m}}(w_2/w_1)=
    \frac{(1-w_2/w_1)}{(1-q^2 w_2/w_1)(1+q^\dcox w_2/w_1)}.
  \end{displaymath}
\end{corollary}

\begin{pf}
  Note that $\omega_{\vo_{m-1},\vo_{m}}(w)=
  \sum_{t\in\N}w^{t}g_{-1}(t)$.  Multiplying both sides
  of~\eqref{eq:A2odd-recurrence} by $w^{t}$ ($w:=w_2/w_1$) and summing
  over non-negative $t$ we get
  \begin{displaymath}
    \left(1-(q^2-q^{\dcox})w-q^{\dcox+2}w^2\right) \sum_{t\in\N}
    w^k g_{-1}(t)= 1-w,
  \end{displaymath}
  from which the result follows.
\end{pf}
Hence
\begin{displaymath}
  \frac{\omega_{\vo_{m-1},\vo_{m}}(w_2/w_1)}
  {\phi_{\vo_{m-1},\vo_{m}}(w_2/w_1)}=
  \frac{(-q^{\dcox+2}w_2/w_1;q^{2\dcox})_\infty
    (w_2/w_1;q^{2\dcox})_\infty} {(q^2 w_2/w_1;q^{2\dcox})_\infty
    (-q^\dcox w_2/w_1;q^{2\dcox})_\infty},
\end{displaymath}
in agreement with~\eqref{eq:boson-2pt-level1-g}.

\subsection{Type $D^{(1)}_n$}
Recall that we have two ground state sequences ($\kappa=0,1$)
\begin{displaymath}
  \bo_m=
  \begin{cases}
    \b{n\delta_{\kappa,0}+1\delta_{\kappa,1}} &\text{for } m\in 2\Z,\\
    z^\kappa\b{n\delta_{\kappa,0}-1\delta_{\kappa,1}} &\text{for } m\in2\Z+1,
  \end{cases}
\end{displaymath}

{}From~\cite{DO}, (up to a factor of a constant power in $w_2/w_1$) we
have in our notation~\eqref{eq:vertex-2pt}
\begin{displaymath}
  \phi_{\vo_{m-1},\vo_m}(w_2/w_1) =
  \frac{(q^{2\dcox+2}w_2/w_1;q^{2\dcox})_\infty
    (q^{3\dcox}w_2/w_1;q^{2\dcox})_\infty}
  {(q^{\dcox+2} w_2/w_1;q^{2\dcox})_\infty
    (q^{2\dcox}w_2/w_1; q^{2\dcox})_\infty}.
\end{displaymath}

By diagram automorphisms, it is sufficient just to consider the case
when $\kappa=0$ and $m$ is odd.
Let $\kappa=0$ and $m\in2\Z+1$.
Define
\begin{displaymath}
  g_j(t):= \dvac{m-1}z^{t}\v{j}\wedge z^{-t}\v{-j}\wedge\vac{m+1}
  \qquad (j\in J;t\in\Z).
\end{displaymath}
Note that $g_{n}(0)=1$, $g_{-j}(0)=0$ ($j\in\range{1,n}$) and
$g_j(t)=0$ ($j\in J$; $t\in\Z_{<0}$).

\begin{proposition}
  $g_{n}(t)$ satisfies the following recurrence relation
  \begin{equation}\label{eq:D1-recurrence}
    g_{n}(t) -(q^2-q^\dcox)g_{n}(t-1)-q^{\dcox+2}g_{n}(t-2)=
    \delta_{t,0} -\delta_{t,1}.
  \end{equation}
\end{proposition}

\begin{pf}
  The proof is like the proof of Proposition~\ref{prop:A2-recurrence},
  but using
  \begin{align*}
    \A_t := &\;Z(t,1) C_{n,-n}
    -q^\dcox Z(t-1,1)C_{-n,n}\\
    & +Z(t-1,0)\bigl((-q)^{n-1}C_{1,-1}
    +(1-q^2)\sum_{j=2}^{n-1}(-q)^{n-j} C_{j,-j}\bigr)\\
    & -(-q)^{n-1}\bigl(Z(t-1,0)(z\tensor z)C_{-1,1}\\
    & \qquad\qquad+(1-q^2)Z(t-1,1)\sum_{j=2}^{n-1}(-q)^{j-1}C_{-j,j}
    \bigr).
  \end{align*}
\end{pf}

The two point function of Fock intertwiners~\eqref{eq:fock-2pt} is
given by
\begin{corollary} Let $m\in2\Z+1$.
  \begin{displaymath}
    \omega_{\vo_{m-1},\vo_{m}}(w_2/w_1)=
    \frac{(1-w_2/w_1)}{(1-q^2 w_2/w_1)(1-q^\dcox w_2/w_1)}.
  \end{displaymath}
\end{corollary}

\begin{pf}
  Note that $\omega_{\vo_{m-1},\vo_{m}}(w)=
  \sum_{t\in\N}w^{t}g_{n}(t)$.  Multiplying both sides
  of~\eqref{eq:D1-recurrence} by $w^{t}$ ($w:=w_2/w_1$) and summing
  over non-negative $t$ we get
  \begin{displaymath}
    \left(1-(q^2-q^{\dcox})w-q^{\dcox+2}w^2\right) \sum_{t\in\N}
    w^k g_{n}(t)= 1-w,
  \end{displaymath}
  from which the result follows.
\end{pf}
Hence
\begin{displaymath}
  \frac{\omega_{\vo_{m-1},\vo_{m}}(w_2/w_1)}
  {\phi_{\vo_{m-1},\vo_{m}}(w_2/w_1)}=
  \frac{(q^{\dcox+2}w_2/w_1;q^{2\dcox})_\infty
    (w_2/w_1;q^{2\dcox})_\infty} {(q^2 w_2/w_1;q^{2\dcox})_\infty
    (q^\dcox w_2/w_1;q^{2\dcox})_\infty},
\end{displaymath}
in agreement with~\eqref{eq:boson-2pt-level1-g}.

\subsection{Type $D^{(2)}_{n+1}$} \label{subsec:d2-2pt}
This type is somewhat special because of the fact that $\Vaff$ is not
irreducible.  The dual Coxeter number is $\dcox=2n$.  We define
$p=q^2$ and $\xi^2=p^\dcox$.

Recall that we have two ground state sequences ($\kappa=0,\phi$)
\begin{displaymath}
  \bo_m = z^{-m} \b{\kappa}.
\end{displaymath}

By a diagram automorphism, it is sufficient just to consider one of
the two cases $\kappa\in\set{0,\phi}$.  We choose $\kappa=0$.

The boson commutator $\gamma_m = [B_m,B_{-m}]$ is given as follows:
\begin{proposition}\label{prop:d2gam}
\begin{equation}
  \gamma_m =
  \begin{cases}
    m \frac{(1+\xi^m)}{(1-2p^m-\xi^m)} &\text{for } m\in2\Z,\\
    m &\text{for } m\in2\Z+1.
  \end{cases}
\end{equation}
\end{proposition}

This corresponds to the following boson two point function
\begin{equation} \label{eq:D2-boson-2pt}
  \theta(w)= (1-w)\frac{(\xi^4w^2;\xi^4)_\infty
    (p^2\xi^2 w^2;\xi^4)_{\infty}}
  {(\xi^2 w^2;\xi^4)_\infty(p^2 w^2;\xi^4)_\infty}.
\end{equation}

{}From Appendix~\ref{app:c} we have
\begin{lemma} \label{lem:d2ph}
Let $w=w_2/w_1$.
\begin{displaymath}
  \phi_{z^{1-m}\v{0},z^{-m}\v{0}}(w_2/w_1)=
  (1+p \xi^2 w^2) \frac{(\xi^6 w^2; \xi^4)_{\infty}
    (p^2\xi^4 w^2; \xi^4)_{\infty}}{(\xi^4 w^2; \xi^4)_{\infty}
    (p^2\xi^2 w^2; \xi^4)_{\infty}}.
\end{displaymath}
\end{lemma}
It is sufficient just to consider just the case $m=0$.
Define
\begin{displaymath}
  g_j(t):= \dvac{-1} z^{t+1}\v{j} \wedge
  z^{-t}\v{-j}\wedge\vac{1}\qquad (j\in J; t\in\Z).
\end{displaymath}
Note that $g_0(0)=1$, $g_{-j}(0)=0$ ($j\in\range{1,n}$) and $g_j(t)=0$
($j\in J$; $t\in\Z_{<0}$).

\begin{proposition}
  $g_0(t)$ satisfies the following recurrence relation
  \begin{eqnarray} \label{eq:D2-recurrence}
    g_0(t) -(p^2+\xi^2)g_0(t-2) +p^2\xi^2g_0(t-4)=\qquad \notag\\
    \qquad \delta_{t,0} -\delta_{t,1} +p\xi^2\delta_{t,2}
    -p\xi^2\delta_{t,3}.
  \end{eqnarray}
\end{proposition}

\begin{pf}
  The proof is like the proof of Proposition~\ref{prop:A2-recurrence},
  but using
  \begin{align*}
    \A_t := & \bigl(Z(t,1)(z\tensor z) +p^{\dcox+1}Z(t-1,3)\bigr)
    C_{0,0}
    \\
    & +Z(t-1,1){[2]}\bigl(-q(-p)^n(z\tensor z)C_{\phi,\phi}\\
    &\qquad\qquad\qquad\quad +(1-p^2)\sum_{j=1}^n (-p)^{n+1-j}
    C_{j,-j}\bigr)\\
    & -Z(t-1,3){[2]}(1-p^2)\sum_{j=1}^n
    (-p)^{n+j}C_{-j,j}.
  \end{align*}
\end{pf}

The two point function of Fock intertwiners~\eqref{eq:fock-2pt} is
given by
\begin{corollary}
  \begin{displaymath}
    \omega_{z\v{0},\v{0}}(w_2/w_1)=
    \frac{(1-w_2/w_1)(1+p\xi^2(w_2/w_1)^2)}
    {(1-p^2 (w_2/w_1)^2)(1-\xi^2(w_2/w_1)^2)}.
  \end{displaymath}
\end{corollary}

\begin{pf}
  Note that $\omega_{z\v{0},\v{0}}(w)= \sum_{t\in\N}w^{t}g_0(t)$.
  Multiplying both sides of~\eqref{eq:D2-recurrence} by $w^{t}$
  ($w:=w_2/w_1$) and summing over non-negative $t$ we get
  \begin{displaymath}
    \left(1-(p^2+\xi^2)w^2-p^2\xi^2w^4\right) \sum_{t\in\N}
    w^k g_0(t)= 1-w +p\xi^2w^2 -p\xi^2w^3,
  \end{displaymath}
  from which the result follows.
\end{pf}
Finally we have $\frac{\omega_{z\v{0},\v{0}}(w)}
{\phi_{z\v{0},\v{0}}(w)}=\eqref{eq:D2-boson-2pt}$, which
proves Proposition~\ref{prop:d2gam}.

\Remark It is possible to work in an irreducible component of
$V_{\aff}$, say $V_{\aff}^{\text{even}}=V_0\tensor\C[z^2,z^{-2}]+
v_{\phi}\tensor z\C[z^2,z^{-2}]$.  On $V_{\aff}^{\text{even}} \otimes
V_{\aff}^{\text{even}}$ the energy function takes only even values.
The condition $H(\bo_m \otimes \bo_{m+1})=1$ for a ground state
sequence $\{\bo_m\}_{m \in Z}$ should then be replaced by $H(\bo_m
\otimes \bo_{m+1})=2$ for all $m\in Z$.

The ground state sequence in $\Baff^{\text{even}}$ is given by
$\bo_m=\b{0}$ for all $m\in Z$. The Fock two-point function can be
shown to be given by
\begin{equation}
\omega_{\v{0},\v{0}}(w)=
\frac{(1-w^2)(1+p \xi^2 w^2)}{(1-\xi^2 w^2)(1-p^2 w^2)},
\end{equation}
where $w=w_1/w_2$. Comparing with Lemma~\ref{lem:d2ph}, we find that
$\gamma_m$ is now given by the same formula as in
Theorem~\ref{thm:boson-comm}
\begin{equation}
\gamma_m = m\frac{1+\xi^{2m}}{1-p^{2m}}.
\end{equation}



\newsection{Higher level examples: level $k$ $A^{(1)}_1$}

\subsection{Cartan datum}
$I=\set{0,1}$. The Dynkin diagram for $A^{(1)}_1$ is
\begin{displaymath}
  \dynkin{0}=\!\!=\!\!= \dynkin{1}\quad.
\end{displaymath}
We have
\begin{eqnarray*}
  \delta &=& \alpha_0 +\alpha_1,\\
  c &=& h_0 +h_1,\\
  (\alpha_i,\alpha_i) &=& 2 \qquad (i\in I).
\end{eqnarray*}

\subsection{Perfect crystal}
Fix $k\in\Zplus$.
Let $J:=\range{0,k}$.  Let $V$ be the
$(k+1)$-dimensional $U'_q(A^{(1)}_1)$-module with the level~$k$ perfect
crystal $B:=\set{\b{j}}_{j\in J}$ with crystal graph:
\begin{displaymath}
  \b{0} \rightleftcrystal{1}{0} \b{1} \rightleftcrystal{1}{0} \b{2}
  \rightleftcrystal{1}{0} \cdots \rightleftcrystal{1}{0} \b{k}\quad .
\end{displaymath}
The elements of $B$ have the following weights
\begin{eqnarray*}
  \wt(\b{j}) &=& (k/2 -j)\alpha_1\\
  &=& (k-2j)(\Lambda_1-\Lambda_0)
  \qquad (j\in J).
\end{eqnarray*}
Let $\v{j}:=G(\b{j})$ ($j\in J$).  The action of
$U'_q(A^{(1)}_1)$ on $\v{j}\in V$ obeys~\eqref{eq:LG}.

\subsection{Energy function}
The energy function $H$ has the following values on $B\tensor B$
\begin{displaymath}
H(\b{i}\tensor\b{j})=\min(i,k-j) \qquad (i,j\in J).
\end{displaymath}
Write $H(i,j)$ for $H(\b{i}\tensor\b{j})$ ($i,j\in J$).

The Coxeter number of $A^{(1)}_1$ is $h=2$.  We take
\begin{displaymath}
  l(z^m\b{j}) = 2m-j\qquad (m\in\Z,j\in J).
\end{displaymath}
$H$ and $l$ satisfy condition~(L).  Note that
$l(z^m\b{j})=l(z^{m+1}\b{j+2})$ ($m\in\Z$ and $j\in\range{0,k-2}$), so
$l$ gives a partial ordering of $\Baff$ for $k>1$.

\subsection{$q$-binomials}
Define the $q$-binomial coefficient $\qbin{m}{n}$ ($m,n\in\Z$) by
\begin{displaymath}
  \qbin{m}{n}=
  \begin{cases}
    \frac{[m][m-1] \cdots [m-n+1]}
    {[n][n-1] \cdots [1]} &: m\geq n\geq 0\\
    0 &: \text{otherwise}.
  \end{cases}
\end{displaymath}
We will often write sums involving $q$-binomial coefficients as sums
over all integers.  The advantage is that we can then freely change
variables without worrying about the range of summation.  The
following result is widely used in the sequel:

\begin{lemma}\label{lem:qbin}
  \begin{enumerate}
  \item 
    For any $\eta\in\C(q)$ and $n\in\Zplus$, we have
    \begin{displaymath}
      \sum_{j \in \Z} (-\eta)^j  \qbin{n}{j} = \prod_{i=0}^{n-1}
      (q^{n-1-2i}-\eta).
    \end{displaymath}
  \item 
    The sum in~(i) vanishes if $\eta=q^m$ with $m$ an
    integer lying in the range $[-n+1,n-1]_2$. Here $[a,a+2b]_2$ means
    $\{a+2i; 0\leq i \leq b \}$.
  \end{enumerate}
\end{lemma}

\subsection{Wedge relations}
Define the vectors $C_{i,j} \in N$ ($i,j \in J$) by
\begin{displaymath}
  C_{i,j}  =
  \begin{cases}
   (z\otimes z)^{-i} e_0^{(i)} f_1^{(j)}\left(\v{0} \otimes \v{0}
   \right) &:  i+j\leq k,\\
   e_1^{(k-i)} f_0^{(k-j)} \left(\v{k} \otimes \v{k} \right) &: i+j> k.
 \end{cases}
\end{displaymath}
Explicitly, we have
\begin{equation} \label{eqn:cc}
  C_{i,j}=
  \begin{cases}
  \sum_{i',j',a,b}  q^{(k-j')(i'-b)+(k-i')a} \qbin{j'}{a}\qbin{i'}{b}
  z^{-a}\v{i'} \otimes z^{-b} \v{j'} &: i'+j' \leq k\\
   \sum_{i',j',a,b}  q^{i'(k-j'-b)+j'a} \qbin{k-i'}{a}\qbin{k-j'}{b}
   z^{-a}\v{i'} \otimes z^{-b} \v{j'} &: i'+j'> k.
 \end{cases}
\end{equation}
The summation in both cases is over $i',j' \in J$
and $a,b \geq 0$ with
\begin{eqnarray*}
i'+j' & = & i+j,\\
a+b & = & H(i,j).
\end{eqnarray*}

\begin{proposition}
  Identify $C_{i,j}$ with $C_{\b{i}, z^{-H(i,j)}\b{j}}$.  Then
  $\set{z^m\tensor z^m\cdot C_{i,j}}_{m\in\Z;i,j\in J}$ with the
  function $l$ satisfy condition~(R) of
  subsection~\ref{subsec:wedge-prod}.
\end{proposition}

\subsection{Fock space}
We have
\begin{eqnarray*}
  \Bmin&=&B,\\
  (\Pcl^+)_k&=&\set{j\Lambda_0^\cl+(k-j)\Lambda_1^\cl}_{j\in J},
\end{eqnarray*}
with the bijections
\begin{eqnarray*} 
  \varepsilon(\b{j}) &=&  (k-j)\Lambda_0^\cl +j\Lambda_1^\cl,\\
  \varphi(\b{j}) &=& j\Lambda_0^\cl +(k-j)\Lambda_1^\cl.
\end{eqnarray*}
Fix $\kappa \in J$ and let $\kappa'=k-\kappa$.  We have
$H(z^{-\ell(k-2)}b_{\kappa} \otimes z^{-\ell(k-2)-\kappa+1} b_{\kappa'})=
H(b_{\kappa}\otimes b_{\kappa'})-\kappa+1=1$.
$H(z^{-\ell(k-2)-\kappa+1} b_{\kappa'} \otimes
z^{-(\ell+1)(k-2)}b_{\kappa})=1$ and so the following is a ground
state sequence ($\ell\in\Z$)
\begin{equation}
  \begin{aligned}
    \bo_{2\ell-1} &= z^{-\ell(k-2)}b_{\kappa},\\
    \bo_{2\ell} &= z^{-\ell(k-2)-\kappa+1} b_{\kappa'},
  \end{aligned}
\label{eqn:ground}
\end{equation}
with
\begin{displaymath}
  \cl(\lambda_{m}) =
  \begin{cases}
    \kappa \Lambda_0 + \kappa' \Lambda_1 &: \text{if $m$ is odd},\\
    \kappa' \Lambda_0 + \kappa \Lambda_1 &: \text{if $m$ is even}.
  \end{cases}
\end{displaymath}
With $\vo_m= G(\bo_m)$, the vacuum vector of $\fock_m$ is
\begin{displaymath}
  \vac{m}= \vo_m\wedge \vo_{m+1}\wedge \vo_{m+2}\wedge\cdots\cdots
\end{displaymath}
with weight $\lambda_m$.

\subsection{Two point functions}
A priori $\gamma_n=[B_n,B_{-n}]$ may depend on the choice of $\kappa$.
However, we find that it is independent of $\kappa$.
\begin{theorem} 
\begin{displaymath}
\gamma_n
 = n \frac{1-q^{4n}}{1-q^{2n}-q^{4n}+q^{2(k+1)n}}.
\end{displaymath}
\end{theorem}
The theorem follows by applying Proposition~\ref{prop:2pt-decomp}
to Proposition~\ref{prop:ph} and Corollary~\ref{corol:om} below.

{}From~\cite{IIJMNT} we have

\begin{proposition} \label{prop:ph}
\begin{displaymath}
  \phi_{\vo_{2\ell-1},\vo_{2\ell}}(w)=
  \frac{(q^{2(k+2)} w; q^4)_{\infty}}
    {(q^4 w; q^4)_{\infty}} \sum_{p=0}^\infty \left(q^{k+2}w\right)^p
    \qbin{\kappa}{p}\qbin{\kappa'}{p},
\end{displaymath}
where $w=w_2/w_1$.
\end{proposition}

Without loss of generality, we can choose $\ell=0$.  Define
\begin{displaymath}
  g_\kappa(t):= \dvac{-1}z^t\v{\kappa}\wedge
  z^{-t+1-\kappa}\v{\kappa'}\wedge\vac{+1}
  \qquad (t\in\Z;j\in J).
\end{displaymath}
Note that $g_\kappa(t)=\delta_{t,0}$ for
$t\in\Z_{\leq0}$ by Theorem~\ref{kern}.

\begin{proposition}\label{prop:key}
  $g_\kappa(t)$ satisfies the following recurrence relation
  \begin{equation}
    \sum_{\alpha\in\Z} \left(-q^{k+1}\right)^{\alpha} \qbin{k}{\alpha}
    g_\kappa(t-\alpha)
    - \left(q^{k+2}\right)^t \qbin{\kappa}{t} \qbin{\kappa'}{t}
    + \left(q^{k+2}\right)^{t-1} \qbin{\kappa}{t-1} \qbin{\kappa'}{t-1}=0
    \label{eqn:rec}
  \end{equation}
\end{proposition}

\begin{corollary} \label{corol:om}
\begin{displaymath}
  \omega_{\vo_{-1},\vo_{0}}(w) =
  \frac{(1-w)}{\prod_{j=1}^k(1- q^{2j} w)}
  \sum_{p=0}^\infty \left(q^{k+2}w\right)^p
  \qbin{\kappa}{p}\qbin{\kappa'}{p},
\end{displaymath}
where $w=w_2/w_1$.
\end{corollary}

\begin{pf}
  We have
\begin{equation}
 \omega_{\vo_{-1},\vo_{0}}(w)=
 \sum_{j \in\Z} \left(\frac{w_2}{w_1}\right)^j g_\kappa(j).
\label{eqn:omom}
\end{equation}
Multiply both sides of (\ref{eqn:rec}) by $w^t$
and sum over all $t \geq 0$.  After relabelling of $t$ and using
(\ref{eqn:omom}) we obtain
\begin{equation}
\omega_{\vo_{-1},\vo_{0}}(w)
\sum_{\alpha\in\Z} \left(-q^{k+1} w\right)^{\alpha} \qbin{k}{\alpha}
  =    (1-w) \sum_{t=0}^\infty
    \left(q^{k+2}w\right)^t \qbin{\kappa}{t} \qbin{\kappa'}{t}.
\label{eqn:gf}
\end{equation}
{}From Lemma \ref{lem:qbin}~(i) we have
\begin{displaymath}
\sum_{\alpha\in\Z} \left(-q^{k+1} w\right)^{\alpha} \qbin{k}{\alpha}  =
\prod_{j=1}^k (1-q^{2j} w),
\end{displaymath}
thereby proving the result.
\end{pf}

Only Proposition~\ref{prop:key} remains to be proved.

\subsection{Proof of recurrence relation}
Let $Z(t,d)$ be the operator defined in~\S\ref{sec6:summary}:
\begin{displaymath}
  Z(t,d)= z^t\tensor z^{d-t} +\delta(2t>d)z^{d-t}\tensor z^t
  -\delta(2t<d)z^t\tensor z^{d-t}\qquad (t,d\in\Z).
\end{displaymath}
For $t \in\Z$, define
\begin{displaymath}
  \A_t :=
  \sum
  \begin{Sb}
    i\in J\\
    \gamma\in\Z
  \end{Sb}
  (-q^{\kappa+1})^{k-i-\kappa} q^{\gamma(k+2)}
  \qbin{i}{\gamma}\qbin{k-i}{\kappa-\gamma} Z(t-\gamma,-i+\kappa'+1)
  C_{k-i,i}.
  \label{eqn:apdef}
\end{displaymath}
We split the proof into three parts.  Define
\begin{eqnarray}
Z^{(1)}(t,d) &=& z^t\tensor z^{d-t},\nonumber\\
Z^{(2)}(t,d) &=& -z^t\tensor z^{d-t}\delta(2t<d),\\
Z^{(3)}(t,d) &=& z^{d-t}\tensor z^t\delta(2t>d).\nonumber
\label{eqn:a123}
\end{eqnarray}
Then $Z(t,d)=Z^{(1)}(t,d)+Z^{(2)}(t,d)+Z^{(3)}(t,d)$.  Let Define
$\A_t^{(i)}$ ($i\in\set{1,2,3}$) by replacing $Z$ by $Z^{(i)}$
in the definition of $\A_t$.  Then $\A_t=\A_t^{(1)}+{\cal
  A}_t^{(2)}+\A_t^{(3)}$.  We will deal with each ${\cal
  A}_t^{(i)}$ separately.  Note that $\A_t^{(i)}\not\in N$
($i\in\set{1,2,3}$), only $\A_t\in N$.

\subsubsection{$\A_t^{(1)}$}
{}From (\ref{eqn:cc}) we obtain
\begin{displaymath}
  C_{k-i,i} = \sum
  \begin{Sb}
    j\in J\\
    b\in\Z
  \end{Sb}
  q^{j^2
    -j(k+i) +k(k-b)} \qbin{j}{k-i-b} \qbin{k-j}{b}
  z^{b+i-k} \v{k-j} \otimes z^{-b} \v{j}.
\end{displaymath}
Substituting into (\ref{eqn:apdef})
and performing a change of variable $b \rightarrow \gamma+k-i-\alpha$,
followed by $i \rightarrow i+ \gamma$ we obtain
\begin{multline}
\A_t^{(1)} = \sum
\begin{Sb}
  i,\gamma,\alpha\in\Z\\
  j\in J
\end{Sb}
\left(-q^{\kappa+1}\right)^{k-i-\kappa-\gamma} q^{2\gamma+k\alpha+
  (i-j+\gamma)(k-j)}
\qbin{i+\gamma}{\gamma}\qbin{k-i-\gamma}{\kappa-\gamma} \\
\times\qbin{j} {\alpha-\gamma} \qbin{k-j}{k-i-\alpha} z^{t-\alpha}
\v{k-j}\tensor z^{-t+\alpha -\kappa+1} \v{j}.
\label{eqn:conj1}
\end{multline}
Let us now argue that {\em only} the $j=k-\kappa$ terms contribute in the above
sum. Recall that our convention for $q$-binomial coefficients implicitly
defines for us the upper and lower limits of summation in formulae like
$\A_t^{(1)}$. For instance, the constraints on $i$ are
\begin{equation}
\max \left(0,j-\alpha\right) \leq i \leq \min\left(k-\kappa,k-\alpha\right).
\label{eqn:lims}
\end{equation}

Let us assume first that $j\leq k-\kappa-1$. The strategy is to recast the sum
over $i$ in (\ref{eqn:conj1}), more specifically,
\begin{equation}
I_j=\qbin{j}{\alpha-\gamma}\sum_{i\in\Z} \left(-q^{k-j-\kappa-1}\right)^i
 \qbin{i+\gamma}{\gamma}\qbin{k-i-\gamma}{\kappa-\gamma} \qbin{k-j}{k-i-\alpha}
\label{eqn:idef}
\end{equation}
into a form such that Lemma \ref{lem:qbin}
 applies. Consider the case $j \leq \alpha \leq
\kappa$, so that according to (\ref{eqn:lims}) we have $ 0\leq i\leq k-\kappa$.
By manipulating the $q$-binomial coefficients we obtain
\begin{eqnarray}
I_j = \frac{[k-j]![j]!}{[k-\kappa]![\kappa-\gamma]![\gamma]!}
  \sum_{i=0}^{k-\kappa}
  \left(-q^{k-j-\kappa-1}\right)^i \qbin{k-\kappa}{i}
  \qbin{i+\gamma}{j-\alpha+\gamma} \qbin{k-i-\gamma}{\alpha-\gamma}.
\notag\\ \label{eqn:fc}
\end{eqnarray}
Now treat the product of the last two $q$-binomial coefficients in $I_j$
together with $\left(-q^{k-j-\kappa-1}\right)^i$ as a polynomial in $q^i$;
the powers of $q^i$ which appear can be seen to
lie in the range $[k-\kappa-1-2j,k-\kappa-1]_2$. In fact, due to the assumption
on $j$ the range is $[1-k+\kappa,k-\kappa-1]_2$. Therefore $I_j$ is a finite
sum of sums for which Lemma \ref{lem:qbin}~(ii) applies and thus vanishes.

For the other three remaining cases (a) $j,\kappa \leq \alpha$,
(b) $j,\kappa \geq \alpha$ and (c) $j \geq \alpha \geq \kappa$ we
use, respectively, the identities for $I_j$:
\begin{eqnarray*}
 I_j & =& \frac{[k-j]![j]!}{[k-\alpha]![\alpha-\gamma]![\gamma]!}
  \sum_{i=0}^{k-\alpha}
    \left(-q^{k-j-\kappa-1}\right)^i \qbin{k-\alpha}{i}
  \qbin{i+\gamma}{j-\alpha+\gamma}\nonumber\\
  & & \quad\times\qbin{k-i-\gamma}{\kappa-\gamma}\\
I_j & = & \frac{[k-j]![j]!}{[j-\alpha+\gamma]![\kappa-\gamma]!
   [k-\kappa-j+\alpha]!}
   \sum_{i=j-\alpha}^{k-\kappa}
   \left(-q^{k-j-\kappa-1}\right)^i  \nonumber\\
  & &\quad\times\qbin{k-\kappa-j+\alpha}
  {i-j+\alpha} \qbin{i+\gamma}{\gamma} \qbin{k-i-\gamma}{\alpha-\gamma} \\
I_j  & =& \sum_{i=j-\alpha}^{k-\alpha}
    \left(-q^{k-j-\kappa-1}\right)^i \qbin{k-j}{i-j+\alpha}
   \qbin{i+\gamma}{\gamma} \qbin{k-i-\gamma}{\kappa-\gamma}.
\end{eqnarray*}
In each case $I_j$ vanishes by application of Lemma \ref{lem:qbin}~(ii).

We have proved that the sum over $j<k-\kappa$ in (\ref{eqn:conj1}) vanishes.
The sum over $j>k-\kappa$ vanishes for similar reasons. Keeping only the
$j=k-\kappa$ term we arrive at
\begin{equation}
\A_t^{(1)} = \sum_{\gamma,\alpha\in\Z}
  \left(-q\right)^{k-\kappa-\gamma} q^{2\gamma + k\alpha}\; I_{k-\kappa} \;
  z^{t-\alpha}\v{\kappa}\tensor z^{ -t+\alpha
    -\kappa+1}\v{\kappa'}
\label{eqn:op}
\end{equation}
where $I_{k-\kappa}$ is given by (\ref{eqn:idef}). Once again, we have
the constraint (\ref{eqn:lims}) and have to treat the four cases
separately. We consider in detail only the case $k-\kappa \leq \alpha
\leq \kappa$, using the form (\ref{eqn:fc}) for $I_{k-\kappa}$. The
other three cases are similar. We proceed as before but now find that
the powers of $q^i$ lie in the range $[-1-k+\kappa,-1+k-\kappa]_2$.
By Lemma \ref{lem:qbin}~(ii) only the term whose power of $q^i$ is
$-1-k+\kappa$ survives. In other words,
\begin{displaymath}
I_{k-\kappa} = \qbin{\kappa}{\gamma} \sum_{i=0}^{k-\kappa} (-q^{-1})^i
    \qbin{k-\kappa}{i} \frac{q^{-(i+\gamma)} q^{-(i+\gamma-1)} \cdots}
  {[k-\kappa-\alpha+\gamma]!}
\frac{q^{k-i-\gamma} q^{k-i-\gamma-1} \cdots}
  {[\alpha-\gamma]!}
  \frac{(-)^{k-\kappa-\alpha+\gamma}}{(q-q^{-1})^{k-\kappa}}.
\end{displaymath}
Applying Lemma \ref{lem:qbin}~(i) and simplifying we find
\begin{displaymath}
I_{k-\kappa}=\qbin{\kappa}{\gamma} \qbin{k-\kappa}{\alpha-\gamma}
     \left(-q^{\kappa+1}\right)^{\alpha}
      \left(-q^{k+1}\right)^{-\gamma}
      \left(-q\right)^{\kappa-k}.
\end{displaymath}
Substituting into (\ref{eqn:op}) we obtain
\begin{equation}
\A_t^{(1)} = \sum_{\gamma,\alpha\in\Z}
 \left(-q^{k+\kappa+1}\right)^{\alpha-\gamma}
 \left(-q^{\kappa+1}\right)^{\gamma} \qbin{\kappa}{\gamma}
 \qbin{k-\kappa}{\alpha-\gamma} z^{t-\alpha}\v{\kappa}\tensor z^{ -t+\alpha
    -\kappa+1}\v{\kappa'}.
\label{eqn:one}
\end{equation}
We now note the identity
\begin{multline*}
\left(\sum_{\beta\in\Z} \left(-q^{k+\kappa+1} x\right)^{\beta}
   \qbin{k-\kappa}{\beta}\right)
\left(\sum_{\gamma\in\Z} \left(-q^{\kappa+1} x\right)^{\gamma}
   \qbin{\kappa}{\gamma}\right)
= \\
\prod_{i=1}^{k}\left(1-q^{2i} x\right) = \sum_{\alpha\in\Z}
\left(-q^{k+1} x \right)^{\alpha} \qbin{k}{\alpha},
\end{multline*}
which follows from the ubiquitous Lemma \ref{lem:qbin}, to
perform the $\gamma$-sum in~\eqref{eqn:one} with the result
\begin{displaymath}
  \A_t^{(1)} =\sum_{\alpha\in\Z}
  \left(-q^{k+1}\right)^{\alpha} \qbin{k}{\alpha}
  z^{t -\alpha} \v{\kappa} \otimes
  z^{-t +\alpha +1 -\kappa} \v{\kappa'}.
\label{eqn:ores}
\end{displaymath}

\subsubsection{$\A_t^{(2)}$}
The only difference between $\A_t^{(1)}$ and $\A_t^{(2)}$
is that the latter has a negative sign and an additional constraint
\begin{equation}
i < \gamma -2t+k-\kappa+1
\label{eqn:con}
\end{equation}
on the sum (denoted by prime) due to the
definition of $Z^{(2)}$:
\begin{multline}
  \A_t^{(2)} = -\sideset{}{'}\sum
  \begin{Sb}
    i,\gamma,\alpha\in\Z\\
    j\in J
  \end{Sb}
  (-q^{\kappa+1})^{k-i-\kappa-\gamma} q^{2\gamma+(i-j+\gamma)(k-j)}
  \qbin{i+\gamma}{\gamma}\qbin{k-i-\gamma}{\kappa-\gamma} \\
  \times\qbin{j} {\alpha-\gamma} \qbin{k-j}{k-i-\alpha} z^{t-\alpha}
  \v{k-j} \tensor z^{-t+\alpha -\kappa+1} \v{j}.
\label{eqn:two}
\end{multline}
Furthermore we are now interested in dropping terms that annihilate
the vacuum.
Using Theorem 3.5 this means that we require
\begin{eqnarray}
\lefteqn{H( z^{-t+\alpha -m(k-2) -\kappa+1} \b{j} \otimes
   z^{-(m+1)(k-2)} \b{\kappa})}\nonumber\\
& = &
t-\alpha+\kappa-k+1+\min(j,k-\kappa) > 0.
\label{eqn:req}
\end{eqnarray}

Let us assume first that $j > k-\kappa$. From (\ref{eqn:req})
we need $\alpha - t <0$.
Now from the last $q$-binomial in (\ref{eqn:two}) and (\ref{eqn:con}) we have
\begin{equation}
j \leq i+\alpha < (\alpha-t)+(\gamma-t)+k-\kappa+1.
\label{eqn:tp}
\end{equation}
Thus we have $k-\kappa < j < \gamma-t + k -\kappa+1$ and so $\gamma-t >0$.
But this means $\gamma > t \geq \alpha$ which contradicts the requirement
$\gamma\leq \alpha$ coming from the third $q$-binomial in (\ref{eqn:two}).

Next assume that $j < k-\kappa$. From (\ref{eqn:req}) we now need
$j \geq k-\kappa+
\alpha-t$ and thus $\alpha-t <0$. But again we have (\ref{eqn:tp}), and so
\begin{displaymath}
j \leq (\alpha-t)+(\gamma-t)+k-\kappa < (\gamma-t)+k-\kappa.
\end{displaymath}
Thus we have $k-\kappa+\alpha-t \leq j < \gamma-t+k-\kappa$ and so
$\alpha<\gamma$ which again  contradicts $\gamma\leq \alpha$.

Hence we must have
\begin{equation}
j=k-\kappa.
\label{eqn:c1}
\end{equation}
According to (\ref{eqn:req}) we need $\alpha-t <0$. But
again (\ref{eqn:tp}) is required, which leads to
$0\leq (\alpha-t)+(\gamma-t)\leq \gamma-t$. Therefore $\alpha\leq t\leq \gamma$
which together with $\gamma\leq \alpha$ from the third $q$-binomial in
(\ref{eqn:two}) makes mandatory
\begin{equation}
\alpha=\gamma=t.
\label{eqn:c2}
\end{equation}
This means that (\ref{eqn:con}) can be rephrased as $i\leq -t+k-\kappa$.
But from the last $q$-binomial in (\ref{eqn:two}), together with
(\ref{eqn:c1}) and (\ref{eqn:c2}) we must have also
\begin{equation}
i= -t+k-\kappa.
\label{eqn:c3}
\end{equation}
Substituting (\ref{eqn:c1}), (\ref{eqn:c2}) and (\ref{eqn:c3}) into
(\ref{eqn:two}) we arrive at
\begin{displaymath}
  \A_t^{(2)} =
  -\left(q^{k+2}\right)^t \qbin{\kappa}{t}
  \qbin{\kappa'}{t}\v{\kappa}\otimes z^{-\kappa+1} \v{k-\kappa} +
  \cdots.
\end{displaymath}

\subsubsection{$\A_t^{(3)}$}
One argues in the same way that
\begin{displaymath}
  \A_t^{(3)} =
  \left(q^{k+2}\right)^{t-1} \qbin{\kappa}{t-1} \qbin{\kappa'}{t-1}
  \v{\kappa} \otimes z^{-\kappa+1} \v{k-\kappa}
  + \cdots.
\end{displaymath}
Let $\A_t^\wedge$ denote the image of $\A_t$ in $\Vaff^{\wedge2}$.
Adding the three parts together, the relation $\dvac{-1}\A_t^\wedge
\wedge\vac{1}=0$ gives us Proposition~\ref{prop:key}.


\appendix
\renewcommand{\thelemma}{\thesection.\arabic{lemma}}
\renewcommand{\theequation}{\thesection.\arabic{equation}}

\setcounter{equation}{0}
\setcounter{lemma}{0}
\newsection{Perfect crystal}\label{perf}

Let $V$ be an integrable finite-dimensional $U'_q(\Gg)$-module with a
perfect crystal base $(L,B)$ of level $l$.  We assume that it has a
lower global base (i.e.~satisfies (G)).  In~\cite{KMN1}, we proved
that the ``semi-infinite tensor product'' $B\otimes B\otimes\cdots$ is
isomorphic to the crystal base $B(\lam)$ of the highest irreducible
module, provided that the rank of $\Gg$ is greater than $2$.  In this
appendix, we prove the same statement for any rank.  In~\cite{KMN1},
the proof is combinatorial, and here it is by the use of a vertex
operator.  Let us take a ground state sequence
$(\cdots,b^\circ_m,b^\circ_{m+1},\cdots)$ in $B_\aff$.  Set
$\vo_k=G(b^\circ_k)$.  For an integral dominant weight $\lam$, we
denote by $V(\lam)$ the irreducible $U_q(\Gg)$-module with highest
weight $\lam$ and highest weight vector $u_\lam$, and by
$\bigl(L(\lam),B(\lam)\bigr)$ its crystal base.

\Proposition
$B\otimes B(\lam_m)\cong B(\lam_{m-1})$.
\enproposition

This proposition implies the following result.

\Proposition
\begin{multline*}
  B(\lambda_m)\cong
  \{(b_m,b_{m+1},\cdots);\;
  b_k\in B_\aff,H(b_k\otimes b_{k+1})=1
  \text{ for any $k\ge m$}\\ \text{and $b_k=b^\circ_k$ for $k\gge m$}\}.
\end{multline*}
\enproposition

\medskip
The following lemma is proved in~\cite{DJO}.

\Lemma $\Hom_{U_q(\Gg)}\bigl(V_\aff\otimes V(\lam_m),V(\lam_{m-1})\bigr)=K$.
\enlemma

Let
$\Phi:V_\aff\otimes V(\lam_m)\to V(\lam_{m-1})$
be a $U_q(\Gg)$-linear homomorphism.
We normalize it by
$$\Phi(\vo_{m-1}\otimes u_{\lam_m})=u_{\lam_{m-1}}.$$

Then the following lemma is also proved in~\cite{DJO} in the dual form.

\Lemma
$\Phi(L_\aff\otimes_AL(\lam_m))\subset L(\lam_{m-1})$.
\enlemma

Let $\bar{\Phi}:(L_\aff\otimes_AL(\lam_m))/q(L_\aff\otimes_AL(\lam_m))
\to L(\lam_{m-1})/qL(\lam_{m-1})$
be the induced homomorphism.

The following two lemmas are easily proved.
\Lemma\label{lem:a1}
Let $M_j$ be an integrable $U_q(\Gg)$-module, and
$(L_j,B_j)$ a crystal base of $M_j$ for $j=1,2$.
Let $\psi:M_1\to M_2$ be a $U_q(\Gg)$-linear homomorphism
sending $L_1$ to $L_2$.
Let $\bar{\psi}:L_1/qL_1\to L_2/qL_2$
be the induced homomorphism.
Set $\tB=\{b\in B_1|\bar{\psi}(b)\in B_2\}$.
Then $\tB$ has a crystal structure
such that
$\iota:\tB\to B_1$ and $\bar{\psi}:\tB\to B_2$
are strict morphism of crystals.
\enlemma

Here a strict morphism means
a morphism commuting with $\te_i$ and $\tf_i$.

\Lemma\label{isom}
Let $\lam$ be a dominant integral weight.
Let $B$ be a semi-regular crystal
(i.e.~$\varepsilon_i(b)=\max\{n\in\BZ_{\ge0}|\te_i^nb\not=0\}$
and $\varphi_i(b)=\max\{n\in\BZ_{\ge0}|\tf_i^nb\not=0\}$).
We assume further that $B$ is connected.
\begin{enumerate}
\item 
If $\psi:B(\lam)\to B$ is a strict morphism such that
$\psi(B(\lam))\subset B$, then $\psi$ is an isomorphism.
\item 
If $\psi:B\to B(\lam)$ is a strict morphism such that
$\psi(B)\subset B(\lam)$, then $\psi$ is an isomorphism.
\end{enumerate}
\enlemma

Let $B'$ be the connected component
of $B_\aff\otimes B(\lam_m)$ containing
$b_{m-1}^\circ\otimes u_{\lam_m}$.
Then $\Phi$ sends $B'$ to $B(\lam_{m-1})$.
Hence $B'$ is a subcrystal of $B_\aff\otimes B(\lam_m)$,
and Lemma \ref{lem:a1} implies
$B'\to B_\aff\otimes B(\lam_m)$ and $B'\to B(\lam_{m-1})$
are strict morphisms.
Moreover any $b\in B'$ is not mapped to $0$ by
the morphism $B'\to B(\lam_{m-1})$.
Hence by Lemma \ref{isom},
$B'$ is isomorphic to $B(\lam_{m-1})$.
Hence we obtain a strict morphism
$B(\lam_{m-1})\to B_\aff\otimes B(\lam_m)$.
Composing it with $B_\aff\to B$,
we obtain a strict morphism
$B(\lam_{m-1})\to B\otimes B(\lam_m)$.

The following lemma is proved in~\cite{KMN1}.
\Lemma
$B\otimes B(\lam_m)$ is connected.
\enlemma

Thus $B(\lam_{m-1})\to B\otimes B(\lam_m)$
is an isomorphism by Lemma \ref{isom}.


\setcounter{equation}{0}
\setcounter{lemma}{0}
\newsection{{Serre relations}\label{sec:serre}}
Let $\tU$ be
the algebra associated
to a symmetrizable Kac-Moody algebra
with the same generators and the defining relations
as the quantized universal enveloping algebra
except the Serre relations.
Let $U_q(\Gg)_i$ be its subalgebra generated by
$e_i$, $f_i$ and $t_i^{\pm1}$.
In this appendix, we prove the following proposition.

\Proposition\label{serre}
Let $M$ be a  $\tU$-module.
Assume that $M$ is an integrable
$U_q(\Gg)_i$-module for every $i$.
Then the action of
$\tU$ on $M$ satisfies
the Serre relations.
\enproposition
Hence $M$ has the structure of a $U_q(\Gg)$-module.

Let $M$ and $N$ be integrable $U_q(\Gg)_i$-modules.
We endow the structure of $U_q(\Gg)_i$-module on
$\Hom(M,N)$
such that $\Hom(M,N)\otimes M\to N$ is $U_q(\Gg)_i$-linear.
Namely for $x\in U_q(\Gg)_i$ with $\Delta(x)=\sum x_{(1)}\otimes x_{(2)}$,
$x$ acts on $\psi\in\Hom(M,N)$ by
$x_{(1)}\psi a(x_{(2)})$.

Recall that an element $\psi$ of $\Hom(M,N)$
is called locally $U_q(\Gg)_i$-finite,
if it is contained in an integrable $U_q(\Gg)_i$-submodule.

\Lemma Let $M$ and $N$ be integrable
$U_q(\Gg)_i$-modules.
Assume that a weight vector $\psi\in \Hom(M,N)$
satisfies
$$f_i^{n+1}\psi=0\qbox{for some $n\ge0$.}$$
Then $\psi$ is locally $U_q(\Gg)_i$-finite.
\enlemma
\proof
Assume $t_i\psi=q_i^m\psi$.
It is enough to show
\beq e_i^{s}\psi=0.\label{aa}\endeq
Here $s=\max(n-m+1,0)$.
In order to see this, we may assume that $M$ is finite-dimensional.
Replacing $N$ with the $U_q(\Gg)_i$-module
generated by $\psi(M)$, we may assume that $N$ is
also finite-dimensional.
Hence $\Hom(M,N)$ is finite-dimensional and hence integrable.
In this case it is a well-known fact
that $f_i^{n+1}\psi=0$ implies (\ref{aa}).
\qed

\noindent
{\it Proof of Proposition \ref{serre}}.
\quad
Let $\ad:\tU\to \End(\tU))$
be a $\tU$-module structure on $\tU$
such that the multiplication $\tU\otimes M\to M$ is $\tU_q$-linear.
We have
\beq
\ad(t_i)(a)&=&t_iat_i^{-1}\\
\ad(e_i)(a)&=&e_ia-t_i^{-1}at_ie_i\\
\ad(f_i)(a)&=&[f_i,a]t_i^{-1}
\endeq
for $a\in\tU$.
Let $X:\tU\to \End(M)$ be the
homomorphism given by the $\tU$-module structure on $M$.
Let $i\not=j$. Since $[f_i,e_j]=0$,
$f_iX(e_j)=0$.
Since $X(e_j)$ has weight $\lan h_i,\alpha_j\ran$
with respect to
$U_q(\Gg)_i$, the preceding lemma implies
\beq e_i^{b}X(e_j)=0,\label{bb}\endeq
where $b=1-\lan h_i,\alpha_j\ran$.
On the other hand
\beqn
e_i^{b}X(e_j)&=&X(\ad(e_i^{b})e_j)\\
&=&X\left(\sum_{k=0}^{b}(-1)^{k}q_i^{-k(b-k)}\left[b\atop k\right]_i
e_i^{b-k}t_i^{-k}e_j(t_ie_i)^{k}\right)\\
&=&X\left(\sum_{k=0}^{b}(-1)^{k}\left[b\atop k\right]_i
e_i^{b-k}e_je_i^{k}\right).
\endeqn
This along with (\ref{bb}) gives the Serre relation
$$X\left(\sum_{k=0}^{b}(-1)^{k}\left[b\atop k\right]_i
e_i^{b-k}e_je_i^{k}\right)=0.$$

By applying the automorphism
$e_i\mapsto f_i$,
$f_i\mapsto e_i$
$q^h\mapsto q^{-h}$$(h\in P^*)$
of $\tU$,
we obtain the other
Serre relations
$$X\left(\sum_{k=0}^{b}(-1)^{k}\left[b\atop k\right]_i
f_i^{b-k}f_jf_i^{k}\right)=0.$$

\setcounter{equation}{0}
\setcounter{lemma}{0}
\newsection{Two-point function for $D^{(2)}_{n+1}$}
\label{app:c}

In this appendix we will sketch the calculation for level~1
$D^{(2)}_{n+1}$, of the two-point function $\Psi(z_1/z_2)=\langle
\Lambda_n | \Phi^{\Lambda_n V_2}_{\Lambda_n} (z_2) \Phi^{\Lambda_n
  V_1}_{\Lambda_n}(z_1)| \Lambda_n\rangle$, for the intertwiner
$\Phi^{\mu V}_{\lambda}(z) : V(\lambda) \longrightarrow V(\mu) \otimes
V_{z}$, by solving the corresponding $q$-KZ equation it must satisfy.
The corresponding calculations for the other cases in this paper have
been done in~\cite{IIJMNT} and~\cite{DO}. For the theoretical
background the appendix in~\cite{IIJMNT} should be consulted. To
conform with their conventions, we will use here the upper global base
and corresponding coproduct~$\Delta_+$, in contrast to the main text
of this paper.

Recall the total order $\succ$ on the index set $J$ defined in
(\ref{eqn:totord}). Extend the natural
definition of minus on $J\setminus\{\phi\}$ to all
of $J$ by defining $-\phi=\phi$. Let
\begin{eqnarray*}
  \overline{j}=j,
  \hspace{10pt} & & \hspace{10pt} \overline{-j}=2n+1-j, \hspace{20pt}
     j=1,\ldots,n \\
 \overline{0}=n, \hspace{10pt} & & \hspace{10pt} \overline{\phi}=2n.
\end{eqnarray*}
Denote, as usual, by $E_{jk}$ the matrix acting on $\{\v{j}\}_{j\in J}$
as $E_{jk} \v{i} =\delta_{ki}
\v{j}$. The R-matrix $\overline{R}(z)$ with normalization
 $\overline{R}(z) \v{1} \otimes \v{1} =  \v{1} \otimes \v{1}$ is then given by
\begin{eqnarray*}
\overline{R}(z) \!&=& \!\sum_{i\neq 0,\phi} E_{ii}\otimes E_{ii}+
\frac{q^2(1-z^2)}{(1-q^4 z^2)} \sum_{i\neq j, -j} E_{ii}\otimes E_{jj}
\nonumber\\
& & +\frac{(1-q^4)}{(1-q^4z^2)}\bigl(\sum_{i \succ j, i\neq -j}
  z^{\alpha_{ij}} E_{ij}\otimes E_{ji} + z^2 \sum_{i \succ j, i \neq -j}
  z^{-\alpha_{ij}} E_{ij}\otimes E_{ji} \bigr)
\nonumber\\
& &  + \sum a_{ij}(z) E_{ij} \otimes E_{-i,-j}
\end{eqnarray*}
where
\begin{eqnarray*}
  \alpha_{ij} &=&
  \begin{cases}
    1 &: \text{if $i=\phi$ or $j=\phi$},\\
    0 &: \text{otherwise},
  \end{cases}\\
  a_{ij}(z)&=&
  \begin{cases}
    (1-z^2)(q^4-\xi^2 z^2)+\delta_{i,-i}(1-q^2)(q^2+z^2) (1-\xi^2 z^2)
    &: i=j,
    \\
    (1-q^4)\bigl\{z^{\alpha_{ij}}(z^2-1)s_{ij}(-q^2)^{\overline{j}-
        \overline{i}}+ \delta_{i,-j}(1-\xi^2 z^2)\bigr\} &: i \succ j,
    \\
    (1-q^4)z^2\bigl\{\xi^2 z^{-\alpha_{ij}}(z^2-1)s_{ij}(-q^2)^
      {\overline{j}-\overline{i}}+ \delta_{i,-j} (1-\xi^2 z^2)\bigr\}
    &: i \prec j,
  \end{cases}
\end{eqnarray*}
and
\begin{alignat*}{4}
  s_{i0} &= -\frac{1}{[2]} \sgn(i)&\quad (i \neq \phi),
  &\qquad\quad &
  s_{0j} &= -[2]\sgn(j) &\quad& (j \neq \phi), \\
  s_{i\phi} &= \frac{1}{[2]} \sgn(i) &\quad (i \neq 0),
  &\qquad &
  s_{\phi j} &= [2]\sgn(j)& & (j \neq 0), \\
  s_{0\phi}&=s_{\phi 0}=-1, &
  &\qquad &
  s_{ij} &= \sgn(i)\;\sgn(j) &\quad &(i,j \neq 0,\phi).
\end{alignat*}
Also we have $\xi^2=q^{4n}$.  The expression for $\overline{R}(z)$ is
given in~\cite{J} in a different basis.

Let $\{\v{j}^*\}_{j\in J}$ be the canonical dual basis of the upper global
base.
The following isomorphism of $U_q({\frak g})$-modules
\begin{eqnarray*}
  C : V_{\xi^{-1}z} & \longrightarrow & (V_z)^{*\;a}\\
     \v{j} & \mapsto &\sgn(j)\; (-q^2)^{\overline{j}-1} \v{-j}^*
\quad (j \in J / \{0,\phi\})\\
     \v{0} & \mapsto & -\frac{1}{[2]}(-q^2)^{\overline{0}-1} \v{0}^*\\
     \v{\phi} & \mapsto & \frac{1}{[2]}\xi^{-1}
          (-q^2)^{\overline{\phi}-1} \v{\phi}^*
\end{eqnarray*}
gives rise to crossing-symmetry for the R-matrix
\begin{equation}
\bigl(R^{-1}(z)\bigr)^{t_1} = \beta(z) (C \otimes 1) R(z \xi^{-1})
    (C \otimes 1)^{-1},
\label{eqn:cross}
\end{equation}
with
\begin{equation}
\beta(z)=q^{-4} \frac{(1-z^2)(1-q^{-4n+4} z^2)}{(1-q^{-4n} z^2)(1-q^{-4}z^2)}.
\end{equation}
The image $R^+(z_1/z_2)=\pi_{V_{z_1}} \otimes \pi_{V_{z_2}}({\cal R}')$
of the modified universal R-matrix ${\cal R}'$ also satisfies
(\ref{eqn:cross}) with $z=z_1/z_2$ and $\beta(z)=1$. Therefore we have
\begin{equation}
R^+(z)=q^{-2}
\frac{(q^4z^2;\xi^4)_{\infty}(\xi^2z^2;\xi^4)_{\infty}^2
      (q^{-4}\xi^4z^2;\xi^4)_{\infty}}
     {(z^2;\xi^4)_{\infty}(q^{-4}\xi^2z^2;\xi^4)_{\infty}
      (q^4\xi^2z^2;\xi^4)_{\infty}(\xi^4z^2;\xi^4)_{\infty}}
\overline{R}(z).\nonumber
\end{equation}
The two-point function satisfies the $q$-KZ equation
\begin{equation}
\Psi ( q^{2(h^{\vee}+k)} z) = R^+(q^{2(h^{\vee}+k)} z)
(q^{-\phi} \otimes 1) \Psi(z),
\label{eqn:qkz1}
\end{equation}
where $k=1$ is the level and $\phi=2\Lambda^{\cl}_n+2\rho^{\cl}$
and, as a consequence, also
\begin{equation}
(\pi_{V_{z_1}} \otimes \pi_{V_{z_2}}) \Delta'(e_i)^{\langle h_i,\Lambda_n
\rangle+1} \Psi(z)=0, \hspace{20pt} \wt \Psi(z)=0.
\label{eqn:qkz2}
\end{equation}
It can be shown that
\begin{eqnarray}
w(z)& =& (1+z^2 q^2 \xi^2)\v{0} \otimes \v{0} + q [2] (-q^2)^n z
   \v{\phi} \otimes \v{\phi} -\nonumber\\
  & &q \sum_{i=1}^n (-q^2)^{n-i}
  (\v{i} \otimes \v{-i} + z^2 q^{4i-2} \v{-i}\otimes \v{i})
\end{eqnarray}
solves (\ref{eqn:qkz2}) and satisfies
\begin{equation}
\overline{R}(q^2\xi^2 z)
(q^{-\phi} \otimes 1) w(z) = q^2 \frac{(1-q^4\xi^2z^2)
(1-\xi^4z^2)}{(1-q^8\xi^4z^2)(1-q^4\xi^6z^2)}
w(q^2\xi^2 z).
\label{eqn:int}
\end{equation}
Letting $\Psi(z)=\phi(z)w(z)$ and substituting (\ref{eqn:int}) into
(\ref{eqn:qkz1}) one gets a scalar $q$-difference equation for $\phi(z)$
which can be solved to obtain the final result
\begin{equation}
\Psi(z) = \frac{(q^4 \xi^4 z^2;\xi^4)_{\infty}(\xi^6z^2;\xi^4)_{\infty}}
      {(q^4\xi^2z^2;\xi^4)_{\infty}(\xi^4z^2;\xi^4)_{\infty}}
  w(z).
\end{equation}


\setcounter{equation}{0}
\setcounter{lemma}{0}
\newsection{The limit $q\rightarrow 1$ for the $U_q(A^{(2)}_{2n})$
Fock space}
\label{app:d}

In this appendix we will show how to recover the known classical (i.e.\ at
$q=1$) Fock space ${\cal F}_{\text{class}}$ for ${\frak{g}}=A^{(2)}_{2n}$
at level 1. This involves reduction
of the Fock space ${\cal F}$ defined for generic $q$ by means of an invariant
inner product on ${\cal F}$. To facilitate the discussion we shall make a
transcription from the semi-infinite wedge description of ${\cal F}$ to
one involving Young diagrams or, synonymously, partitions
(the so-called ``combinatorial description'').

Define the following subspace of $V_{\aff}$:
\begin{equation}
V_{\aff}^+ = z^{-1} {\Bbb Q}[z^{-1}] \otimes V  \oplus
{\Bbb Q} \langle v_{-1},\ldots, v_{-n}, v_0 \rangle.
\end{equation}
In any normally ordered pure wedge in ${\cal F}$ it is clear that only
bases in $V_{\aff}^+$ appear as components. Recall the single-valued
function $l$ on $\Baff$ in~\eqref{eq:a2even-l-func}.  To the normally
ordered pure wedge $u=G(u_1) \wedge G(u_2) \wedge \cdots \wedge G(u_k)
\wedge v_0 \wedge v_0 \wedge \cdots$ let us associate the sequence
$Y(u)= [ -l(u_1), -l(u_2), \ldots, -l(u_k), 0, 0, \ldots ]$, whose
tail of zeroes we shall ignore.  Now, $-l$ takes non-negative values
on $V_{\aff}^+$. Also, the sequence $Y(u)$ is non-increasing because
of the normal-ordering rules. Furthermore, the only integers allowed
to repeat belong to $h{\Bbb N}$, where $h=2n+1$, because of the rule
$v_i \wedge v_i =0$ if $i \neq 0$.  Thereby we have the identification
\begin{equation}
{\cal F} \simeq {\Bbb Q}(q)\langle Y \rangle_{Y \in \DP_{h}},
\end{equation}
where $\DP_k$ is the set of Young diagrams whose rows are allowed to repeat
only if their length is $0 \mod k$. In this notation, $\DP_{\infty}$ is the
set of Young diagrams with {\em no} repeating rows, i.e., the set of
Distinct Partitions.

The action of $\Uq$ on ${\cal F}$ can be transcribed to the Young
diagram setting. The generators $t_i$ act diagonally, of course, while
$f_i$ (respectively $e_i$) act by adding (removing) one box in the
following manner.  Let the Young diagram $Y$ be denoted by
$[y_1,\ldots,y_m]$.  For $y\in\Zplus$, let $\alpha_Y(y)$ denote the
number of occurences of $y$ in $Y$. Define the functions $\beta_i$ for
$i=0,1,\ldots,n$ by
\begin{eqnarray*}
\beta_0(y) & = & \begin{cases}
       \pm 4 &: y \in h{\Bbb Z} \pm n\\
       0 &: \text{otherwise} \end{cases} \nonumber\\
\beta_i(y) & = & \begin{cases}
       \pm 2 &: y \in h{\Bbb Z} \mp n \pm (i-1)\\
       \mp 2 &: y \in h{\Bbb Z} \mp n \pm i\\
       0 &: \text{otherwise} \end{cases}
\hspace{20pt} (i=1,\ldots,n-1)\nonumber\\
\beta_n(y) & = & \begin{cases}
       \pm 2 &: y \in h{\Bbb Z} \mp 1\\
       0 &: \text{otherwise} \end{cases}
\end{eqnarray*}
Then the action of $\Uq$ on $Y$ is given explicitly by
\begin{align*}
t_i \cdot Y =&\; q^{\sum_j \beta_i(y_j) + \delta_{i,n}} \; Y \\
f_i \cdot Y =& \sum
  \begin{Sb}
    y_k \in h{\Bbb N} + n \pm i\\
    y_{k-1}\neq y_k+1
  \end{Sb}
q^{\sum_{j>k} \beta_i(y_j)}
\; [y_1,\ldots,y_k+1,\ldots,y_m] \qquad\qquad (i \neq n)\\
f_n \cdot Y =& \sum
  \begin{Sb}
    y_k \in h{\Bbb N} -1\\
    y_{k-1}\neq y_k+1
  \end{Sb}
      q^{\sum_{j>k} \beta_i(y_j)+1}
            \; [y_1,\ldots,y_k+1,\ldots,y_m] \nonumber\\
            & +\sum
              \begin{Sb}
                y_k \in h{\Bbb N}\\
                y_{k-1}\neq y_k+1,y_k
              \end{Sb}
      q^{\sum_{j>k} \beta_i(y_j)+1}
            \left(1-(-q^2)^{\alpha_Y(y_k)}\right) \;
            [y_1,\ldots,y_k+1,\ldots,y_m] \nonumber\\
            & +\delta(y_m \neq 1) [y_1,\ldots,y_m,1]\\
e_i \cdot Y =& \sum
  \begin{Sb}
    y_k \in h{\Bbb N} + n +1 \pm i\\
    y_{k+1}\neq y_k-1
  \end{Sb}
  q^{-\sum_{j<k}
    \beta_i(y_j)} \; [y_1,\ldots,y_k-1,\ldots,y_m] \qquad\qquad (i \neq n)\\
e_n \cdot Y =& \sum
\begin{Sb}
  y_k \in h{\Bbb N} +1\\
  y_{k+1}\neq y_k-1
\end{Sb}
q^{-\sum_{j<k}
    \beta_i(y_j)} \; [y_1,\ldots,y_k-1,\ldots,y_m] \nonumber\\
            & +\sum
            \begin{Sb}
              y_k \in h{\Bbb N}\\
              y_{k+1}\neq y_k-1,y_k
            \end{Sb}
 q^{-\sum_{j<k}
    \beta_i(y_j)-1} \; \left(1-(-q^2)^{\alpha_Y(y_k)} \right)\;
    [y_1,\ldots,y_k-1,\ldots,y_m].
\end{align*}
Note that all Young diagrams appearing on the right-hand side belong
to $\DP_h$. In other words, the corresponding pure wedges are already
normally ordered. The factors $\left(1-(-q^2)^{\alpha_Y(y_k)} \right)$
come from normal ordering and summing up Young diagrams which arise
when $Y$ has repeated rows. Note also that the vacuum vector is the
empty Young diagram $\emptyset$ and $f_n \cdot \emptyset = [1]$. This
combinatorial description is in the same spirit as that for
$U_q(A^{(1)}_n)$ given in~\cite{MM}.

Let us now introduce an inner product $(\;,\;)$ on ${\cal F}$. We shall
require that the normally ordered pure wedges, or equivalently Young
diagrams in $\DP_h$, form an orthogonal basis with respect to
$(\;,\;)$. We shall also require that
with respect to $(\;,\;)$ the adjoints of the generators satisfy
\begin{eqnarray*}
  f_i^{\dagger} &=& q_i \;e_i \;t_i, \\
  e_i^{\dagger} &=& q_i \;f_i \;t_i^{-1}, \\
  t_i^{\dagger} &=& t_i.
\end{eqnarray*}
These conditions are natural for a $\Uq$-module $V$ because then on
the module $V \otimes V$ with induced inner product given by $(v_1
\otimes v_2, u_1 \otimes u_2)=(v_1,u_1)(v_2,u_2)$ we have
$\Delta(f_i)^{\dagger} = q_i \Delta(e_i) \Delta(t_i)$, etc.  Using the
explicit description of the $\Uq$ action on ${\cal F}$ one can show
that the squared norm of an arbitrary Young diagram $Y$ in $\DP_h$ is
given by
\begin{equation}
  || Y ||^2 = (Y,Y) = \prod_{y\in h\Zplus}
  \prod_{i=1}^{\alpha_Y(y)}
  \bigl( 1- (-q^2)^i \bigr).
  \label{eqn:inn}
\end{equation}
{}From calculations for small $k$ we conjecture that the boson operators
satisfy
\begin{displaymath}
    B_{-k}^{\dagger} = B_k.
\end{displaymath}

As at the end of \S\ref{subsec:classical-limit}, we denote by ${\cal
  F}^{\Bbb Q}$ the ${\Bbb Q}[q,q^{-1}]$-vector space generated by the
pure wedges.  Set ${\cal F}_1 = {\cal F}^{\Bbb Q}/(q-1){\cal F}^{\Bbb
  Q}$. Then the action of $\Uq$ on ${\cal F}$ induces an action of
$\Uq$ on ${\cal F}_1$. The inner product $(\;,\;)$ on ${\cal F}^{\Bbb
  Q}$ induces a ${\Bbb Q}$-valued inner product on ${\cal F}_1$, which
we also denote by $(\;,\;)$. The adjoint of operators in $\frak{g}$ is
then given by $e_i^{\dagger}=f_i$, $f_i^{\dagger}=e_i$ and
$h_i^{\dagger}=h_i$.  Define the subspace ${\cal F}_0 = \{ u \in {\cal
  F}_1: (u,{\cal F}_1) =0 \}$.  The reduced Fock space ${\cal
  F}_{\text{red}}= {\cal F}_1 / {\cal F}_0$ is a $U(\frak{g})$-module.
{}From (\ref{eqn:inn}) we note that ${\cal F}_0$ is the ${\Bbb Q}$-span
of Young diagrams with {\em some} repeated rows. It follows that
${\cal F}_{\text{red}}$ is the ${\Bbb Q}$-span of Young diagrams in
$\DP_{\infty}$. This is isomorphic to the well-known classical Fock
space ${\cal F}_{\text{class}} \simeq {\Bbb Q}[x_k]_{k \in {\Bbb
    N}_{\text{odd}}}$~\cite{KKLW,DJKM}.  In fact, the action of the
generators on ${\cal F}_{\text{red}}$ and at $q=1$ reduces to a known
classical action~\cite{JY}. Furthermore we recover the known
decomposition of ${\cal F}_{\text{class}}\simeq {\Bbb Q}[x_k]_{k \in
  {\Bbb N}_{\text{odd}}\setminus h{\Bbb N}} \otimes {\Bbb
  Q}[x_{hk}]_{k \in {\Bbb N}_{\text{odd}}}$ as a $U({\frak{g}})
\otimes {\Bbb Q}[H_-]$-module. Here we identify bosons $x_{hk}\sim
B_{-k}$ for $k \in {\Bbb N}_{\text{odd}}$. The even boson commutators
$\gamma_{k}$ for $k \in {\Bbb N}_{\text{even}}$ have a pole at $q=1$.
After appropriately rescaling we find that such $B_{-k}$ act as 0 on
${\cal F}_{\text{red}}$ at $q=1$.

In most of the cases considered in this paper, the boson commutator $\gamma_k$
has a pole at $q=1$ for some $k$. We take this to indicate that similar
Fock space reductions to the one considered in this Appendix are necessary
to recover any known classical Fock spaces.


\newpage
\ifx\undefined\bysame
\newcommand{\bysame}{\leavevmode\hbox to3em{\hrulefill}\,}
\fi

\bigskip
\end{document}